\documentclass[reqno]{amsart}

\usepackage{latexsym}
\usepackage{graphicx}

\newcommand{\T}{{\mathbf T}^m}

\newcommand{\sm}{\setminus}
\newcommand{\szego}{Szeg\"o }

\newcommand{\Si}{\Sigma}
\newcommand{\inv}{^{-1}}
\newcommand{\kahler}{K\"ahler }

\newcommand{\wt}{\widetilde}
\newcommand{\wh}{\widehat}
\newcommand{\PP}{{\mathbb P}}
\newcommand{\N}{{\mathbb N}}
\newcommand{\R}{{\mathbb R}}
\newcommand{\C}{{\mathbb C}}

\newcommand{\Z}{{\mathbb Z}}
\newcommand{\B}{{\mathbb B}}
\newcommand{\CP}{\C\PP}
\renewcommand{\d}{\partial}
\newcommand{\dbar}{\bar\partial}
\newcommand{\ddbar}{\partial\dbar}

\newcommand{\E}{{\mathbf E}}

\newcommand{\half}{{\frac{1}{2}}}
\newcommand{\Log}{{\operatorname{Log\,}}}
\newcommand{\vol}{{\operatorname{Vol}}}

\newcommand{\codim}{{\operatorname{codim\,}}}
\newcommand{\SU}{{\operatorname{SU}}}
\newcommand{\FS}{{{}_{\rm FS}}}

\renewcommand{\phi}{\varphi}

\newcommand{\acal}{\mathcal{A}}

\newcommand{\ccal}{\mathcal{C}}
\newcommand{\dcal}{\mathcal{D}}
\newcommand{\ecal}{\mathcal{E}}
\newcommand{\fcal}{\mathcal{F}}

\newcommand{\hcal}{\mathcal{H}}

\newcommand{\lcal}{\mathcal{L}}

\newcommand{\ncal}{\mathcal{N}}
\newcommand{\ocal}{\mathcal{O}}
\newcommand{\pcal}{\mathcal{P}}
\newcommand{\qcal}{\mathcal{Q}}
\newcommand{\rcal}{\mathcal{R}}
\newcommand{\scal}{\mathcal{S}}
\newcommand{\tcal}{\mathcal{T}}

\newcommand{\al}{\alpha}
\newcommand{\be}{\beta}
\newcommand{\ga}{\gamma}
\newcommand{\Ga}{\Gamma}
\newcommand{\La}{\Lambda}
\newcommand{\la}{\lambda}
\newcommand{\ep}{\varepsilon}
\newcommand{\de}{p} 
\newcommand{\De}{\Delta}
\newcommand{\om}{\omega}

\newcommand{\di}{\displaystyle}
\newtheorem{theo}{{\sc Theorem}}[section]
\newtheorem{maintheo}{{\sc Theorem}}
\newtheorem{cor}[theo]{{\sc Corollary}}
\newtheorem{maincor}[maintheo]{{\sc Corollary}}
\newtheorem{lem}[theo]{{\sc Lemma}}
\newtheorem{prop}[theo]{{\sc Proposition}}

\newenvironment{rem}{\medskip\noindent{\it Remark:\/}}{\medskip}
\newenvironment{defin}{\medskip\noindent\it Definition: }{\medskip}

\title[Random polynomials with prescribed Newton polytope] {Random
polynomials with prescribed Newton\\ polytope, I }

\author{Bernard Shiffman}
\author{Steve Zelditch}
\address{Department of Mathematics, Johns Hopkins University, Baltimore, MD
21218, USA}
\email{shiffman@math.jhu.edu, zelditch@math.jhu.edu}

\thanks{Research partially supported by NSF grant
DMS-0100474 (first author) and DMS-0071358 (second author).}

\date{March 12, 2002}

\begin{document}

\begin{abstract} The Newton polytope $P_f$ of a polynomial $f$ is
well known to have a strong impact on its zeros, as in the
Kouchnirenko-Bernstein theorem on the number of simultaneous zeros of $m$
polynomials with  given Newton polytopes.  In this article, we
show that
$P_f$ also has a strong impact on the distribution of zeros of one
or several polynomials. We equip the
space of (holomorphic) polynomials of degree
$\leq N$ in $m$ complex variables with its
usual $\SU(m + 1)$-invariant Gaussian measure and then
consider the conditional measures $\gamma_{| NP}$ induced on the
subspace of  polynomials whose Newton
polytope $P_f \subset NP$. When $P =  \Sigma$, the unit simplex,
then it is obvious and well-known that the expected distribution
of zeros $Z_{f_1, \dots, f_k}$ (regarded as a current) of $k$ polynomials
$f_1,\dots,f_k$ of degree $N$ is uniform relative to the Fubini-Study form. Our
main results concern the conditional  expectation 
$\E_{|N P}(Z_{f_1, \dots, f_k})$ of zeros of $k$ polynomials with Newton polytope
$NP\subset Np\Sigma$ (where $p=\deg P$); these results are
asymptotic as the scaling factor $N\to\infty$.  We show that $\E_{|N P}(Z_{f_1,
\dots, f_k})$ is asymptotically uniform on the inverse image
${\mathcal A}_P$  of the open scaled polytope $\frac{1}{p}P^\circ$ via the moment
map
$\mu:\CP^m\to \Sigma$ for projective space. However, the zeros have an exotic
distribution outside of
${\mathcal A}_P$ and when $k = m$ (the case of the
Kouchnirenko-Bernstein theorem) the percentage of zeros outside  ${\mathcal A}_P$
approaches 0 as $N\to\infty$. 
\end{abstract}

\maketitle

\tableofcontents

\section*{Introduction}

It is well known that the Newton polytope $P_f$ of a holomorphic
polynomial $f(z_1, \dots, z_m)$ of degree $p$ has a crucial influence on
its value distribution and in particular on its zero set.  Even the number of
simultaneous zeros   in $\C^{*m}:=(\C\sm\{0\})^m$
of $m$ generic polynomials $f_1, \dots, f_m$
depends on their Newton polytopes $P_{f_j}$ \cite{Be, Ku1, Ku2}. Our
purpose in this paper  is to demonstrate that the Newton polytope
of a polynomial $f$ also has a crucial influence on
its {\it mass density} $|f(z)|^2 dV$ and on the {\it
spatial distribution} of
{\it zeros} $\{f = 0\}$.
We will show that there is a {\it classically allowed
region} region where the mass almost surely concentrates  and a
{\it classically forbidden region} where it almost surely is
exponentially decaying. The classically allowed region is the inverse
image $\mu_{\Sigma}^{-1}(\frac{1}{\de}P)$ of the (scaled) polytope
$\frac{1}{\de}P$ under the standard moment map $\mu_\Si$ of $\CP^m$. The
simultaneous zeros of $m$ generic polynomials $f_1, \dots, f_m$ in $\C^{*m}$
almost surely concentrates (in the limit of high degrees) in the classically
allowed region, giving a kind of quantitative version of the
Bernstein-Kouchnirenko theorem. The Newton polytope has an equally
strong (though different) impact for $k < m$ polynomials $f_1,
\dots, f_k$.  The  image of the
zero set of $f_1,
\dots, f_k$ under the moment map is (up to a
logarithmic re-parametrization) known as an {\it amoeba\/} in the
sense of \cite{GKZ, M}. Results on the expected distribution of amoebas can be obtained
from our results on the expected zero current; for example, as a consequence of
Corollary~\ref{subtler}, for polytopes
$\frac{1}{\de}P$ with vertices in the interior of the standard simplex $\Sigma \subset
\R^m$, there is also a {\it very forbidden region} which the
amoeba almost surely avoids. These patterns in zeros are statistical---they hold for
random polynomials with prescribed Newton polytope---and are asymptotic as the degree of
the polynomials tends to infinity.

      To state our problems and results precisely, let us recall
      some definitions. Let $$f(z_1, \dots, z_m) = \sum_{\alpha \in \N^m:
|\alpha| \leq p} c_{\alpha} \chi_\al(z_1, \dots, z_m),\;\;\;\;\
\chi_\al(z) = z_1^{\alpha_1} \cdots z_m^{\alpha_m}\;$$ be a
polynomial of degree $p$ in
      $m$ complex variables. By the {\it support} of $f$ we mean the set
    \begin{equation} \label{SUPPORT} S_f = \{\alpha\in\Z^m: c_{\alpha} \not
= 0\},\end{equation}
    and by its   {\it Newton polytope} $P_f$ we mean the integral polytope
\begin{equation}\label{NEWTONP}  P_f : = \hat{S}_f =
\mbox{the convex hull in $\R^m$ of } \;\; S_f. \end{equation}
Throughout this paper we will make the simplifying assumption that
$P$ is a Delzant polytope \cite{De, G}. Non-Delzant polytopes are
equally important in the study of polynomials, and we plan to
 extend our results to more general convex integral
polytopes in  the sequel \cite{SZ2}.  We also plan to study a parallel
tunneling phenomenon for critical points.

Our aim is to study the statistical patterns in polynomials with a fixed
Newton polytope.  We therefore define  a natural  Gaussian measure
$d\gamma_P$ on the space of polynomials with support contained in
a fixed Newton polytope $P$. Since our purpose is to compare zero
sets and masses as the polytope $P$ varies, we  view the polytope
as placing a condition on the Gaussian ensemble of all polynomials
of degree $p$ equal to that of $P$, i.e. we define  the
measures $d \gamma_P$ as conditional probability measures of one
fixed Gaussian measure.
To define this ensemble, we homogenize the
polynomial $f$ to obtain a homogeneous polynomial  $F$ of degree
$p$ in $m + 1$ complex variables. Such polynomials may
be identified  with holomorphic sections $F\in H^0(\CP^m,
{\mathcal O}(p))$ of the $p^{\rm th}$ power of the hyperplane
section bundle ${\mathcal O}_{\CP^m}(1)$.
  We  denote by
\begin{equation} \label{SUBSPACE} H^0(\CP^m,
{\mathcal O}(p), P) = \{F \in  H^0(\CP^m, {\mathcal O}(p)): P_f \subset P\}
\end{equation}
the space of homogeneous polynomials $F$ of
degree $d$ whose associated inhomogeneous form $f(z_1, \dots, z_m)
= F(1, z_1, \dots, z_m)$  has Newton polytope $P_f$ contained in $P$.

Let us recall that  $H^0(\CP^m, {\mathcal
O}(p))$ carries  the standard $\SU(m+1)$-invariant  inner product
$$\langle F_1, \bar F_2\rangle = \int_{S^{2m+1}} F_1 \bar F_2\, d\sigma\;,$$
where $d\sigma$ is Haar measure on the $(2m+1)$-sphere $S^{2m+1}$. (See
(\ref{IP}) for an alternate description in terms of the Fubini-Study form on
$\CP^m$  and the Fubini-Study metric  on ${\mathcal O}(1)$.) An orthonormal basis of
$H^0(\CP^m, {\mathcal
O}(p))$ is given by $\left\{\|\chi_{\alpha}\|\inv \chi_{\alpha}\right\}_{|\al|\le
p}$, where $\|\cdot\|$ denotes the norm in $H^0(\CP^m,\ocal(\de
))$. (Note that $\|\chi_{\alpha}\|$ depends on $\de$.) The corresponding
$\SU(m+1)$-invariant  Gaussian measure $\gamma_\de$ is defined by
\begin{equation}
\label{G} d
\gamma_\de (s) = \frac{1}{\pi^{k_p}}e^{-|\la|^2} d\la,\;\;\;\; s =
\sum_{|\alpha|\le p} \la_{\alpha}
\frac{\chi_{\alpha}}{\|\chi_{\alpha}\|}\;,\end{equation}
where $k_p=\#\{\al: |\al|\le p\}= 
{m+p\choose p}$.  
Thus, the coefficients $\la_{\alpha}$ are
independent complex Gaussian random variables with mean zero and variance
one. We then endow the space $H^0(\CP^m, {\mathcal O}(p), P)$
with the associated  {\it conditional probability measure\/} $\gamma_\de|_P$:
 \begin{equation} \label{CG} d \gamma_{\de|P} (s) =
\frac{1}{\pi^{\# P}}e^{-|\la|^2} d\la,\quad s = \sum_{\alpha
\in P} \la_{\alpha} \frac{\chi_{\alpha}}{\|\chi_{\alpha}\|}\;,  \end{equation}
where the coefficients $\la_{\alpha}$ are again independent complex Gaussian random
variables with mean zero and variance one. (By a slight abuse of notation, we
let $\sum_{\al\in P}$ denote the sum over the lattice points $\al\in
P\cap\Z^m$; $\# P$ denotes the cardinality of $P\cap\Z^m$.) We observe
that,
 as a subspace of $H^0(\CP^m, {\mathcal O}(p))$,  $H^0(\CP^m, {\mathcal
O}(p), P)$ inherits the inner product $\langle s_1,s_2
\rangle$ and that $\gamma_{|P}$ is the induced  Gaussian measure.
Probabilities relative to $\gamma|_P$  can be  considered as conditional
probabilities; i.e. for any event $E$,
$$\mbox{Prob}_{\gamma}\{ f \in E | P_f = P\} = \mbox{Prob}_{\gamma|_P} (E).$$

We now wish to consider the expected distribution of mass and
zeros of polynomials with fixed Newton polytope. It turns out to
involve the moment map $\mu_{\Sigma}$ of projective space.  We use
a notation for this moment map which reflects its dependence on
the
 standard unit simplex $\Si$ in $\R^m$ with vertices at
$(0,\dots,0),\ (1,0,\dots, 0),\ (0,1,\dots,0),\dots,(0,\dots,0,1)$. The `projective
moment map'
is then given by
$$\mu_\Si(z)=\left(\frac{|z_1|^2}{1+\|z\|^2},\dots,
\frac{|z_m|^2}{1+\|z\|^2}
\right)\;.$$
Via this moment map we define:

\begin{defin} Let $P \subset \R_+^m$ be an integral polytope. The
 {\it classically allowed region\/} for polynomials in $H^0(\CP^m, {\mathcal
O}(p), P)$
is the set
$$\acal_P:=\mu_\Si\inv\left(\frac{1}{p}P^\circ\right) \subset \C^{*m}$$
(where
$P^\circ$ denotes the interior of $P$), and the {\it classically forbidden
region\/} is its complement $\C^{*m}\sm\acal_P$. \end{defin}

Our first result concerns the simultaneous  zero set of $m$
independent polynomials in $m$ variables. B\'ezout's theorem tells
us that $m$ generic homogeneous polynomials $F_1,\dots,F_m$ of
degree $p$ have exactly $p^m$ simultaneous zeros; these
zeros all lie in $\C^{*m}$, for generic $F_j$. In fact, one immediately
sees (by uniqueness of Haar probability measure) that the expected
distribution of zeros  is uniform with respect to the Fubini-Study projective
volume form, when the ensemble is given the $\SU(m+1)$-invariant measure
$d\ga_p$.

According to Kouchnirenko's theorem
\cite{Ku1,Ku2,A,Au}, the number of joint zeros in $\C^{*m}$ of $m$ generic
polynomials
$\{f_1,\dots, f_m\}$ with given Newton polytope $P$ equals $ m!
\vol(P)$, where $\vol(P)$ is the Euclidean volume of $P$. For
example, if $P=p\Si$, where $\Si$ is the standard unit
simplex in $\R^m$, then $\vol
(p\Si)=p^m\vol(\Si)=\frac{p^m}{m!}$, and we get
B\'ezout's theorem.  (More generally, the $f_j$ may have different
Newton polytopes, in which case the number of zeros is given by
the Bernstein-Kouchnirenko formula as a `mixed volume.')

Now consider  $m$ independent random
polynomials with Newton polytope
$P$, using the conditional probability (\ref{CG}). We let
$\E_{|P} (Z_{f_1,
\dots, f_m})$ denote the expected density of their simultaneous zeros.
It is the measure on $\C^{*m}$ given
by
\begin{equation}\label{EZ0}\E_{|P} (Z_{f_1, \dots, f_m})(U) =
\int d\ga_{p|P}(f_1)\cdots\int d\ga_{p|P}(f_m)\;\big[\#\{z\in U:
f_1(z)=\cdots=f_m(z)=0\}\big]\;,\end{equation}
for  $U\subset
\C^{*m}$, where the integrals are over
$H^0(\CP^m,\ocal(p),P)$.  In fact, $\E_{|P} (Z_{f_1, \dots, f_m})$ is
an absolutely continuous measure given by a $\ccal^\infty$ density
(see Prop.~\ref{EZsimult}).  Our first result is the  surprising fact that as
the polytope $P$ expands, these zeros
 are expected to concentrate in the classically
allowed region and are (asymptotically) uniform there:

\begin{maintheo}\label{probK} Suppose that $P$ is a Delzant polytope in $\R^m$.
Then
$$\frac{1}{(N\de)^m}\E_{|NP} (Z_{f_1, \dots, f_m})\to
\left\{\begin{array}{ll}\om_\FS^m 
\ \ & \mbox{\rm on \ }\acal_P\\[10pt]
0 & \mbox{\rm on \ }\C^{*m}\sm \acal_P\end{array}
\right.\ ,$$ in the distribution sense; i.e., for any open
$U\subset\C^{*m}$, we have
$$\frac{1}{(N\de)^m}\E_{|NP}\big(\#\{z\in U:
f_1(z)=\cdots=f_m(z)=0\}\big)\to m!\vol_{\CP^m}(U\cap \acal_P)\;.$$
\end{maintheo}

In fact, our methods imply that the convergence of the zero current on the classically
allowed region is exponentially fast in the sense that
$$\E_{|NP} (Z_{f_1, \dots, f_m})= (Np)^m \om_\FS^m  +O\left(e^{-\la N}\right)
\ \  \mbox{\rm on \ }\acal_P\;,$$ for some positive continuous function $\la$ on
$\acal_P$.

Theorem \ref{probK} is a special case of a
general result stated below (Theorem~\ref{simultaneous}) on
the expected distribution of simultaneous zeros of $k$ random
polynomials with given Newton polytopes $P_1,\dots,P_k$.  Although
we only discuss expected behavior of zeros here, it is possible to
show (using methods similar to those of  \cite{SZ})
that distribution of zeros is self-averaging: i.e., almost all
polynomials exhibit the expected behavior in an asymptotic sense.
For the sake of brevity,  we concentrate here on only the
essentially new aspects of our problem, namely the dependence of
the expected behavior on $P$.

The following illustration\footnote{We would like to thank Aaron Carass for
providing the illustrations in  Figures \ref{allowed-square}, \ref{f-fan}, \ref{L},
\ref{F2} and
\ref{Fn}.} shows the classically allowed region (shaded) and the classically 
forbidden region (unshaded) when $P$ is the unit square in
$\R^2$ (and
$\de=2$):

\begin{figure}[htb]
\centerline{\includegraphics*[bb= 1.9in 7.0in 4in 9.6in]{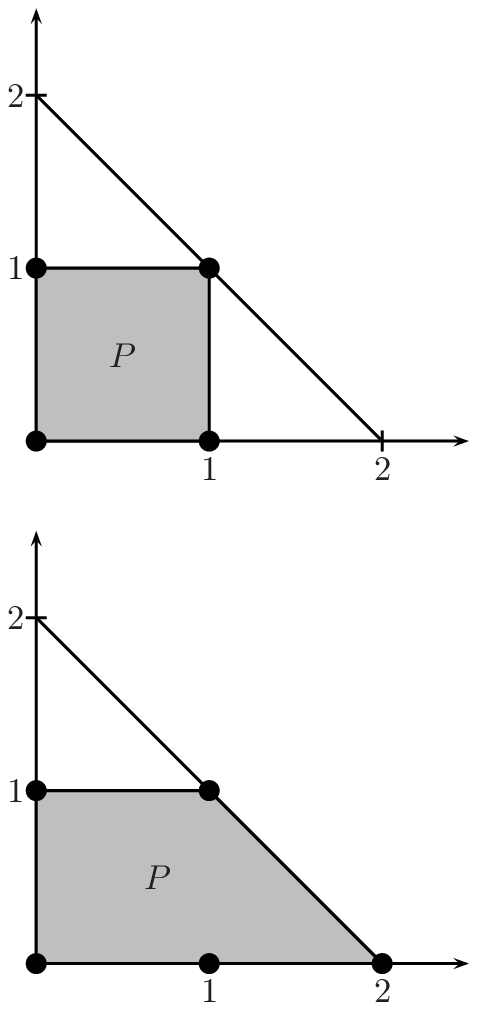}
\includegraphics*[bb= 1.6in 6.7in 5in 9.6in]{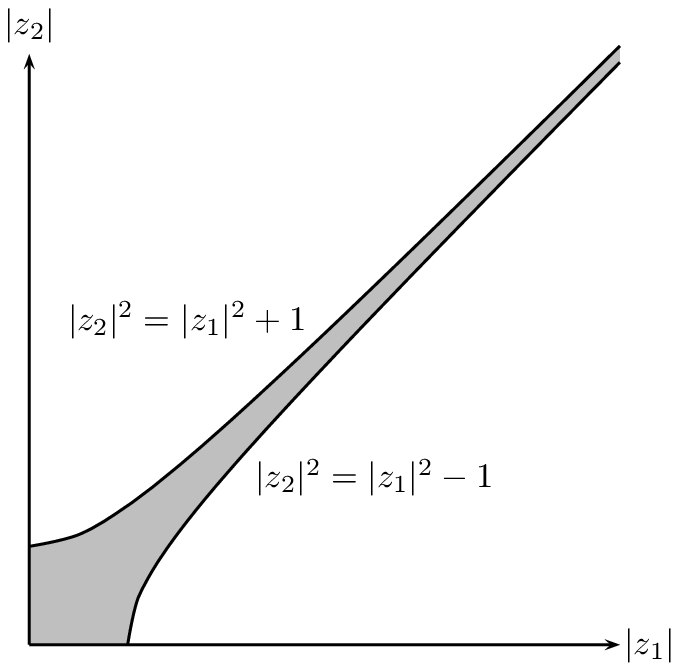}}
\caption{The classically allowed region for $P=[0,1]\times [0,1]$}
\label{allowed-square}\end{figure}

The terminology is taken from the semiclassical analysis of ground states of
Schroedinger
operators $H_h = -\hbar^2 \Delta + V$ on $\R^n$. The well-known `Agmon
estimates' of ground states
(cf.\ \cite{Ag}) show that $\lcal^2$-normalized ground states or low-lying
eigenfunctions
$H_\hbar \phi = E \phi$ of $H$ are concentrated as $\hbar \to 0$  in the classically
allowed region
$C_E = \{x \in \R^n: V(x) \leq E\}$ and have exponential decay $|\phi(x)|^2
= O(e^{- d(x, C_E)/\hbar})$ as $\hbar \to 0$ in the complement. Here, $d(x, C_E)$
is the `Agmon distance'
to the classically allowed region, i.e. the distance in the metric
$\sqrt{V(x) - E}\, dx$.
In our setting, the Hamiltonian is $\bar{\partial}^* \bar{\partial}$ on
$\lcal^2$-sections of powers ${\mathcal O}(N p)$ of the hyperplane section bundle,
$\hbar = 1/N$,
and the ground states are the holomorphic sections $H^0(\CP^m, {\mathcal
O}(p), P)$.
Replacing $|\phi(x)|^2$ is the expected  `mass density' of the polynomials
in this subspace.  By
definition, the mass density of a polynomial $f$ at a point $z\in \C^{*m}$
is the Fubini-Study norm
\begin{equation}\label{FSnorm}|f(z)|_\FS := \frac {|f(z)|}{(1+\|z\|^2)^{\delta/2}}
\qquad (\delta=\deg f)\;.\end{equation}

Our next result describes the asymptotics of the expected mass density of
$\lcal^2$-normalized polynomials with Newton polytope $NP$, i.e.\ the expected
density $\E_{\nu_{NP}}(|f(z)|^2_\FS)$ with respect to the Haar probability
measure, denoted by $\nu_{NP}$,
on the unit sphere in $H^0(\CP^m,\ocal(Np),NP)$.  One easily shows (see \S
\ref{s-mass})  that
\begin{equation}\label{2measures}\E_{\nu_{NP}}(|f(z)|^2_\FS) =
\frac{1}{\# (NP)}\E_{|NP}
\left(|f(z)|^2_\FS\right) =\frac{1}{\# (NP)}\Pi_{|NP}(z,z)\;,\end{equation}
where $\Pi_{|NP}$ is the \szego kernel of $H^0(\CP^m, {\mathcal O}(Np), NP)$. Here, the
term `\szego kernel' of a space $\scal$ of holomorphic  sections of a line bundle refers
to the kernel for the orthogonal projection to
$\scal$ from the space
$\lcal^2(\d D)$ of square-integrable functions on  the boundary of the  associated unit
disk bundle
$D$; i.e., it is of the form
$\Pi(x,y)=\sum
_j s_j(x) \overline{s_j(y)}$, where
$\{s_j\}$ is an orthonormal basis of $\scal$. Precise definitions and the  relationship
to the classical
\szego kernel for a strictly pseudoconvex domain are described in \S \ref{s-szego}.
Note that $$\# (NP) =\dim H^0(\CP^m,\ocal(Np),NP) =
 \vol(P) N^m +\dots\,,$$
which is a polynomial given by the Riemann-Roch formula (see \S \ref{s-mass}).

We shall show  that the expected mass density $\E_{\nu_{NP}}(|f(z)|^2_\FS)$ is
asymptotically uniform with respect to Fubini-Study measure in the classically
allowed region (as if there were no constraint at all); while in the forbidden
region the mass decays exponentially. Thus in the semiclassical limit
$N\to\infty$, all the mass concentrates in the classically allowed
region. To describe the behavior in the forbidden region, we
subdivide it into regions corresponding to the faces of the
polytope $\frac{1}{\de}P$ that lie in the interior $\Si^\circ$ of
the simplex.  To describe these regions, we recall that the fan
associated to $\frac{1}{\de}P$ (or $P$) is the collection of
polyhedral cones $C_F\subset \R^m$ normal to the faces $F$ of
$\frac{1}{\de}P$ (see \S \ref{s-fans}). For each face
$F\subset\Si^\circ$, we let
\begin{equation} \rcal_F = \{e^\tau \cdot z: \tau\in C_F ,\ \mu_\Si(z) \in
F\}\;.
\label{RF}\end{equation}
Note that if $F$ is the open face $\frac{1}{\de}P^\circ$, then
$\rcal_F$ is the classically allowed region $\acal_P$.  We show in
\S
\ref{s-normalbundle} that the classically forbidden
region decomposes into the (disjoint) forbidden subregions
$$\textstyle\left\{\rcal_F: F
\mbox{ is a  face of}\
\frac{1}{\de}P,\ F\subset \Si^\circ,\ F\neq \frac{1}{p}P^\circ\right\}\;.$$
We call $\rcal_F$ the {\it flow-out\/} of the face $F$. Boundary points of the
$\rcal_F$ are called {\it transition points.} Points of $\d \acal_P$ are always
transition points, but there may be others. 

For the case where $P$ is the unit
square, then there are two interior faces whose flow-outs
are the connected components of the classically forbidden region
(see Figure~\ref{allowed-square} above). In this case, the set of transition points
equals $\d \acal_P$.  We shall give an example with an interior vertex,
where the forbidden region is connected but decomposes into flow-outs of 3 faces
(see Figure~\ref{Fn} in \S \ref{s-ex3}); in this case there are transition points
not in $\d
\acal_P$. 

\begin{maintheo}\label{MASS}  Suppose that $P$ is a Delzant polytope in $\R^m$.
Then the expected
mass density  of random $\lcal^2$ normalized polynomials with Newton polytope
$NP$  has 
$\ccal^\infty$ asymptotic expansions  of the form:
$$ \E_{\nu_{NP}}\left(|f(z)|^2_\FS\right) \sim \left\{  \begin{array}{l} 
\frac{\prod_{j=1}^m (pN+j)}{\# (NP)} =
c_0 + c_1N^{-1} + c_2N^{-2} + \cdots,  \quad \mbox{for }\ z
\in \acal_P=
\mu_\Si^{-1}(\frac{1}{p}P^\circ)
\\[10pt] N^{-s/2} e^{-N b(z)}\left[c_0^F(z) +
c_1^F(z)N\inv+ c_2^F(z)N^{-2}+\cdots\right], \quad \mbox{for } z\in
\rcal_F^\circ
\end{array} \right. $$where 
$c_0^F$ and $b|_{\rcal_F^\circ}$ are positive $\ccal^\infty$ functions on
$\rcal_F^\circ$, and
$s=\codim F$ (for each face
$F\subset\Si^\circ$).  Furthermore, the remainder estimates in the expansions are
uniform on compact subsets of the open regions $\rcal_F^\circ$ and of $\acal_P$.
\end{maintheo}

We shall use (\ref{2measures}) to restate Theorem \ref{MASS} in terms of asymptotic
formulae for the conditional \szego kernel on the diagonal for which we give more
precise asymptotics in Proposition~\ref{SZEGO}. In \cite{SZ2} we will extend the
result to the non-Delzant case as well. Heuristically, the mass concentration in
$\acal_P$ can be understood as follows:  the mass of a monomial $m_{\alpha}$
concentrates (exponentially) on the torus $\mu_\Si^{-1}(\alpha)$. The
constraint $\alpha \in NP$ thus concentrates all the mass of
polynomials with fixed Newton polytope $N P$ in $\acal_P$, with
exponentially small errors, and the mass there is
uniformly distributed with small errors.   

In particular,  $$ \E_{\nu_{NP}}\left(|f(z)|^2_\FS\right) \longrightarrow
\left\{ 
\begin{array}{ll} 
\frac{p^m}{\vol(P)}\,,  \quad  &\mbox{for }\ z
\in \acal_P\\[8pt] 0\,, & \mbox{for } z\in
\C^{*m}\sm\overline{\acal_P}\end{array} \right.\ , $$
as illustrated  by the following graphs,
plotted using Maple, for the case where $P$
is the unit square. (Recall Figure \ref{allowed-square} for the depiction of
$\acal_P$ for this case.)

\begin{figure}[htb]
\centerline{\includegraphics*[bb= 164 295 447 530]{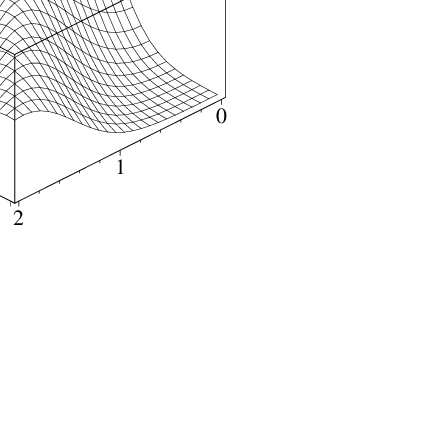}\hspace{-2cm}
\includegraphics*[bb= 164 295 447 530]{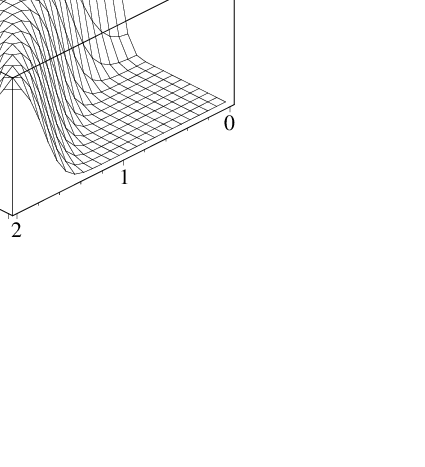}}
\caption{$\frac{1}{4}\E_{\nu_{NP}}\left(|f(z)|^2_\FS\right)$ for
$P=[0,1]\times [0,1]$}
\label{mass-graph}\end{figure}

The decay function
$b$ in Theorem
\ref{MASS} is analogous to the Agmon action to the allowed region.  We shall show in
Proposition~\ref{SZEGO} that $b$  is  $\ccal^1$ (but not $\ccal^2$) on all
of $\C^{*m}$, and we derive several formulas for $b(z)$
((\ref{b}), (\ref{b-action}), (\ref{CRITVAL})). To describe our
first formula, we associate with each point
$z\in\C^{*m}\sm\acal_P$ a unique point $\xi\in \d \acal_P$ of the
form $\xi=e^{\tau/2}\cdot z$, where  $-\tau$  is in the (real)
normal cone to the convex set $P$ at the point $p\mu_\Si(\xi)\in\d
P$. Here $r \cdot z = (r_1z_1,\dots,r_mz_m)$ denotes the $\R_+^m$
action on $\C^{*m}$. (The normal cone to a point $q\in \d P$,
described in \S \ref{s-fans}, is the element of the fan associated
to the face of $P$ containing $q$.)  We write $\tau_z=\tau$,
$q(z)=p\mu_\Si(\xi)$; these vectors are then given by the
conditions:
\begin{itemize}
\item
$\
\de\mu_\Si(e^{\tau_z/2}\cdot z) =q(z)\in\d P$;
\item $\ -\tau_z$ is in the normal cone to $P$ at $q(z)$.\end{itemize}
The existence and uniqueness of $q(z)$ and $\tau_z$ are stated in Lemma
\ref{claim2-3}.

Our first formula for $b$ is
:\begin{equation} \label{b}b(z)= \langle q(z), \tau_z\rangle + \de
\log\left(\frac{1+\|z\|^2}{1+\|e^{\tau_z/2}\cdot z\|^2}\right) \qquad
\mbox{for }\ z\in\C^{*m}\sm\acal_P\;.\end{equation} We will obtain (\ref{b}) from the
integral formula,
\begin{equation}\label{b-action} b(z) =
\int_0^{\tau_z}\left[q(e^{\sigma/2}\cdot z) - \de\mu_\Si(e^{\sigma/2}\cdot
z)\right] \cdot d\sigma \end{equation} (where $\int_0^{\tau_z}$ denotes the integral
over any path in $\R^m$ from $0$ to $\tau_z$),  which could be interpreted as an action,
thereby bringing the results closer to classical Agmon estimates.

We also show that
$\frac{1}{N}\log
\E_{|NP}\left(|f(z)|^2_\FS\right)
\to - b(z)$ uniformly on compact subsets of $\C^{*m}$
(Proposition~\ref{convergence}).

Our final results concern the zero set of one  or more random polynomials.  We shall
consider not only the expected volume density of the zero set, but also the more precise
description involving the expected zero current. The zero current $Z_{f_1, \dots, f_k}$
of $k$ polynomials (or more generally,
$k$ holomorphic functions)
$f_1,\dots,f_k$ is the current of integration over the
zero set
$$|Z_{f_1, \dots, f_k}|:=\{z\in\C^{*m}:f_1(z)=\cdots =f_m(z)=0\}\;.$$ This
current is given by

$$(Z_{f_1, \dots, f_k}, \phi) := \int_{|Z_{f_1, \dots, f_k}|}\phi\;,\qquad
\mbox{for a test form }\
\phi\in
\dcal^{m-k,m-k}(\C^{*m})\;,$$ where integration is over the set of regular
points of $|Z_{f_1, \dots, f_k}|$ (provided that it is of codimension
$k$ and without multiplicity, which is almost surely the case).

We denote by
 $\E_{|P }(Z_{f})=\E_{\ga_\de|P }(Z_{f})$ the
conditional expectation of the zero current of a random polynomial $f
\in H^0(\CP^m, {\mathcal O}(\de), P)$  with Newton polytope $P$.  In fact,
$\E_{|P }(Z_{f})$ is actually a smooth $(1,1)$-form on $\C^{*m}$ (Proposition \ref{EZ}).

Let us recall what happens when $P=\de\Si$ (cf.\ \cite{SZ}).
By the uniqueness of Haar measure, the expected zero current $\E
(Z_{f})$ taken over all polynomials of degree $p$
is given by $p\om_\FS$, where $\om_\FS=\frac{i}{2\pi}
\d\dbar \log\|z\|^2$ is the Fubini-Study \kahler form on $\CP^m$.
Thus the expected distribution of zeros, as well as the tangent to the zero
varieties, is uniform over $\CP^m$. We now consider how the expectation changes
if we add the condition that $P_f = P$.

\begin{maintheo} \label{main}  Let $P$ be a Delzant polytope.  Then there exists
a closed semipositive $(1,1)$-form $\psi_P$ on $\C^{*m}$ with
piecewise $\ccal^\infty$ coefficients such that:
\begin{enumerate}
\item[i)\ ]  $N^{-1}\E_{|NP} (Z_f)\to \psi_P$ \
in  $\lcal^1_{\rm{loc}}(\C^{*m})$.
\item[ii)\ ] $\psi_P =\de\om_\FS$ on the classically allowed region
$\mu_\Si\inv(\frac{1}{p}P^\circ)$.

\item[iii)\ ] On each region $\rcal_F^\circ$, the $(1,1)$-form
$\psi_P$ is $\ccal^\infty$ and has constant rank equal to $\dim F$; in
particular, if
$v\in\Si^\circ$ is a vertex of $\frac{1}{\de}P$, then
$\psi_P|_{\rcal_v^\circ}=0$.
\end{enumerate}
\end{maintheo}

As a corollary, we obtain some statistical results on the so-called `tentacles'
of amoebas in dimension $2$ (see \S \ref{AM}). Roughly speaking, the (compact) amoeba
of a polynomial $f(z_1,z_2)$ is the image of the Riemann surface $Z_f$ under the 
moment map
$\mu_{\Sigma}$ on
$(\C^*)^2$, and the tentacles are the ends of the amoeba. Certain tentacles must end
at vertices of the triangle $\Si$ while others are `free' to end anywhere along the
boundary of
$\Si$. In Corollary \ref{amoeba}, we will prove
that (in the limit $N \to \infty$)  almost all of the free tentacles of  typical amoebas tend to end in the classically
allowed portion of $\d\Si$.

We call the form $\psi_P$ in Theorem \ref{main} the {\it limit expected zero
current.\/} By
$\lcal^1_{\rm{loc}}$ convergence in (i), we mean
$\lcal^1_{\rm{loc}}$ convergence of the coefficients. (Recall that $\E_{|NP}
(Z_{f_1, \dots, f_k})$ is a $(k,k)$-form with smooth coefficients.) If we write
$\psi_P=\sqrt{-1}\sum
\psi_{jk}(z)dz_j
\wedge d\bar z_k$, then
$\big(\psi_{jk}(z)\big)$ is a semi-positive Hermitian matrix, for
non-transition points $z$.  By the rank of $\psi_P$ at $z$, we mean the rank of
the matrix $\big(\psi_{jk}(z)\big)$.  Note that if $\rcal_F$ and $\rcal_{F'}$ are
adjoining regions (i.e., have a common codimension 1 interface), then $F$ and
$F'$ are of different dimensions, so $\psi_P$ must be discontinuous along the
interface. We remark that as an element of $\dcal'{}^{1,1}(\C^{*m})$, the
current $\psi_P$ is closed and positive.

The form $\psi_P$ not only encodes the expected (normalized) density of the
zero set, but also the expected density of tangent directions to the zero set.
In the course of the proof of Theorem \ref{main}, we will show that in the
forbidden region, the limit tangent directions are restricted.  In particular,
as the polytope expands, the tangent spaces to typical zero sets through a
point $z^0\in\rcal_F^\circ$ approach spaces containing the tangent space (at
$z^0$) of the `normal flow' $\{e^{\tau+i\theta}\cdot z^0:
\tau,\theta\in T_F^\perp\subset\R^m\}$. Precise formulations of this fact are
given in Theorem
\ref{more} and (\ref{limtangents}). Thus, while the expected distribution of
zero densities is absolutely continuous, the expected distribution of zero
tangents is singular.

We can also consider $k$ independently chosen random polynomials
$$f_j\in H^0(\CP^m,\ocal(\de_j) ,P_j)\;,\qquad 1\le j\le k$$
($k\le m$), and we let $\E_{|P_1,\dots,P_k}(Z_{f_1,\dots,f_k})$ denote the
expected zero current with respect to the probability measure
$\ga_{\de_1|P_1}\times
\cdots\times
\ga_{\de_k|P_k}$ on the product space. If $P_1=\cdots = P_k=P$,
then we also write
$E_{|P,\dots,P}(Z_{f_1,\dots,f_k}) = E_{|P}(Z_{f_1,\dots,f_k})
$.

\begin{maintheo} Let $P_1,\dots,P_k$ be Delzant polytopes.  Then
$$N^{-k}\E_{|NP_1,\dots,NP_k}(Z_{f_1,\dots,f_k}) \to
\psi_{P_1}\wedge\cdots\wedge\psi_{P_k} \quad
\mbox{in\ \ }
\lcal^1_{\rm{loc}}(\C^{*m})\;,
\quad \mbox{as\ \ } N\to \infty\;.$$ \label{simultaneous}\end{maintheo}

We see from Theorem \ref{main}(iii) that the zero set $|Z_f|$ of a
polynomial with polytope $NP$ almost surely  creeps into the classically
forbidden region
$\mu_\Si\inv(\Si\sm\frac{1}{p}P)$ in the semiclassical  limit $N\to\infty$.  Indeed the
expected volume of the zero set, or more generally the simultaneous zero set of
$k$ polynomials, has the
following exotic distribution law:

\begin{maincor} \label{zerovolumes} Let $P_1,\dots,P_k$ be Delzant polytopes.
Then for any open set $U\subset\!\subset\C^{*m}$,
$$\frac{1}{N^k}\E_{|NP_1,\dots,NP_k} \vol (|Z_{f_1,\dots,f_k}|\cap U) \to
\frac{1}{(m-k)!}\int_U
\psi_{P_1}\wedge\cdots\wedge\psi_{P_k} \wedge\om^{m-k}_\FS\;.$$
\end{maincor}

By the volume of $|Z_{f_1,\dots,f_k}|\cap U$, we mean the $(2m-2k)$-dimensional
volume in $\CP^m$, which is given by
\begin{equation}\label{vol} \vol (|Z_{f_1,\dots,f_k}|\cap U)=
\int_{|Z_{f_1,\dots,f_k}|\cap U}
\textstyle\frac{1}{(m-k)!}\om^{m-k}
=\left(Z_{f_1,\dots,f_k}, \textstyle\frac{1}{(m-k)!}\om^{m-k}\right)\;.
\end{equation}  Corollary \ref{zerovolumes} follows immediately from Theorem
\ref{simultaneous} and (\ref{vol}).

In addition, the following consequence of Corollary
\ref{zerovolumes} and part (iii) of Theorem \ref{main} says that
there are subtler `very forbidden regions' in the case where the polytope has
faces of codimension
$\ge 2$ in the interior of $\Sigma$, namely the regions comprising the normal
flow-out of these faces:

\begin{maincor} \label{subtler} The expected zero current $N^{-k}\E_{|NP}
(Z_{f_1,
\dots, f_k})$ tends to $0$ at all points of each forbidden subregion $\rcal_F$
with
$\dim F <k$.  In particular, $N^{-1}\E_{|NP}(Z_f)$ tends to zero at all
points of the flow-out $\rcal_v$ of a vertex
$v\in\frac{1}{\de}P\cap \Si^\circ$. \end{maincor}

We now say a few words about the proofs. The key result is Theorem
\ref{MASS} on the mass of polynomials with Newton polytope $P$, which we formulate
more precisely in Proposition~\ref{SZEGO}
using the conditional \szego
kernel. As will be shown in \S
\ref{DZ}, we can derive the expected distribution of zeros from our asymptotic
formula for the conditional \szego kernel $\Pi_{|P N}(z, z)$.

 As  in  the proof in
\cite{A,Au} of Kouchnirenko's theorem, a natural approach to our
results on polynomials with prescribed Newton polytope
 is through the \kahler  toric variety $(M_P, \omega_P)$ associated to
$P$. We recall that $M_P$ carries a natural line bundle $L_P$ with
$c_1(L_P) =  \omega_P$. The main  connection is that there
exists a natural identification
\begin{equation} H^0(\CP^m, {\mathcal O}(Np), P) \simeq H^0(M_P, L_P^N)
\end{equation} between polynomials with prescribed Newton polytope and
holomorphic sections of powers $L_P^N$ of $L_P$. The torus $\T =
\{(e^{i\phi_1},\dots, e^{i\phi_m})\}$ acts on $H^0(M_P, L_P^N)$ with character
\begin{equation}
\chi_{NP}(e^{i \phi}) = \sum_{\alpha \in NP}
e^{i \langle \alpha, \phi \rangle}\;,\qquad e^{i\phi}=
(e^{i\phi_1},\dots, e^{i\phi_m})\;. \end{equation} We then have
the simple expression
\begin{equation} \label{CSCH} \Pi_{|N P} (z,z) =  \int_{\T }
\Pi^{\CP^m}_{N\de} (t\cdot z,z) \overline{\chi_{NP}(t)} dt.
\end{equation}
To obtain asymptotics as $N \to \infty$, we  derive  an
oscillatory integral formula.  First we
combine (\ref{CSCH}) with the formula
\begin{equation} \label{SZCH} \chi_{NP}(e^{i \phi}) = \int_{M_P} \Pi_{N
}^{M_P}(e^{i\phi} \cdot w, w)\,  d\vol_{M_P}(w),
\end{equation}  where $\Pi^{M_P}$ is the \szego kernel of $M_P$,
i.e.  the orthogonal projection onto  $H^0(M_P, L_P^N) $ with
respect to the volume form determined by $\omega_P$.
  In order to obtain an  exponential
decay rate $e^{-Nb(z)}$ of Theorem \ref{MASS}, we will need a very
precise understanding of $\Pi_{N }^{M_P}$. This we obtain in \S \ref{SKTV} by
making a special construction of the \szego
kernel of a toric variety which is perhaps of independent
interest. Our formula for $\Pi^{M_P}$ involves the
application of a certain {\it Toeplitz Fourier multiplier} to the
classical \szego kernel for the sphere, pulled
back by a monomial embedding.  After showing in Proposition~\ref{TOEP} that this
Toeplitz-Fourier multiplier is a Toeplitz operator in the sense of
\cite{BG}, we obtain an oscillatory integral  formula
for $\Pi_{|N P}$ (see  (\ref{szego1}) and (\ref{szego3})) in terms of an integral
over $\T\times M_P$ with complex phase of positive type. 

After that, the proof
of Theorem~\ref{MASS} follows from the method of stationary phase for complex
oscillatory integrals. For the case where
$z$ is in the classically allowed region, we easily
find that the critical manifold is given by $\phi=0$ and $\mu_P(w) =p\mu_\Si(z)$,
where $\mu_P$ is the moment map of $M_P$. Since the phase vanishes and has
nondegenerate  normal Hessian along the critical manifold, we
immediately obtain an asymptotic expansion.  The case where $z$ is in the
classically forbidden region is more subtle. Since
$p\mu_\Si(z)$ lies outside of $P$ for this case, the phase has no critical points. 
To complete the analysis, we must deform the contour to pick up critical
points.  The main part of the proof of Theorem~\ref{MASS} involves finding a
contour for which there is a nondegenerate critical manifold where the phase has
maximum real part. In particular, we consider the complexification $\C^{*m}$ of $\T$
and deform $\T$ to a contour of the form $(\log |\zeta_1|,\dots, \log|\zeta_m|)
=\tau\in\R^m$. We then show that the phase takes its maximal real part on a
critical manifold if and only if $\tau$ is the unique vector $\tau_z$
used in formula (\ref{b}) for the decay rate
$b(z)$; indeed
the maximal real part of the phase is $-b(z)$.  Finally we prove that
nondegeneracy of the normal Hessian holds if (and only if) $z$ is not a transition
point, so that we obtain  the asymptotic expansion of Theorem~\ref{MASS}. 

We would like to remark here that other formulae for the lattice sums 
$\chi_{P}(t)$ are well known from the work of Khovanskii-Pukhlikov \cite{KP},
Brion-Vergne \cite{BV} and others.
An alternative approach is to replace (\ref{SZCH}) with  the formula
from
\cite[Theorem~3.12]{BV}:
 \begin{equation} \chi_P(e^{i \phi}) = \mbox{Todd}(\partial/\partial h)
\left. \left(\int_{P(h)} e^{i \langle x, \phi \rangle} dx\right)\right|_{h = 0},
\end{equation}
where
$P(h) = \{x: \langle u_j, x\rangle + a_j + h_j \geq 0,\ 1 \leq j \leq n\}$, and 
Todd$(\partial/\partial
 h)$ is a certain infinite order differential operator known as a {\it  Todd operator}.
Upon dilating the polytope, one obtains
\begin{equation}\label{CHARACTER} \chi_{N P}(e^{i \phi}) =  N^m  \;
\mbox{Todd}(N\partial/\partial h)
\left. \left(\int_{P(h)} e^{i N\langle x, \phi \rangle} dx\right)\right|_{h = 0}.
 \end{equation}
The asymptotics of mass and zeros  in the classically allowed
region follow very easily from this expression for $\chi_{N P}$,
and with further effort one can prove the result in the forbidden
region.  Both
approaches seem to us quite interesting and to have their own
advantages. In the sequel \cite{SZ2}, we will use the approach
based on (\ref{CHARACTER}), which seems more efficient when the
toric varieties have singularities.

We end the introduction with some  comments on the relation of our
results to other work on polynomials with a fixed Newton polytope.
Such polynomials are  called {\it sparse} in the literature, and
methods of algebraic (including toric) geometry have recently
been applied to the computational problem of locating zeros of
systems of  such sparse polynomials (e.g., see \cite{HS, MaR, S, 
R,V}).  Our results give information on the expected location of
zeros when the polynomials are given the conditional measure
$\ga_{| N P}$ on $H^0(\CP^m, {\mathcal O}(N), P)$. To our knowledge, the asymptotic
concentration of zeros of sparse systems in the classically
allowed region has not been observed before. It obviously pays
most to search for zeros of systems of $m$ polynomials in $m$
variables in the allowed region.

It should be noted that this asymptotic pattern reflects the choice of the
conditional measures
$\ga_{| N P}$ on $H^0(\CP^m, {\mathcal O}(N), P)$.
In previous work \cite{SZ}, we obtained the  expected
and almost everywhere distribution of zeros of several holomorphic sections
of positive line bundles $L \to M$ over general \kahler manifolds.
Those
results specialize to our current setting where $M = M_P, L = L_P$ and
also give results on the distribution of zeros of polynomials in  $H^0(\CP^m, {\mathcal O}(N), P)$. However,
in that paper we defined different Gaussian probability measures $\ga_N^{M_P}$ on
$H^0(M_P, L_P^N)$, namely those induced by the $\lcal^2$ inner product induced
by a hermitian metric on $L_P$  and its curvature form $\omega_P$. In
Proposition~\ref{toriczero} of \cite{SZ}, we showed that  the expected distribution of zeros relative to the $\ga_N^{M_P}$ satisfies
$  \frac{1}{N^k } \E_{\ga_N^{M_P}}  (Z_{f_1,
\dots, f_k}) = \om_P^k +
O(\frac{1}{N})$.
The measures  $\ga_N^{M_P}$ have more recently been
studied by  \cite{MaR},  who also obtain (among other things) the formula for the expected
distribution.
Our present
measures $\ga_N^{M_P}$ are quite singular relative to $\ga_{| N P}$ in the limit
as $N \to \infty$.   Our point of view is that the conditional measures $\ga_{| N P}$ on
$H^0(\CP^m, {\mathcal O}(N), P)$ are the  natural ones when comparing polynomials in
 $H^0(\CP^m, {\mathcal O}(N), P)$ as $P$ varies. The toric variety $M_P$ is a
technical device for studying $\Pi_{| NP}$ but does not figure into the statements
of our main results.
The difference in the two measures is that in the conditional case, we fix the
$\lcal^2$ norms of the monomials $z^{\alpha}$ in advance (as their $\CP^m$ norms),
while in the $\gamma^{M_P}_N$ case the norms vary as $P$ and $M_P$ vary. The latter
would seem to create complicated  biases  towards some monomials and away from
others as
$P$ varies and make it difficult to understand what such a comparison is measuring.

\section{Background on toric varieties and moment polytopes}

\subsection{Newton polytopes and toric varieties}

The space   $H^0(\C\PP^m,\ocal(p), P) $ of polynomials with Newton
polytope $P$ may be identified with the space $H^0(M_P, L_P)$ of
holomorphic sections of a line bundle $L_P \to M_P$ over a
  toric variety $M_P$. Under our running assumption that $P$ is Delzant (see below),
  the variety $M_P$ is smooth. Recall that a {\it toric variety\/} is a complex
algebraic variety
$M$ containing the complex torus $$\C^{*m}= (\C\sm\{0\})\times\cdots\times
(\C\sm\{0\})$$ as a Zariski-dense open set such that the
group action of $\C^{*m}$
 on itself extends to a $\C^{*m}$ action on $M$. In the Delzant case, $M_P$ can be
given the structure of a symplectic manifold such that the restriction of the action
to the underlying real torus
$$\T=\{(\zeta_1,\dots\,\zeta_m)\in \C^{*m}: |\zeta_j|=1, 1\le j\le m\}$$ is a
Hamiltonian action (see \S \ref{s-torus}).

We now
review two (well-known) constructions of $M_P$, by fans, resp.\  monomial
embeddings,  to establish notation and clarify the properties we will be
needing. For further details, see  \cite{Da, F,O}. We point out here that
there is also a geometric approach due to Delzant \cite{De} through {\it symplectic
reduction\/}; see
\cite{De, Au, G}.

\subsubsection{Fans}\label{s-fans} 
By a {\it convex integral polytope\/}, we mean the convex hull in $\R^n$ of a finite
set  in the lattice $\Z^m$. 
A convex integral polytope $P$ with $n$ {\it
facets} (i.e. codimension-one faces) can be defined by linear
equations
$$\ell_i(x): = \langle x, u_i \rangle +a_i \geq 0, \;\;\; (i =  1,
\dots, n), $$ where $u_i \in \Z^m$ is the primitive interior
normal to the $i$-th facet.
The polytope $P$ is
called {\it Delzant\/}
 if each vertex  is the intersection of
exactly $m$ facets  whose  primitive normal vectors generate the lattice $\Z^m$.  

For each point $x\in P$, we consider the {\it normal cone to $P$ at $x$\/},
$$C_x: =
\{u\in
\R^m:
\langle u,
x\rangle =\sup _{y\in P}  \langle u,y\rangle\}\;,$$
which is a closed convex polyhedral cone. We decompose $P$ into a finite
union of {\it faces\/}, each face being an equivalence class under the equivalence
relation $x\sim y  \iff C_x = C_y$.  For each face $F$, we let $C_F$ denote the
normal cone of the points of $F$.  {\bf Note that by our convention, the faces are
disjoint sets.}  We shall use the term {\it closed face\/} to refer to the closure of
a face of $P$. 

Each face of dimension $r$ ($0\le r\le m$) is an open polytope in an
$r$-plane in $\R^m$; i.e., the $0$-dimensional faces are the vertices of $P$, the
1-dimensional faces are the edges with their end points removed, and so forth. The  facets $F_i$ and their normal cones are given by: 
$$\bar F_i=\{x\in P:
\ell_i(x)=0\}\;,\qquad C_{F_i}=\{-tu_i:t\ge 0\}\;.$$
The
$m$-dimensional face is the interior  $P^\circ$ of the polytope with normal cone
$C_{P^\circ}=\{0\}$.

A convex integral polytope $P$ determines the fan $\fcal_P:=\{C_F:F\
\mbox{is a face of} \ P\}$.  A {\it fan\/} $\fcal$ in $\R^m$ is a collection of
closed convex rational polyhedral cones such that a closed face of a cone in
$\fcal$ is an element of $\fcal$ and the intersection of two cones in $\fcal$ is
a closed face of each of them.  An example of an integral polytope and its fan is
given in the illustration below.

\begin{figure}[htb]
\centerline{\includegraphics*[bb= 1.0in 7.5in 5in 9.6in]{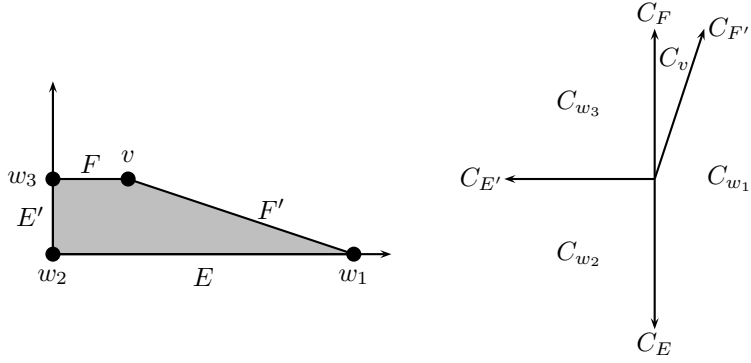}}
\caption{A convex polytope and its fan}
\label{f-fan}\end{figure}

There is a natural way to construct a toric variety from a fan by gluing
together the affine varieties arising from the cones in the fan.  As we will not
use this construction here, we refer the reader to \cite{F} for the details. 
We note that a fan  whose union is all
of
$\R^m$ defines a complete toric
variety, which is projective if and only if $\fcal$ is the fan of a convex
polytope $P$.  If $P$ is Delzant, then the toric variety of the fan
$\fcal_P$ is a smooth projective manifold.  The polytope in Figure \ref{f-fan} is
Delzant if $w_3=(0,1)$, and the corresponding toric variety $M_P$ is a Hirzebruch
surface (see \cite{F}).  We shall study this example in \S \ref{s-ex3}.

\subsubsection{Monomial embeddings}\label{s-monomial}

An alternative definition of  $M_P$ is  through $\T $-equivariant
monomial embeddings. We let
$P\cap\Z=\{\al(1),\al(2),\dots\al(\#P)\}$. If we fix $c = (c_1, \dots, c_{\#
P})\in(\C^*)^{\#P}$, we get a map
$$\Phi_P^c = \big[c_{\al(1)}\chi_{\al(1)},\dots
,c_{\al(\#P)}\chi_{\al(\#P)}\big]:\C^{*m}\to\CP^{\#P-1}
\;;$$
i.e.,
$$ \Phi_P^c(z)= \big[c_{\al(1)}z^{\al(1)},\dots
,c_{\al(\#P)}z^{\al(\#P)}\big]\;.$$
The closure of the image is a toric
variety $M_P^c \subset \CP^{\# P - 1}$.  If $P$ is Delzant, then $\Phi_P^c$ is
injective and we identify $\C^{*m}$ with
its image (the `open orbit') in $M_P^c$.

We refer to the resulting
embedding
\begin{equation} \label{TAM} \Phi_P^c : M_P^c \hookrightarrow \CP^{\# P - 1}
\end{equation} as a monomial embedding, and  we define the line
bundle
\begin{equation} L_P := \Phi_P^{c*} {\mathcal O}(1). \label{LP}\end{equation}
Furthermore,
$H^0(M_P^c, L_P^c)= \Phi_P^{c*} H^0(\CP^{\#P -1}, {\mathcal O}(1))$; i.e., the
sections are spanned by the monomials $\chi_\al$. It follows that
\begin{equation} H^0(M_P^c, L_P) \simeq H^0(\CP^m, {\mathcal O}(p), P)\;.
\end{equation} We further have
(see
\cite{F}):
\begin{equation} H^0(M_P^c, L_P^N) \simeq \Phi_P^{c*}H^0(\CP^{\#P -1}, {\mathcal
O}(N))= H^0(\CP^m, {\mathcal O}(Np), N P)\;.
\end{equation}

Recall that
$H^0(\CP^{\#P-1},\ocal(1))$ has as a basis the linear coordinate functions
$\la_j:\C^{\#P}\to\C$, $1\le j\le {\#P}$.  The {\it Fubini-study metric\/} $h_\FS$ on
$\ocal(1)$ is given by $$|\la_j|_\FS ([\zeta])=\frac{|\zeta_j| }{\|\zeta\|}\qquad
(\zeta\in\C^{\#P})\;,$$ which has curvature form $\om_\FS=\frac{i}{2\pi}\ddbar \log
\|\zeta\|^2$. We  endow
$L_P^c$ with the Hermitian metric
 $h_P^c:= \Phi_P^{c*} h_{\FS}$ of curvature $\omega_P^c$ given on
$\C^{*m}$ by
\begin{equation} \om_P^c=\Phi_P^{c*}\om_\FS=\frac{\sqrt{-1}}{2\pi}\partial
\bar{\partial}
\log
\sum_{\alpha
\in P} |c_{\alpha}|^2 |z^{\alpha}|^2. \end{equation}

Each monomial $\chi_{\alpha}$  with $\alpha \in P$
corresponds to a section of $H^0(M_P^c, L_P^c)$ and vice versa. To
explicitly define this correspondence, we make the
identifications (recalling (\ref{LP})):
\begin{equation}\label{identify}\chi_{\al(j)}\equiv
c_{\al(j)}\inv\Phi_P^c{}^*\zeta_j\in H^0(M_P^c, L_P^c)\;, \quad 1\le j \le
\#P\;.\end{equation}

So far, we have not specified the constants $c_\al$.  For
studying our phenomena, the choice of constants defining the toric variety $M_P$ is not
important. However,  when our polytope $P$ is the full simplex $p\Si$, we shall use the
special choice
$c_\al = {p\choose \al}^{1/2}$, where ${p\choose \al}$ is the multinomial coefficient
defined in the next section.  (We use $c_\al= 1$ in \S \ref{SKTV},
but except for Proposition \ref{PQ}, which is not used in our proof, any choice of
$c$ will work.)

\begin{rem}
For a fixed polytope $P$, the symplectic
manifolds
$(M_P^c, \omega_P^c)$ are equivariantly symplectically equivalent, by Delzant's theorem
\cite{De,G,Au}.  The $M_P^c$ are also equivariantly holomorphically equivalent.
However, in general the maps defining the equivalences are not simultaneously
holomorphic and symplectic (see
\cite{G}).\end{rem}

\subsection{\szego kernels }\label{s-szego} We first describe a general
\szego kernel for a positive Hermitian  line bundle $(L,h)$ on a compact complex
manifold $M$. Let
${\mathcal S}_N\subset H^0(M,L^N)$ denote a subspace of the (finite dimensional)
vector space of holomorphic sections of the
$N^{\rm th}$ power of $L$. We give
$\scal_N$ an inner product.  An important example is the `induced inner
product':
\begin{equation}\langle s_1,\bar s_2 \rangle_{\scal_N} = \int_M \left\langle
s_1(z),\overline{s_2 (z)}\right\rangle_{h_N}d\vol_M(z)\;,\quad\quad
s_1,s_2\in
\scal_N\;.\label{inner}\end{equation}

By the `\szego projector' $\Pi_{\scal_N}$, we mean the orthogonal
projection onto ${\mathcal S}_N$. As  in \cite{Z,SZ,  BSZ}, we
lift  our  \szego kernels on the associated principal $S^1$ bundle
$X \to M$, defined via a Hermitian metric $h$ as follows:   Let
$L^* \to M$ denote the dual line bundle to $L$ with dual metric
$h^*$, and put $X = \{v \in L^*: \|v\|_{h^*} =1\}$. We let
$r_{\theta}x =e^{i\theta} x$ ($x\in X$) denote the $S^1$ action on
$X$.  We then identify sections $s_N$ of $L^N$ with equivariant
functions $\hat{s}_N$ on $X$ by the rule
\begin{equation} \label{sNhat}\hat{s}_N(\lambda) = \left( \lambda^{\otimes
N},
s_N(z)
\right)\,,\quad
\la\in X_z\,,\end{equation} where $\lambda^{\otimes N} = \lambda \otimes
\cdots\otimes
\lambda$; then $\hat s_N(r_\theta x) = e^{iN\theta} \hat s_N(x)$, and we
regard
elements of $\scal_N$ as equivariant functions on $X$.

If $\{S^N_j\}$ denotes
an orthonormal basis of
${\mathcal S}_N$, then the projector $\Pi_{\scal_N}$ is given by the kernel
\begin{equation}\label{szego}\Pi_{\scal_N}(x,y) = \sum_{j = 1}^{k_N} \hat
S_j^N(z)
\overline{\hat S_j^N(y)}\;.\end{equation}
(Note that the usual \szego kernel
for the strictly pseudoconvex boundary $X$ is given by $\sum_{N=1}^\infty
\Pi_{\scal_N}$, where
 ${\mathcal S}_N= H^0(M,L^N)$.)

We now describe  three  different
sequences of
\szego kernels on toric varieties, which play a crucial role in  our main
results:

\subsubsection{The projective \szego kernels}

As a first example, we let $L$ be the hyperplane section
bundle $\ocal(1)$ on $\CP^m$. Recall that the space $H^0(\CP^m,\ocal(p))$ of holomorphic
sections of $L^p=\ocal(p)$ consists of homogeneous polynomials
$$F(\zeta_0, \dots,
\zeta_m) = \sum_{|\la| = p} C_\la \zeta^\la
\qquad (\zeta^\la = \zeta_0^{\la_0} \cdots \zeta_m^{\la_m})$$
in $m+1$ variables.  We can identify $F$  with the (non-homogeneous) polynomial
 in $m$ variables, $$f(z_1,\dots,z_m)=F(1,z_1,\dots,z_m) = \sum
_{|\al|\le p} c_\al z^\al \qquad (z^\al=z_1^{\al_1} \cdots  z_m^{\al_m}),$$
where $c_\al=C_{\hat\al^p}$,
$\hat\al^p=(p-|\al|,\al_1,\dots,\al_m)$.  Inversely, a polynomial $f(z)$ of degree $\le
p$ can be identified with an element of $H^0(\CP^m,\ocal(p))$.

 We give $\CP^m$ the Fubini-Study \kahler form
given in homogeneous coordinates $(\zeta_0,\dots,\zeta_m)$ by
$\om_\FS= \frac{i}{2\pi}\ddbar \log\|\zeta\|^2$, and we   give $\ocal(p)$ the
{\it Fubini-Study metric\/}:
$$|F(\zeta)|_\FS= |F(\zeta)|/\|\zeta\|^p, \quad\mbox{for }\ F\in
H^0(\CP^m,\ocal(p))\;.$$ Identifying  $F$
with the polynomial $f(z)=F(1,z_1,\dots,z_m)$, the Fubini-Study norm can be written
$$|f(z)|_\FS=|f(z)|/(1+\|z\|^2)^{p/2} \qquad (z\in \C^m)\;.$$
 We
equip the space
${\mathcal S}_\de =  H^0(\CP^m, \ocal(p))$ of all homogeneous polynomials of
degree $p$ with the  inner product:
\begin{equation}\label{IP}\langle f, \bar g \rangle = \int _{\CP^m}\left\langle
F,\overline{G}\right\rangle_\FS \,d\vol_{\CP^m}=  \frac{1}{m!}
\int _{\C^m}\frac{\langle
f(z),\overline{g(z)}\rangle}{(1+\|z\|^2)^p}
\,\om_\FS^m(z),\quad f,g\in  H^0(\CP^m, \ocal(\de)).\end{equation}
(We use here the Riemannian volume $d\vol_{\CP^m}=\frac{1}{m!}\om_\FS^m$; note that the
total volume of $\CP^m$ is $\frac{1}{m!}$, using our conventions.)

Under the
$\T $ action, we have the weight space decomposition
 $$H^0(\CP^m, {\mathcal O}(p)) = \bigoplus_{|\alpha|\le p}  \C
\chi_{\alpha}\;,$$ where we recall that $\chi_\al(z) = z_1^{\alpha_1} \cdots
z_m^{\alpha_m}$.
The monomials  $\{
\chi_{\alpha}\}$ are orthogonal but not normalized. Any choice of norming
constants
$\{r_{\alpha} \in \C\}$ will give a monomial basis $\{  r_{\alpha}
\chi_{\alpha}\}$ for
$H^0(\CP^m, {\mathcal O}(p))$. We shall choose
$r_\al=\frac{1}{\|\chi_\al\|}$,  where $$\|\chi_\al\| =\sqrt{ \langle
\chi_\al, \chi_\al\rangle} =
\left[\frac{p!}{(p+m)!{p\choose\al}}\right]^\half\;,\quad\quad
{p\choose\al}:= \frac{p!}{(p-|\al|)!\alpha_1!\cdots
\alpha_m!}\;,
$$ is the Fubini-Study $\lcal^2$ norm of $\chi_{\alpha}$
given by (\ref{IP}).
(See \cite[\S 4.2]{SZ}; the extra factor $m!$ in
\cite{SZ} is due to the use of $\om^m$ instead of $\frac{\om^m}{m!}$
for the volume form.)
This choice provides an orthonormal basis for $H^0(\CP^m, {\mathcal O}(p))$ given by
the monomials
$$\frac{1}{\|\chi_\al\|}\,\chi_\al= \sqrt{\frac{(p+m)!}{p!} {p\choose
\al}}\ \chi_\al\ ,\qquad |\al|\le p\;.$$

We  can identify
$L^*=\ocal_{\CP^m}(-1)$ with
$\C^{m+1}$ with the origin blown up, and the circle bundle $X\subset L^*$ is
unit sphere
$S^{2m+1}\subset \C^{m+1}$.  We let
$\wh\chi_\al:S^{2m+1}\to \C$ denote the equivariant lift of $\chi_\al\in
H^0(\CP^m,\ocal(p))$. We note that
\begin{equation*}\wh\chi_\al(x) =  x^{\hat\al^p}
\;,\end{equation*}
and hence the \szego kernel
$\Pi_p^{\CP^m}$ for the orthogonal projection is
given by:
\begin{equation}\label{szego-proj}\Pi_{p}^{\CP^m}(x,y)=\sum_{|\al|\le
p}\frac{1}{\|\chi_\al\|}\wh
\chi_\al(x) \overline {\wh \chi_\al(y)} =\frac{(p+m)!}{p!}\sum_{|\al|\le p} {p\choose
\al} x^{\hat
\al^p}
\bar y^{\hat\be^p} =\frac{(p+m)!}{p!}\langle x,\bar y\rangle^p
\;,\end{equation} for $x,y\in S^{2m+1}$. (The sum $\sum_{p=0}^\infty \Pi_{p}^{\CP^m}$ is
the usual \szego kernel for the sphere.)

\subsubsection{The conditional \szego kernels associated to a
polytope $P$}

In this case, the relevant space of polynomials is the subspace
$H^0(\CP^m, {\mathcal O}(p), P)\subset H^0(\CP^m,
\ocal(p))$ of polynomials with Newton polytope $P$. Here, we may choose any $p\ge
\deg P:=\max\{|\al|:\al\in P\}$, but we normally choose $p=\deg P$.  In
the conditional
\szego kernel, we retain the Fubini-Study inner product, denoted  $\langle,
\rangle_{|P}$, on this subspace. Hence this example is very similar to the previous one.
 The main difference is that under the  $\T $
action, we have the weight space decomposition
 $$H^0(\CP^m, {\mathcal O}(p), P) = \bigoplus_{\alpha \in P}  \C
\chi_{\alpha}.$$

\begin{defin}The  conditional \szego kernel $\Pi_{|P}$ is the kernel for the orthogonal
projection to
$H^0(\CP^m,\ocal(p),P)$ with respect to the
induced Fubini-Study  inner product:
\begin{equation}\label{Sz2} \Pi_{|P}(x,y) = \sum_{\alpha \in P}
\frac{1}{\|\chi_\al\|^2}
\wh\chi_\al(x) \overline{\wh \chi_{\alpha}(y)} \;.\end{equation}\end{defin}

When defining the term `random polynomial with fixed Newton
polytope $P$', we wish to use an $\lcal^2$-norm on monomials which is
defined independently of $P$. This explains why the conditional
\szego kernel is the essential one in our problem.  The conditional \szego kernel can be
written explicitly as
\begin{equation}\label{Sz2e} \Pi_{|P}(z,\theta;w,\phi) = \frac{e^{ip(\theta-\phi)}
\sum_{\alpha
\in P}{p\choose\al}z^\al \bar w^\al}{(1+\|z\|^2)^{p/2} (1+\|w\|^2)^{p/2}}
\;,\end{equation} where $(z,\theta)$, resp.\ $(w,\phi)$, are local coordinates for $x$,
resp.\ $y$, in $X$; i.e., $$x=e^{i\theta}(1+\|z\|^2)^{-1/2}(1,z_1,\dots,z_m)\;.$$

\subsubsection{The intrinsic \szego kernels of a toric variety}

We now consider \szego kernels which are defined by intrinsic
inner products associated to the toric variety $M_P=M_P^c$. (To simplify notation, we
now drop the $c$.) We regard them as a technical device for obtaining properties of the
conditional \szego kernels.

We thus put ${\mathcal S}_N = H^0(M_P, L_P^N)$, equipped  with the
inner product $\langle, \rangle_{M_P}$ derived from a Hermitian
inner product $h_P$ on $L_P$. 
As on all $S^1$ bundles associated to line bundles,  the $S^1$
action on $X_P$ gives the decompositions
$$\begin{array}{c}\lcal^2(X_P) = \bigoplus_{N = 1}^{\infty}\lcal^2_N(X_P),\;\;
f\in\lcal^2_N \iff f(e^{i \theta} \cdot x) = e^{ i N \theta}
f(x).\end{array}$$ Restricting to the Hardy space, we obtain
$$\;\;\; \\[8pt]
\hcal^2(X_P) = \bigoplus_{N = 1}^{\infty} \hcal^2_N(X_P),\;\;\; \hcal^2_N
=\lcal^2_N \cap \hcal^2.$$ Recall that we have the canonical identifications
$\hcal^2_N(X_P) \equiv H^0(M_P, L_P^N)$ given by identifying a section $s_N$ with the
equivariant function $\hat s_N$ on $X_P$.  We let $\langle,\rangle_{X_P}$ denote the
corresponding inner product on $\hcal^2_N(X_P)$,  and we consider the orthogonal
projectors $\Pi_N^{M_P}(x,y)$ relative to this inner product.

We let $\wh\chi_\al^P$ denote the equivariant lift of
$\chi_\al$ to $X_P$.
Since
$\langle,\rangle_{X_P}$ is invariant under the
$\T $ action on  $X_P$, the
$\{\wh \chi^P_\al\}$ are orthogonal:  indeed,
$$\langle \wh \chi^P_\al,\overline{\wh \chi^P_\be}\rangle_{X_P}=
\langle t^*\wh \chi^P_\al,\overline{t^*\wh \chi^P_\be}\rangle_{X_P} =t^\al\bar
t^\be \langle
\wh \chi^P_\al,\overline{\wh \chi^P_\be}\rangle_{X_P}\;,
\quad \mbox{for all } \  t\in \T \;,$$ and
hence $\langle \wh \chi^P_\al,\overline{\wh \chi^P_\be}\rangle_{X_P} = 0$ unless
$\al=\be$. We then have:
\begin{equation} \label{Sztor} \Pi_N^{M_P}(x,y)  = \sum_{\alpha \in NP}
\frac{1}{\|\wh\chi_\al\|^2_{X_P}}
\wh\chi_\al^P(x) \overline{\wh \chi_{\alpha}^P(y)} \;, \end{equation} where
$\|\wh\chi_\al\|_{X_P}$
is the
$\lcal^2$ norm of
$\wh\chi_\al$ with respect to the inner product $\langle, \rangle_{X_P}$.

In particular, we note that $M_\Si=M_{p\Si}=\CP^m$, and we have the circle bundles
$X_{p\Si}\to
\CP^m$, for $p\ge 1$. When $p = 1$, $X_\Si = S^{2m + 1}$ while for
$p>1$, it is the lens space $X_{p\Si} = S^{2m + 1}/ \{e^{2 \pi i/p}\}$.
The latter statement follows from the fact that homogeneous polynomials of
degree $p$ are well defined on (and separate points of) the quotient by the cyclic group
of
$p^{\rm th}$ roots of unity.

\subsubsection{Powers of a line bundle}

Our main results are asymptotics results in powers of a line
bundle or (equivalently) dilates of the polytope $P$. We therefore
take powers of the line bundles in the previous three examples
obtaining  the spaces  $ H^0(\CP^m, {\mathcal O(Np)}, N P)$ and
$H^0(M_P, L_P^N) \simeq H^0(M_{NP}, L_{NP})$. Hermitian metrics on
the original bundles induce inner products  on these spaces.

To summarize, we have the following three sequences of \szego kernels:
\begin{itemize}

\item $\Pi_{N p}^{\CP^m}\ =\ $  the orthogonal projection to $ H^0(\CP^m,
{\mathcal O(Np
))}$;

\item $\Pi_{|N P}\ =\ $ the orthogonal projection to $ H^0(\CP^m, {\mathcal
O(Np )}, N
P)$;

\item $\Pi^{M_P}_N\ =\ $ the orthogonal projection to $H^0(M_P, L_P^N) $.
\end{itemize}

\noindent We note that for the case $P=\de\Si$, the three sequences coincide.
However, if
$P\,\raisebox{-3pt} {$\buildrel{\subset}\over{\scriptstyle\neq}$}\, \de\Si$,
then the second and third sequences are quite different.

\subsection{Moment maps and torus actions}\label{s-torus}

The group $\C^{*m}$ acts on $M_P^c$ and the subgroup $\T $ acts in a
Hamiltonian fashion. Let us recall the formula for its moment map
$\mu_{P}^c: M_P^c \to \R^m$,  restricted to the open orbit
$\C^{*m}$. This  moment map is the composition
$$\mu_P^c: M_P^c\subset \CP^{\#P -1} \buildrel{\mu_0}\over \to
\R^{\#P } \buildrel{A}\over\to \R^m\;,$$ where
$$\mu_0([z_1, \dots, z_{\#P}]) = \frac{1}{\|z\|^2} (|z_1|^2, \dots,
|z_{\#P}|^2),$$ and $A$ is the linear projection is given by the column
vectors
$(\alpha^1, \dots,
\alpha^{\#P})$.  Hence we have
\begin{equation}\label{muP} \mu_P^c = \frac{1}{\sum_{\alpha \in P}
|c_{\alpha}|^2 |\chi_{\alpha}|^{2 }} \sum_{\alpha \in P} |c_{\alpha}|^2
|\chi_{\alpha}|^{2 }\alpha \;. \end{equation} For any $c$, the image of
$M^c_P$
under
$\mu_P^c$ equals $P$.

For $\al\in P$, we again lift
$\chi_\al\in H^0(M_P^c, L_P^c)$ to an equivariant function $\wh\chi_\al^P$ on the circle
bundle $X_P^c\to M_P^c$,
and we write
\begin{equation}\label{mhat}\wh m_{\al(j)}^P:= c_{\al(j)} \wh
\chi_{\al(j)}^P=
\zeta_j\circ\iota_P\end{equation} where $\iota_P :X_P^c \to
S^{2d +1}$ is the lift of the  embedding  $M_P^c \hookrightarrow
\CP^{\#P-1}$. (Of course,
$\wh m_\al^P$ depends on
$c$, which we omit to simplify notation.) We also consider the monomials
$$m_\al:=c_\al \chi_\al$$ so that $\wh m_\al^P$ is the equivariant lift of $m_\al$ to
$X_P^c$.  In terms
of local coordinates $(z,\theta)$ on
$\pi\inv(\C^{*m})\subset X_P^c$, we have
\begin{equation}\label{equivariantm}
\wh m_\al^P(z,\theta)= \frac{e^{i\theta}c_\al
z^\al}{\left(\sum_{\be\in
P}|c_\be|^2|z^\be|^2\right)^{1/2}}\;.
\end{equation}
Noting that
\begin{equation}\label{sum=1}\sum_{\al\in P}|\wh m_\al^P|^2 =\sum_{j=1}^{\#P}|\zeta_j\circ\iota_P|^2 \equiv 1\;,\end{equation}
we obtain the formula:
\begin{equation}\label{muP2} \mu_P^c(z) = \sum_{\alpha \in P}
|\wh m_\al^P(z)|^2\alpha\;.\end{equation}
(We write $|\wh m_\al^P(z)|=|\wh m_\al^P(z,\theta)|$, since the absolute value is
independent of $\theta$.)

In the case where the polytope is the full simplex $p\Si$,  we shall use the
special choice
$$c^*_\al= {p\choose
\al}^\half=\left(\frac{p!}{(p+m)!}\right)^\half\;\|\chi_\al\|\inv\;,$$ so that
\begin{equation}\label{special} \mu_{p\Si}(z):=\mu_{p\Si}^{c^*}
(z)=  \frac{1}{\sum_{|\alpha|\le p}
{p\choose\al}|z^{\alpha}|^{2 }} \sum_{|\alpha|\le p}
{p\choose\al}|z^{\alpha}|^{2 }\alpha
=\frac{p}{1+\sum|z_j|^2}(|z_1|^2,\dots ,|z_m|^2)\,,\end{equation}
where the last equality follows by differentiating the identity
$(1+\sum x_j)^p=\sum_{|\al|\le p} {p\choose\al}x^\al$.  Note that this choice gives us
the scaling formula
$$\mu_{p\Si}=p\mu_\Si.$$
Furthermore, \begin{equation}\label{mchi}|\wh m_\al^{p\Si}(z)|=\left[
\frac{p!}{(m+p)!}\right]^\half \frac{|\wh
\chi_\al(z)|}{\|\chi_\al\|}\;,\end{equation} where, by abuse of notation, we regard
$|\wh\chi_\al|$ and $|\wh m_\al^{p\Si}|$ as functions on $\C^{*m}$, since they are
invariant under the circle action. Recalling that $M_{p\Si}\equiv S^{2m+1}/\Z_p$, we
have the more precise formula
$$\wh m_\al^{p\Si}(x')=\left[
\frac{p!}{(m+p)!}\right]^\half\frac{\wh \chi_\al(x)}{\|\chi_\al\|}\;,$$ where $x'\in
M_{p\Si}$ is the equivalence class of $x\in S^{2m+1}$,  and therefore by
(\ref{szego-proj}),
\begin{equation}\label{msum}\sum_{|\al|\le p}\wh m_\al^{p\Si}(x') \overline {\wh m_\al^{p\Si}(y')} =\langle x,\bar
y\rangle^p \;.\end{equation}

We now recall a fact on zeros of eigensections which will be
useful later on.

\begin{prop}\label{toricdivisors} Let $\al\in P$.  Then for all $z\in M_P$, we have
$$ \chi_\al(z)=0 \ \Leftrightarrow  \ \mu_P(z)\in
\bigcup
\left\{\bar F: F \ \mbox{\rm is a facet} \ , \al\not\in \bar
F\right\}\;,$$
\end{prop}

\begin{proof} We first verify the implication ($\Leftarrow$): Let
$z\in\mu_P\inv(F)$,
$\al\not\in\bar F$. The moment map (\ref{muP2}) exhibits $\mu_P(z)$
as a convex combination of $\alpha \in P$. Points on $F$ can
however only be convex combinations of lattice points $\alpha \in
\bar F$. Hence the coefficients of $\alpha \notin \bar F$ must
vanish.  Therefore $\chi_\al(z)=0$.

Since the image under the moment map of the support of $\mbox{Div}\,\chi_\al$ is a
union of facets of $P$, to prove the reverse implication, it suffices to show
that for if $F$ is a facet of $P$, then $\chi_\al \not\equiv 0$ on $\mu_P\inv(F)$
whenever
$\al\in
\bar F$. But since $\mu_P$ is surjective, we see from (\ref{muP2}) that this
is true whenever $\al$ is a vertex of $\bar F$.  (Otherwise, points of $F$
near $\al$ would not be in the image.) Now suppose that $\al\in\bar F$ is
arbitrary. We can write
$n_0\al = \sum_1^m n_j v_j$, where the $v_j$ are the vertices of $\bar F$,
$n_j\in \N$ and $\sum n_j=n_0\ge 1$.  Then
$\chi_\al^{n_0} =\chi_{v_1}^{n_1}\cdots \chi_{v_m}^{n_m}\not\equiv 0$ on $\mu_P\inv(F)$.
\end{proof}

\begin{rem}  Recall that the polytope $P$ is defined by the equations
$\ell_j\ge 0$, where
$\ell_j(x)=\langle u_j,x\rangle +a_j$, with $u_j$ the primitive inward
normal to the closed facet
$\bar F_j=\{\ell_j=0\}$.  A more precise statement of Proposition \ref{toricdivisors}
is
$$\mbox{Div}\, \chi_\al = \sum_j \ell_j(\al)D_j\;,  \mbox{where } \ D_j=
\mu_P\inv(\bar {F_j})\;,$$ for all
$\al\in\Z^m$. In particular, letting $\al=0$, we see that $L_P$ is the line
bundle of the
divisor
$D_P=\sum a_jD_j$.\end{rem}

\medskip We now review how the action of the real torus ${\bf T}^m$ lifts from
$M_P^c$ to
$X_P^c$ and combines with the $S^1$ action to define a ${\mathbf
T}^{m+1}$ action on $X_P^c$.
We use the monomial embedding approach of \S \ref{s-monomial}: Recall that under
the monomial embedding
$$\Phi_P^c:\C^{*m}\hookrightarrow M_P^c\hookrightarrow \CP^{\#P-1},\qquad
z\mapsto \big[c_{\al(1)}z^{\al(1)},\dots
,c_{\al(\#P)}z^{\al(\#P)}\big]\;,$$ the $\T$ action on $M_P^c \subset \CP^{\#P-1}$
is given by \begin{equation}\label{actionM}e^{i\phi} \cdot
[\zeta_1,\dots,\zeta_{\#P}]=
\big[e^{i\langle\al(1),\phi\rangle}\zeta_1,\dots,
e^{i\langle\al(\#P),\phi\rangle}\zeta_{\#P}\big]\;.\end{equation}
The action (\ref{actionM}) lifts to an action  on $L_P\inv$:
\begin{equation}\label{actionX}e^{i\phi} \cdot \zeta = \big( e^{i\langle\al(1),\phi\rangle}\zeta_1,\dots,
e^{i\langle\al(\#P),\phi\rangle}\zeta_{\#P}\big) \;.\end{equation}
Since the circle bundle $X_P^c\subset S^{2\#P-1}$ is invariant
under this action, (\ref{actionX}) also gives a lift of the action
(\ref{actionM}) to $X_P^c$.

We also have the standard circle action on $X_P^c$:
\begin{equation}\label{circle} e^{i\theta}\cdot \zeta =  e^{i\theta}
\zeta\;,\end{equation} which commutes with the $\T$-action (\ref{actionX}). Combining
(\ref{actionX}) and (\ref{circle}), we then obtain a
${\mathbf T}^{m+1}$-action on $X_P^c$:
\begin{equation}\label{actionX+}(e^{i\theta}, e^{i\phi_1},\dots, e^{i\phi_m})\bullet
\zeta =  e^{i\theta}(e^{i\phi}\cdot \zeta)\;.\end{equation}

The torus action on $X_P^c$ can be quantized to define an
action of the torus as unitary operators on
$\hcal^2(X_P^c)= \bigoplus_{N=0}^{\infty}\hcal^2_N(X_P^c)$.  Specifically,
we let $\Xi_1,\dots\Xi_m$ denote the differential operators on $X_P^c$
generated by the $\T$ action:
\begin{equation}\label{Xi}(\Xi_j \hat S)
(\zeta)=\frac{1}{i}\frac{\d}{\d\phi_j}  \hat
S(e^{i\phi}\cdot\zeta)|_{\phi=0}\;,\quad \hat S\in \ccal^\infty(X_P^c)\;.
\end{equation}

\begin{prop} For  $1\le j\le m$,
\begin{itemize}

\item[\rm (i)] $\ \Xi_j:\hcal^2_N(X_P^c) \to \hcal^2_N(X_P^c)$;

\item[\rm (ii)] The lifted monomials $\hat\chi_{\alpha}^P \in \hcal^2_N(X_P^c)$ satisfy
$\Xi_j
\hat\chi_{\alpha}^P =
\alpha_j \hat\chi_{\alpha}^P$ ($\al\in NP$).

\end{itemize}
\end{prop}

\begin{proof}Item (i) follows from the fact that the $\T$ action is holomorphic and
commutes with
$\frac{\d}{\d\theta}$. For the case $N=1$, (ii) follows immediately from (\ref{mhat})
and (\ref{actionX}).  For $\al\in NP$, $N>1$, we write $\al=\be^1\cdots\be^N$ with
$\be^k\in P$, and the conclusion then follows from the first case and the product
rule.\end{proof}

Furthermore, we recall that
\begin{equation}\label{dtheta} \frac{\d}{\d\theta} :\hcal^2_N(X_P^c) \to
\hcal^2_N(X_P^c)\;,\qquad \frac{1}{i}\frac{\d}{\d\theta}\hat s_N = N\hat s_N
\quad \mbox{for }\ \hat s_N\in \hcal^2_N(X_P^c)\;.\end{equation}

\begin{rem}  The vector fields $\Xi_j$ can be constructed geometrically as follows (see
\cite{G}): Let $\xi_j=\frac{\partial}{\partial\phi_j}$ ($1\le j\le m$) denote the
Hamiltonian vector fields generating the $\T$ action on $M_P^c$.  There is a natural
contact 1-form
$\al$ on $X_P^c$ determined by the Hermitian connection;  a key property of $\alpha $
is that $d \alpha = \pi^*
\omega$  (see \cite{Z, SZ}).
We use $\al$ to define the horizontal lifts of the
Hamilton vector fields $\xi_j$:
$$\pi_* \xi^h_{j} = \xi_j,\;\;\; \alpha(\xi^h_j) = 0.$$
The vector fields $\Xi_j$ are then given by: \begin{equation*}  \Xi_j =
\xi^h_j + 2
\pi i
\langle \mu_P^c\circ\pi, \xi_j^* \rangle \frac{\partial}{\partial \theta}
=\xi^h_j + 2\pi i
(\mu_P^c\circ\pi)_j\, \frac{\partial}{\partial \theta}.
\end{equation*} \end{rem}
(Here, $\xi_j^* \in \R^m$ is the element of the Lie algebra of $\T$
which acts as  $\xi_j$ on $M_P$.)

\section{\label{SKTV} \szego kernels on toric varieties}

As mentioned in the introduction, our analysis of the conditional
\szego kernel is by means of the integral formula (\ref{SZCH})
relating it to the intrinsic \szego kernels $\Pi_N^{M_P}$ of
$M_P.$
The purpose of this section is to give a special construction of
the \szego kernels of a toric variety  which is sufficiently
precise for our applications.

Our first step, Proposition \ref{PQ}, is to  give an
exact formula for $\Pi_N^{M_P}$ as the composition of a certain
{\it Toeplitz-Fourier multiplier} denoted $({\mathcal P}{\mathcal Q})^{-1}$
and the pull back of the Fubini-Study kernel under a monomial
embedding.
The pulled-back Fubini-Study kernel is very simple to analyze and
is the raison d'etre of our method.  The Toeplitz-Fourier multiplier
requires more work. In Proposition \ref{TOEP}, we prove that this operator
is a ${\bf T}^{m+1}$-invariant  {\it Toeplitz operator} of order
$m $ on $X_P$ (modulo a smoothing operator). We apply Proposition
\ref{TOEP} to obtain the goal of this section, Lemma~\ref{aha},
which allows us to use the method of stationary
phase to obtain the mass asymptotics of Theorem
\ref{MASS}.

\subsection{The exact formula}\label{s-exact}

Our exact formula for the \szego kernel of a toric
variety involves two ingredients: the first is the kernel
\begin{equation}
\wt \Pi_1^{M_P}(x,y):=\frac{1}{\delta+1} \Pi_1^{\CP^{\delta}}(\iota_P(x),\iota_P(y))=
\langle\iota_P(x),\overline{\iota_P(y)}\rangle = \sum_{\alpha
\in P} \wh m_\al^P(x) \overline{\wh m_\al^P(y)}\;,\label{Pi1}\end{equation} where
$\delta=\#P-1$ and
$\iota_P :X_P \to S^{2\delta+1}$ is the lift of the monomial embedding  $M_P
\hookrightarrow
\CP^{\delta}$, as described in \S \ref{s-torus}.We shall choose the constants
$c_\al = 1$ for all $\al\in P$, so that $\wh m_\al^P = \wh\chi_\al^P$.  (Except for
Proposition~\ref{PQ}, all our results hold without change for arbitrary $\{c_\al\}$.)
For simplicity of notation, we shall write
$$\Pi_1^N(x,y)=\big[\wt \Pi_1^{M_P}(x,y)\big]^N.$$  We note that
$\Pi_1^N(x,x)=1.$

Second, we introduce  two sequences of functions ${\mathcal P}_N$
and ${\mathcal Q}_N$ on $N P$ for each $N = 1, 2, \dots:$

\begin{enumerate}

\item [(i)] The `partition function' ${\mathcal P}_N(\alpha)  = \# \{(\beta_1,
\dots, \beta_N): \beta_j \in P, \beta_1 + \cdots + \beta_N =
\alpha \},$ where $ \alpha \in N P. $
 It counts the number of ways of writing a lattice
point in $N P$ as a sum of $N$ lattice points each in $P$.

\item [(ii)] The norming function: $ {\mathcal Q}_N(\alpha):=  \int_{X_P}
|\wh\chi_\al^P(x)|^2d\vol_{X_P}(x).$

\end{enumerate}

In our application, we only consider the product
$${\mathcal P}_N {\mathcal Q}_N(\alpha) = \int_{M_P} \int_{{\bf
T}} \Pi_1^{N}(z, e^{i \phi} z) e^{- i \langle \alpha, \phi \rangle
} d \phi dV(z). $$

We observe that each sequence of  functions can be re-defined as a
single function on the homogenized lattice cone
$\Lambda_P=\bigcup_{N = 1}^{\infty} \wh{NP}$. By definition,
$\wh{NP} \subset \Z^{m+1}$ is obtained by  homogenizing  $\alpha
\in N P $ to $\wh\alpha^{Np}:=(Np-|\al|,\al_1,\dots,\al_m)$. We
then define $\pcal,\qcal:\La_P\to \R_+$ by:
  $${\mathcal P}(\wh\alpha^{N}):= {\mathcal
P}_N(\alpha), \;\;\; {\mathcal Q}(\wh\alpha^{N}) = {\mathcal
Q}_N(\alpha)\;.$$

Since the monomials $\qcal_N(\al)^{-\half}\wh \chi^P_\al$ form a complete orthonormal
basis of $\hcal^2(X_P)$, they define the eigenvalues of a special kind of
operator on $\hcal^2(X_P)$.

\begin{defin} An operator $\Pi F \Pi $ on $\hcal^2(X_P)$ will be called a
{\it Toeplitz-Fourier multiplier} if it satisfies the following
(equivalent) conditions:

\begin{itemize}

\item $F$ may be expressed as a function $F(D)$ of the commuting
system of operators $D = (\Xi_1, \dots, \Xi_m,
\frac{\partial}{\partial \theta})$.

\item Its eigenfunctions are the monomials
$\wh\chi^P_{\alpha}$.

\end{itemize}

\end{defin}

Thus, we put:

\begin{equation} {\mathcal P}(D) \wh \chi_\al^P = {\mathcal P}(\wh\alpha^N)
\wh \chi_\al^P,\;\;\; {\mathcal Q}(D) \wh \chi_\al^P = {\mathcal Q}
(\wh\alpha^N) \wh \chi_\al^P,\;\;\; \alpha \in N P.
\end{equation}

The following  explicit {\it factorization formula} exhibits  the \szego
kernels $\Pi_N^{M_P}$ as the composition of two  simple operators. It is
the basis for our analysis of the analytic continuation $\Pi_N^{M_P}$.

\begin{prop}\label{PQ} We have:
$$\Pi_N^{M_P} = ({\mathcal P} {\mathcal Q})^{-1} \Pi_1^N. $$ \end{prop}

\begin{proof}

First, we have by definition, $$\Pi_N^{M_P}(x, y ) = \sum_{\alpha
\in N P} \frac{1}{{\mathcal Q}_N(\alpha)}  \wh\chi_\al^P(x)
\overline{\wh\chi_\al^P(y)}\;.$$

On the other hand, by definition of the partition function, we
also have
$$ \Pi_1^N (x,y) = \sum_{\alpha \in N P} {\mathcal
P}_N(\alpha)  \wh\chi_\al^P(x)
\overline{\wh\chi_\al^P(y)}\;.
$$
We note that the $N^{\rm th}$ power of $\Pi_1^{M_P}$ gives the sum over
the correct set of exponents but does not have the correct
normalizing coefficients. We need to divide each term by
${\mathcal P}_N(\alpha) {\mathcal Q}_N(\alpha)$ to adjust the
coefficients. That is just what the Proposition claims.
\end{proof}

Since
$$({\mathcal P} {\mathcal Q})^{-1} \Pi_1^N = \Pi_N^{M_P} ({\mathcal P} {\mathcal Q})^{-1} \Pi_N^{M_P} \Pi_1^N,$$
we see that  $\Pi^{M_P}_N$ is the composition of the
Toeplitz-Fourier multiplier $\Pi_N^{M_P} ({\mathcal P} {\mathcal
Q})^{-1} \Pi_N^{M_P}$ with the simple kernel $\Pi_1^N$.

\subsection{Toeplitz operator approach to $\Pi_N^{M_P}$}

Our next step is to prove that the Toeplitz-Fourier multiplier $\Pi_N^{M_P} ({\mathcal
P} {\mathcal Q})^{-1}
\Pi_N^{M_P}$ is a more familiar kind of operator, namely a Toeplitz operator.

We begin with some   background on symbols.
We consider symbols $\sigma$  of the form  $\sigma(z, D_{\theta})$,
where $\sigma(z, N)$  is a semiclassical symbol.
 Here,  we say  $\sigma \in S^k(M
\times \N)$ is a semiclassical symbol of order $k$ if
$$\sigma_N(z)=\sigma(z,N) = N^k \sum_{j = 0}^{l-1} a_j(z) N^{-j} +
r_N^l(z),\;\; \mbox{with}\;\; \|D_z^n r_N^l (z)\| \leq C_{nl} N^{k
- l}\qquad (l\ge 1), $$ where $a_j\in\ccal^\infty(M)$. It is a
{\it smoothing symbol\/} if $\|D_z^n \sigma_N (z)\| \leq
C_{n}N^{-l} $ for all $n,l$. The {\it Toeplitz operator\/}
associated to a symbol $\sigma$ is the operator $\Pi {\sigma}(z,
D_\theta)\Pi $, where  $D_{\theta}$ denotes the symbol of
$\frac{\partial}{\partial \theta}$ (the $S^1$ generator).  Its
symbol is the polyhomogeneous function on the symplectic  cone
$\Sigma = \{(x, r \alpha_x) : r > 0\} \subset T^*X$ given by
$$\sigma(z, p_{\theta})  \sim \sum_{j = 0}^{\infty} a_j(z)
D_{\theta}^{k-j} \;. $$ Since $\sigma$ commutes with the $S^1$
action, we have $\Pi \sigma \Pi = \sum_N \Pi_N \sigma_N \Pi_N$.

To
simplify notation, we shall write $\Pi_N= \Pi_N^{M_P}$. The goal of this section is to
prove that $\Pi_N (\pcal\qcal )^{-1} \Pi_N $ is a semi-classical Toeplitz operator in
the following sense:

\begin{prop}\label{TOEP} There exists a symbol $\sigma_N$ of
order $m$
with principal symbol equal to $1$ and a smoothing operator $R_N$ so that
$$\Pi_N (\pcal\qcal )^{-1} \Pi_N = \Pi_N \sigma_N \Pi_N + R_N\;.$$ \end{prop}

As a corollary of Proposition \ref{TOEP}, we obtain the following lemma, which is
the motivation for
 our approach to $\Pi_{|N P}$ through $\Pi_N^{M_P}$. It  provides an
essential ingredient in our derivation  of the mass asymptotics.

\begin{lem}\label{aha}  There exists a symbol $\sigma_N \in S^{m}(M, \N)$
with principal symbol equal to $1$ and a smoothing operator $R_N$ such that
for $t\in \T $, we have
\begin{eqnarray*}\int_{X_P} \Pi_N(x, t \cdot
x) \,d\vol_{X_P} (x) &=& \int_{X_P} \sigma_N(x) \Pi_1^{N}(x, t
\cdot x) \,d\vol_{X_P} (x)\\&&+ \int_{X_P} \int_{X_P} R_N(y,x)
\Pi_1^{N}(x, t \cdot y) \,d\vol_{X_P} (x)\,d\vol_{X_P} (y).
\end{eqnarray*}
\end{lem}

\medskip\noindent{\it Proof (assuming Proposition \ref{TOEP}):\/} Applying
$\Pi_1^N=\Pi_N (\pcal\qcal )
\Pi_N$ to the identity of Proposition
\ref{TOEP}, we obtain
$$\Pi_N  = (\Pi_N \sigma_N \Pi_N + R_N) \Pi_1^N. $$
Now let $T_{t}$ denote the translation operator $f(x)\mapsto
f(t\inv\cdot x)$ on $\lcal^2(X_P)$.
Since $[T_t, \Pi_N] = 0$, we then have
$$T_t\Pi_N= T_t\Pi_N(\sigma_N+R_N)\Pi_1^N =
\Pi_NT_t(\sigma_N+R_N)\Pi_1^N\;.$$ Therefore,
\begin{eqnarray*} &&\int_{X_P} \Pi_N( t\inv \cdot x,
x) \,d\vol_{X_P} (x)\ =\ \mbox{Trace\,}   T_t \Pi_N \ =\
\mbox{Trace\,}\Pi_NT_t(\sigma_N+R_N)\Pi_1^N\\ &&\qquad\qquad=\
\mbox{Trace\,}T_t\sigma_N\Pi_1^N +
\mbox{Trace\,}T_tR_N\Pi_1^N\\&&\qquad\qquad=\ \int_{X_P}
\sigma_N(t\inv\cdot x) \Pi_1^{N}(t\inv\cdot x,x) \,d\vol_{X_P} (x)
\\&&\qquad\qquad\qquad+ \int_{X_P}\int_{X_P} R_N(t\inv\cdot x,y)
\Pi_1^{N}(y,x) \,d\vol_{X_P} (x) \,d\vol_{X_P} (y)\;.
\end{eqnarray*}
The conclusion follows by the change of variables $x\mapsto t\cdot
x$.\qed

\subsubsection{Boutet de Monvel-Sj\"ostrand paramet\-rix}

In proving Proposition \ref{TOEP} and other statements about Toeplitz
operators,  we need some a priori facts
about $\Pi^{M_P}$ which follow from the Boutet de Monvel-Sj\"ostrand
parametrix construction. In this section, we provide the relevant background.

The Boutet de Monvel-Sj\"ostrand para\-metrix is a complex
oscillatory integral formula for the \szego kernels of a strictly
convex CR manifold modulo smoothing kernels \cite{BS}. It can be
used to obtain asymptotics of projection kernels onto holomorphic
sections of powers of a positive line bundle over a \kahler
manifold (see \cite{Z,BSZ}). We recall here the basic properties
of the \szego kernel for $(M_P, L_P)$. Unfortunately, the
remainder term is not sufficiently precise for the tunneling
theory we have in mind. In the next section, we give a special
construction of the \szego kernels of toric varieties which is
sufficiently precise. The construction uses the general properties
established in this section.

Now consider any polarized algebraic manifold $(M, L)$, and assume
$L$ is very ample. Choose a basis $\{S_0, \dots, S_{d}\}$ with $d
= \dim H^0(M, L)$ and let $\Phi$ denote the associated embedding
into projective space. That is, we write $S_j = f_j e_L$ relative
to a local frame $e_L$ and put $\Phi(z) = [f_0(z), \dots,
f_{d}(z)]$.  Recalling that $L=\Phi^*\ocal(1)$, we can equip $L$ with the metric
$\Phi^*h_\FS$.  We also give
$M$ the metric
$\omega =
\Phi^*
\omega_{\FS}$, and we let $\Pi_N$ denote the orthogonal projection
onto $H^0(M, L^N)$ with respect to these metrics. Also, let $\Pi=
\sum_{N=1}^{\infty} \Pi_N$ denote the \szego kernel.

It was proved by Boutet de Monvel and Sj\"ostrand \cite{BS} (see
also the Appendix to \cite{BG})  that $\Pi$ is a complex Fourier
integral operator of positive type,
\begin{equation}\label{FIO}  \Pi \in I_c^0(X \times X, {\mathcal C}) \end{equation}
associated to a positive canonical relation ${\mathcal C}$. For
definitions and notation concerning complex FIO's we refer to
\cite{MS, BS, BG}. The real points of ${\mathcal C}$ form the
diagonal $\Delta_{\Sigma \times \Sigma}$ in the square  of the
symplectic cone
\begin{equation} \label{SYMCONE} \Sigma = \{r \alpha_x: r > 0, x
\in X\},
\end{equation}
where $\alpha$ is the connection form.
 We refer to \cite{BG} (see Lemma 4.5 of the Appendix).
   Moreover, in \cite{BS}  a parametrix is constructed for $\Pi,$
from which it follows (see \cite{Z,BSZ}) that
\begin{equation} \label{oscint}\Pi_N(x,y) \sim N \int_0^{\infty}
\int_0^{2\pi} e^{  N ( -i\theta + t  \psi( r_{\theta} x,y))}
s(r_{\theta} x,y, Nt )\, d\theta\, dt \,,\end{equation}  where
$s(x,y,t ) \sim \sum_{k = 0}^{\infty} t ^{m -k} s_k(x,y)\in
S^m(X\times X\times \R_+)$ is  a classical symbol of order $m$.
Here, `$\sim$' means modulo a rapidly decaying term (i.e., a term
whose $\ccal^j$ norms are $O(N^{-k})$ for all $j,k$).

To describe the phase in (\ref{oscint}), we let $e_L$ be a
nonvanishing holomorphic section of $L$ over an open $U\subset M$,
and  consider the  analytic extension' $a(z,w)$ of $a(z,z):=a(z)=
\|e_L(z)\|_h^{-2}$  in $U\times U$ such that
$a(w,z)=\overline{a(z,w)}$ on $U\times U$. Using coordinates
$(z,\theta)$ for the point $x=e^{i\theta} a(z)^\half e_L(z)\in X$,
we have
\begin{equation}\psi(x_1, x_2) = -1 + e^{i (\theta_1 - \theta_2)}
\frac{a(z_1, {z}_2)}{\sqrt{a(z_1)}
\sqrt{a(z_2)}}\;.\label{psi-general}\end{equation}

Again assuming that  the metric  $\om$ on $M$ is the pull-back of
$\om_\FS$, we claim  that  the phase
equals:
\begin{equation}  \psi_(x,y)= -1 +
\langle \iota(x),\overline{\iota(y)}\rangle\;,\label{psiX}\end{equation}
where $\iota:X\to \C^{2d+1}$ is the lift of $\Phi$ given by
$$\iota(z,\theta)=e^{i\theta}
\left(\sum|f_j(z)|^2\right)^{-\half}\big(f_0(z),\dots,f_d(z)\big)\;.$$

To see this, let us recall the \szego kernel of the hyperplane
section bundle $\ocal(1)\to \CP^d$ over projective space with the
Fubini-Study metric. We take $U=\{z_0\ne 0\}\approx \C^d$, and we
consider the local frame $e =z_0$. Using the local coordinates
$[1,z_1,\dots,z_m]\mapsto (z_1,\dots,z_m)\in\C^d$, we then have
$a^{\CP^d}(z)=\|e(z)\|^{-2}=1+\sum_{j=1}^d |z_j|^2$. Hence
$a^{\CP^d}$ has the real-analytic extension
\begin{equation} \label{a-projective} a^{\CP^d}(z,w)=1+\sum_{j=1}^d z_j\bar
w_j\;,\end{equation}
and (\ref{psiX})  follows. In the case of a toric variety $M_P^c$,
(\ref{psiX}) becomes
$$\psi(x,y) = -1 +\sum_{\al\in
P}\wh m_\al^P(x)\overline{\wh m_\al^P(y)}\;. $$

\subsubsection{Complex Fourier integral operators}\label{pf-TOEP}

We now use the theory of  complex Fourier integral operators
\cite{MS} (or alternatively,  Toeplitz operators \cite{BG})  to
prove that the Toeplitz Fourier multiplier $\Pi ({\mathcal
P}{\mathcal Q})^{-1} \Pi$ is a Toeplitz operator of order $m$
modulo a smoothing operator (Proposition~\ref{TOEP}). Recall that
a Toeplitz operator in the sense of \cite{BG} is an operator of
the form $\Pi A \Pi$ where $A \in \Psi^m(X)$ for some $m$. Here
and in the following, $\Psi^m(X)$ denotes the space of
pseudodifferential operators of order $m$ on $X$. A smoothing
Toeplitz operator is a smoothing operator of the form $\Pi R \Pi$.

It follows then that $\Pi_N = T_N \Pi_1^N$ where $T_N$ is a
Toeplitz operator. Since it is of some independent interest and
requires no extra work, we prove this as a corollary of a general
result valid for any polarized algebraic manifold $(M, L)$. We
assume $L$ is very ample, and let $\iota: M \to \CP^{d}$ denote
the holomorphic embedding induced by $H^0(M, L)$. Let $\Pi_{\iota}
(x, y) = \Pi(\iota(x), \iota(y))$ denote the pullback to $X$ of
the Fubini-Study \szego kernel of $\CP^{d}.$ Also, let $\Pi^M$ or
more simply $\Pi$ denote the \szego kernel for $M$ with respect to
$\iota^* \omega_{\FS}.$

About complex Fourier integral operators $A$ (or Toeplitz
operators), all we need to know are the following basic facts:
\begin{itemize}

\item $A \in I^m_c(X \times X, {\mathcal C})$ possesses a principal symbol $\sigma_A$ which is a
half-density (times a Maslov factor) along the underlying
canonical relation ${\mathcal C}$. (In the Toeplitz case, it is a
symplectic spinor.)

\item  We can compose operators in $I^m_c(X \times X, {\mathcal C})$ on the
left and right by elements $B \in \Psi^k(X)$ and $\sigma_{AB} =
\sigma_A \sigma_B = \sigma_{BA}$. (The same is true of Toeplitz
operators.)

\item If $A \in I^m_c(X \times X, {\mathcal C})$ and if $\sigma_A = 0$, then
$A \in I^{m-1}_c(X \times X, {\mathcal C})$ (and also for
Toeplitz operators).

\item $\Pi$ and $\Pi_{\iota}$ are elliptic in that their symbols
are nowhere vanishing. See \cite{BS} for the symbol in the complex
FIO sense and \cite[\S 11]{BG} for the symbol in the Toeplitz
sense.

\end{itemize}

\begin{lem}\label{PIOTA}  Let $\Pi$ denote the \szego projector
associated to the metric $\iota^* \omega_{\FS}$.  Then, there
exist $A \in \Psi^m(X)$ with $[A, D_{\theta}] = 0$ such that $\Pi
\sim  \Pi  A  \Pi \Pi_{\iota} =  \Pi  A \Pi_{\iota}  $ modulo
smoothing operators $\Pi R \Pi.$
\end{lem}

\begin{proof} As mentioned above (\ref{FIO}),   $\Pi$ is a complex
Fourier integral operator associated to a positive canonical
relation $C$, whose real points form the isotropic relation
$\Delta_{\Sigma \times \Sigma} \subset T^*X \times T^*X$ (the
diagonal).
 We  observe that also $\Pi_{\iota} \in
I^*(X \times X, {\mathcal C})$. This follows immediately from the
fact, show by (\ref{psiX}), that $\Pi$ and $\Pi_{\iota}$ are
complex Fourier integral distributions with precisely the same
phase functions. Since the underlying canonical relation  is
parametrized by the phase, they both belong to the same class of
Fourier integral distributions.

Now, the principal symbol $\sigma_{\Pi}$ of $\Pi$, viewed as a
complex Fourier integral distribution, is a nowhere vanishing
1/2-density on ${\mathcal C}$ which is computed in
\cite[Prop.~4.8]{BS}. Alternatively, viewed as a Toeplitz operator
in the sense of \cite{BG} (see Chapter 11), its symbol is an
idempotent symplectic spinor .
 Similarly, the principal symbol $\sigma_{\Pi_{\iota}}$ of $\Pi_{\iota}$
is the pull back under $\iota$ of the nowhere vanishing  symbol of
$\Pi^{\CP^d}$.

By our normalization, $\Pi_{\iota}$ has order $- m$ (since its
amplitude is a constant independent of $N$). We therefore begin by
seeking $A_0 \in \Psi^m(X)$ such that $[A_0, D_{\theta}] = 0$ and
such that  $\Pi- \Pi A_0 \Pi_{\iota}$ is of order $-1$. We first
find $a_0 \in C^{\infty}(M)$ such that $\sigma_{\Pi} = a_0
\sigma_{\Pi_{\iota}}$ and choose $A_0$ so that $ [A_0, D_{\theta}]
= 0$ and so that  $\sigma_{A_0} = a_0.$ Existence of such an $A_0$
follows by ellipticity of $\Pi_{\iota}$ and by averaging; see also
\cite[Prop.~2.13]{BG}.   Thus,  the principal symbol of order $0$
of $\Pi - \Pi A_0 \Pi_{\iota}$ equals zero, i.e. $\Pi - \Pi A_0
\Pi_{\iota} \in I^{-1}(X \times X, {\mathcal C}).$ We denote its
principal symbol by $\sigma_{-1}.$ We then seek $A_{-1} \in
\Psi^{m}(X)$ so that $ [A_{-1}, D_{\theta}] = 0,$ and so that $\Pi
- \Pi A_0 \Pi_{\iota} - \Pi A_{-1} D_{\theta}^{-1} \Pi_{\iota} \in
I^{-2}(X \times X, {\mathcal C}).$ Here, we note that $\Pi
D_{\theta} \Pi$ is an elliptic Toeplitz operator; since $\Pi
\Pi_{\iota} = \Pi_{\iota}$, the expressions $ D_{\theta}^{-1}
\Pi_{\iota}$ are well-defined.
 It suffices to choose $a_{-1} = \sigma_{A_{-1}} \in
C^{\infty}(M)$ so that $a_{-1} \sigma_{\Pi_{\iota}} =
\sigma_{-1}.$ We continue in this way to obtain $a_{-j} \in
C^{\infty}(M) $  always using that $\sigma_{\Pi_{\iota}}$ is
nowhere vanishing. By a Borel summation argument, we can find $A
\in \Psi^{m}(X)$ with the above commutation properties so that
$\Pi A \Pi - \sum_{j = 0}^M \Pi A_{-j} D_{\theta}^{-j} \Pi$ is a
Toeplitz operator of order $-M - 1.$ Then
$$ \Pi - \Pi A \Pi_{\iota} \in I^{-\infty}(X \times X, {\mathcal C}).
$$
\end{proof}

We remark that an  alternative to the observation that the two
complex FIO's have the same phase function is that, by the  choice
of the \kahler form $\iota^*\omega_{\FS}$,   $\iota$ is a
symplectic (as well as holomorphic)
 embedding $M
\subset \CP^d$. Hence, the  pull back operator $\iota^*$ carries
the class of FIO's in the class of $\Pi^{\CP^d}$ to those in the
class of $\Pi$.

\subsubsection{Semi-classical Toeplitz operators}

To complete the proof of Proposition \ref{TOEP}, we need to relate
the complex FIO $\Pi A \Pi_{\iota} $ to semi-classical Toeplitz
operators. The relevant analysis already exists in \cite{Gu2}.

We will need the following asymptotic formula for symbols. In
the language of Berezin-Toeplitz operators, we are essentially
computing the  Berezin
transform between covariant and contravariant symbols of a
Toeplitz operator (see \cite{RT} for background).

\begin{lem} \label{BEREZIN} Let $\sigma_N \sim N^k \sum_{j = 0}^{\infty} s_{-j} N^{-j}$ be  a   semiclassical  symbol of order $k$.  Then
there exists a complete asymptotic expansion
$$ N^{-m}(\Pi_N\sigma_N\Pi_N)(z,z)=  \sum_{j = 0}^{l-1} b_{-j}(z) N^{k-j} +
r_N^l(z) \qquad (l\ge 1)\;,$$ where $b_0 = s_0,b_{-1} = \Delta s_0 + s_{-1}, \dots$,
and in general where $ b_{-j}$ is a sum of differential operators applied to $s_0,
s_{-1}, \dots, s_{-j}.$ Also,
 $ \|D_z^n r_N^l (z)\| \leq
C_{nl} N^{k - l}$. \end{lem}

\begin{proof} Apply the method of stationary phase using
(\ref{oscint})--(\ref{psi-general}) exactly as in the proof of
\cite[Theorem~1]{Z} (which is the case $\sigma=1$).\end{proof}

We then have:

\begin{lem} \label{PAM} Let $(M, L)$ be a polarized algebraic manifold as
above. Then there exists a symbol $a_N \in S^{-m}(M, \N)$ with
principal symbol equal to $1$ so that
$$   \Pi_1^N = \Pi_N a_N \Pi_N  + R_N\quad \mbox{where }\
\|R_N\|_{\rm HS} = O(N^{-k})\;\;\; \forall k\;. $$ Here,
$\|R_N\|_{\rm HS}^2 = \mbox{\rm Trace}\, R_N^* R_N.$
\end{lem}

\begin{proof}
 By Lemma \ref{PIOTA},  there exists $A \in \Psi^{-m}(X)$
such that  $\Pi_{\iota} = \Pi A \Pi.$ Since $[A, D_{\theta}] = 0$,
it follows by \cite{Gu2} that there exists a symbol $a_N$ with
$\|\Pi_1^N- \Pi_N a_N \Pi_N\|_{\rm HS}=O(N^{-\infty})$.    We may
determine $a_N$ by using Lemma \ref{BEREZIN}. Indeed, we have
\begin{equation}\label{diag} \Pi_1^N(z, z) = 1  \sim \Pi_N a_N \Pi_N
(z,z)\end{equation} modulo functions $r_N$ which tend to zero rapidly in
$\ccal^k(M)$. It follows that $a_0 = 1$ (and then the rest of the coefficients may
be determined recursively, e.g. $ a_{-1} = 0, a_{-2} = - D_2 a_0$, and so on).
\end{proof}

\subsubsection{Completion of the proof of Proposition \ref{TOEP}}
 In the toric case, $\Pi_N (\pcal\qcal )^{-1} \Pi_N$ is the inverse
of\break $\Pi_N (\pcal\qcal ) \Pi_N = \Pi_1^N$ on $\hcal^2(X_P)$. From  Lemma
\ref{PAM} we may write $\Pi_N (\pcal\qcal ) \Pi_N  \sim \Pi_N a_N \Pi_N$
modulo smoothing operators.

We note that we may invert $\Pi a \Pi$ in the class of Toeplitz operators; i.e.,  there
exists a symbol
$\sigma_N$ such that
\begin{equation}\label{INVERSE}  \Pi_N \sigma_N \Pi_N \circ \Pi_N a_N \Pi_N \sim
\Pi_N \end{equation} modulo smoothing operators. Such an inverse
symbol exists since $a_0 = 1$. The algebraic formalism in which the inverse
is calculated is that of $*$-products of semiclassical symbols. We
recall that composition of Toeplitz operators defines a $*$-product on semiclassical
symbols by the formula
\begin{equation} \Pi_N a_N \Pi_N \circ \Pi_N b_N \Pi_N \sim \Pi_N
a_N * b_N \Pi_N. \end{equation} The formula for $*$ may be worked
out directly from the parametrix (\ref{oscint}) and the inverse can be computed from
this formula (see
\cite{Gu2, S}).

By (\ref{INVERSE}) we obtain a
symbol $\sigma_N$ (with principal symbol equal to $1$) such that
$$\Pi_N \sigma_N \Pi_N\Pi_N (\pcal\qcal ) \Pi_N \sim \Pi_N\;.$$
Multiplying both sides by $\Pi_N
(\pcal\qcal )^{-1} \Pi_N$, we conclude the proof.
\qed

\begin{rem}
We emphasize that the distinguishing features of  the toric case in
Proposition \ref{TOEP} and Lemma
\ref{PAM} are the exact factorization and the fact that the operator $\Pi ({\mathcal
P}{\mathcal Q})\Pi$ mediating between $\Pi$ and $\Pi_{\iota}$ is
invertible. This is due to the fact that sections of $L$ generate
the ring $\oplus_{N=1}^{\infty} H^0(M, L^N)$ in the toric case. In
general they do not and the exact representation $\Pi = B
\Pi_{\iota}$ in the toric case  is only valid modulo smoothing
operators. \end{rem}

\subsubsection{Alternate proof of Lemma \ref{PAM}}
One could avoid using the calculus of complex Fourier integral operators
or Toeplitz operators in the proof of Proposition \ref{TOEP} by further developing
the calculus of semi-classical Toeplitz operators directly from the Boutet de Monvel-
Sj\"ostrand parametrix as follows:

We first simplify the expression (\ref{oscint})
by using the complex
method of stationary phase  to eliminate the integrals in  the parametrix. The
critical point set of  the phase $\Phi (\theta, \lambda; x, y):= -i \theta + \lambda \psi(r_{\theta} x, y)$ is
given by
\begin{equation} \left\{ \begin{array}{l}d_{\lambda} \Phi (\theta, \lambda; x, y) =
\psi(r_{\theta} x, y) = 0 \iff  e^{i \theta} \langle \iota(x), \iota( y) \rangle= 1,\\
\\  d_{\theta} \Phi (\theta, \lambda; x, y) = -i + d_{\theta} \lambda \psi(r_{\theta}
x, y) = 0 \iff  \lambda e^{i \theta} \langle \iota(x), \iota( y) \rangle = 1.
\end{array} \right. \end{equation} It is easy to see by the Schwartz inequality that a
real critical point exists if and only if $x = y,$ in which case
$\theta = 0, \lambda = 1$ and we obtain the familiar expansion
along the diagonal. When $x \not= y$ we deform the  contour to
$|\zeta| = e^{\tau}$ ($\zeta = e^{i \theta + \tau}$)  so that  $
e^{i \theta + \tau} \langle \iota(x), \iota( y) \rangle = 1$. This
is possible as long as $\langle \iota(x), \iota( y) \rangle \not=
0$, as happens near the diagonal, where the parametrix is valid.
Because the phase is linear in $\lambda$ it is clear that the
critical point is  non-degenerate if and only if $\langle
\iota(x), \iota( y) \rangle  \not= 0$ and that the Hessian
determinant equals $|\langle \iota(x), \iota( y) \rangle |^2.$ On
the critical set the phase equals $- i \theta + \tau = \log
\langle \iota(x), \iota( y) \rangle,$ hence  we have
\begin{equation} \label{oscintN} \Pi_N(x, y) =
e^{N \log \langle \iota(x), \iota( y) \rangle} S_N(x, y) + W_N(x,y), \end{equation}
where $S_N(x,y) \sim \sum_{k = 0}^{\infty} N^{m - k} S_k(x,y),$ and
where $W_N(x,y)$ is a smooth uniformly rapidly decaying function. Here, we have
absorbed the remainder in the parametrix construction as well as the remainder
in the stationary phase expansion of the parametrix in $W_N$. Note that the
first term may be smaller than the second outside a tubular  neighborhood
of radius $N^{-1/2}$ of the diagonal.

As above, we use Lemma \ref{BEREZIN} to find a symbol $a_N$ so that (\ref{diag})
holds and hence  the kernels  $\Pi_1^N$ and $\Pi_N a_N \Pi_N$  agree on the diagonal
modulo smoothing symbols.
Next, we note that both $\Pi_1^N(x,y)$ and $\Pi_N a_N \Pi_N(x,y)$ are complex
oscillatory functions with  common phase
\begin{equation} \label{PHASE} \Psi(z, w) =  \log \langle \iota(x), \iota(y) \rangle.
\end{equation} In the case of $\Pi_1^N$, this  follows from (\ref{Pi1}). Indeed, we
simply have:
\begin{equation} \label{PiN} \Pi_1^N(x,y) = e^{N  \log \langle \iota(x), \iota(y)
\rangle }. \end{equation} In the case of $\Pi_N a_N \Pi_N$, we apply the method of
complex  stationary phase to the integral formula
\begin{equation} \Pi_N a_N \Pi_N(x,y) \sim \int_{M_P} e^{N \Psi(u; x,y)} S_N(x,u) a_N(u)
S_N(u, y) dV(u), \end{equation}
 coming from (\ref{oscintN}),  with $$ \Psi(u; x,y) = \langle \iota(x), \iota( u)
\rangle + \langle \iota(u), \iota( y) \rangle.$$ It follows that
 there exists an amplitude
$A_N$ defined near the diagonal such that
\begin{equation}\label{oscintS} \Pi_N a_N \Pi_N(x, y) = A_N(x,y) e^{N\log \langle
\iota(x),
\iota( u)
\rangle} + V_N(x,y),\end{equation} where $V_N$ is a new smoothing operator.

Recalling that $X\subset L^*$, we extend $A_N$ to 
$L^*\times L^*$ so that it is of the form
$$A_N(z,\la; w, \la') = (\la\bar\la')^N \tilde A_N(z,w)\;,$$
where we use a local holomorphic frame to write $x = (z, \la),\ y = (w,
\la')\in L^*$.  We note that $\tilde A_N(z,w)$ is  holomorphic in $z$ and
anti-holomorphic in $w$ near the diagonal.  To see this, we first conclude from the
construction in \cite{BS} using the $*$-product that the symbol $s(x,y,t) \sim
\sum_{k = 0}^{\infty} t ^{m -k} s_k(x,y)$ in (\ref{oscint}) extends to a
symbol on $L^*\times L^*$ that  is holomorphic in $x$ and anti-holomorphic in
$y$.  It follows by the stationary phase method described above that the same is true
for the symbol
$S_N(x,y)$ in (\ref{oscintN}), and then by (\ref{oscintS}) that the same is also
true for $A_N$ as claimed.  (Note that in terms of coordinates $(z,\theta)$ on $X$,
the function $z\mapsto A_N(z,\theta;w,\theta')
=e^{iN(\theta-\theta')}\sqrt{a(z)a(w)}\,\tilde A_N(z,w)$ is not holomorphic in $z$.)

The amplitude of $\Pi_1^N - \Pi_N a_N
\Pi_N$ becomes $(\la\bar\la')^N (1 - \tilde{A}_N(z,w))$, where we can write
$1-\tilde A_N \sim B_0 +B_1 N^{-1}+B_2 N^{-2}+\cdots$.  Since $\Pi_1^N - \Pi_N a_N
\Pi_N(z,z)=O(N^{-\infty})$ by our choice of $a_N$ above, it follows that
$B_j(z,z)=0$ (for all $j$).  Since $B_j(z,w)$ vanishes on the diagonal and is
holomorphic in $z$ and anti-holomorphic in $w$, it must be identically 0.  Hence,
$\Pi_1^N -
\Pi_N a_N
\Pi_N(x, y)$ is a smoothing operator.\qed

\bigskip

\section{Mass asymptotics} \label{s-mass}

In this section, we prove a precise asymptotic formula for the conditional \szego
kernel on the diagonal (Proposition~\ref{SZEGO}), which yields the
mass asymptotics of Theorem~\ref{MASS}.  

First, we discuss formula (\ref{2measures}) for
the expected mass density. 
 Let $P$ denote a Delzant polytope in $\R^m$. Recalling the definition (\ref{CG}) of
the conditional probability measure $\gamma_{\de|P}$
on the space $H^0(\CP^m,\ocal(p),P)$ of polynomials with Newton
polytope $P$,
we see that
the expected value of the mass density with respect to $\gamma_{\de|P}$ is given by:
\begin{equation*}
\E_{|P}\left(|f(z)|^2_\FS\right) =
\sum_{\alpha,\be
\in P} \frac{\E(\la_{\alpha}\bar \la_\be
)\chi_{\alpha}(z)\overline{\chi_{\be}(z)}} {\|\chi_{\alpha}\|\|\chi_{\be}\|
(1+\|z\|^2)^{\de/2}}
\;.\end{equation*} Since the $\la_\al$ are independent complex random variables with
variance 1 (i.e., $\E_{|P}(\la_{\alpha}\bar \la_\be)=\delta_\al^\be$), we have:
\begin{equation}\label{Eszego}
\E_{|P}\left(|f(z)|^2_\FS\right)
= \sum_{\alpha \in P} \frac{
|\chi_\al(z)|_\FS^2}{||\chi_{\alpha}||^2}=\Pi_{|P}(z,z)\;.\end{equation}
It then follows by expressing the Gaussian in
spherical coordinates that $$\E_{\nu_{P}}(|f(z)|^2_\FS)  =\frac{1}{\#
P}\E_{|P}\left(|f(z)|^2_\FS\right)= \frac{1}{\# P}\Pi_{|P}(z,z)\;.$$ Replacing $P$ with
$NP$, we obtain formula (\ref{2measures}).

The number of lattice points in the polytope $NP$ is given by
the Riemann-Roch formula (see \cite{F}):
\begin{eqnarray}\# (NP) &=&\dim H^0(\CP^m,\ocal(Np),NP)\ =\
\dim H^0(M_P,L_P^N)\nonumber \\&=&\chi(M_P,L_P^N)\ =\ \sum_{k=0}^m
\deg\big[c_1(L_P)^k\cup
\mbox{Todd}_{m-k}(M_P)\big]
\frac{N^k}{k!}\nonumber \\&=&
 \vol(P) N^m +\dots + \deg \mbox{Todd}_m(M_P)\,.\nonumber\end{eqnarray}
Here, we used the Kodaira vanishing theorem and the fact that
\begin{equation*}\vol(P)=\vol(M_P)=\int_{M_P}\frac{1}{m!}\om_P=
\frac{1}{m!}\deg [c_1(L_P)^m]\;.\end{equation*} (This immediately yields
Kouchnirenko's Theorem, since $\deg [c_1(L)^m]$ equals the number of points in the
intersection of $m$ generic divisors of a very ample line bundle $L$; see
\cite{A}.) The holomorphic Euler characteristic $\chi(M_P,L_P^N)$ is also known as the
Hilbert polynomial of the polarized manifold 
$(M_P,L_P)$.   Combining (\ref{2measures}) and the Riemann-Roch formula, we obtain
\begin{equation}\label{Eszego2}\E_{\nu_{NP}}(|f(z)|^2_\FS) 
=\frac{1}{\chi(M_P,L_P^N)}\Pi_{|NP}(z,z)\;.\end{equation}

The mass asymptotics of Theorem~\ref{MASS} is an immediate consequence of
(\ref{Eszego2}) and the following asymptotic expansion of the conditional \szego kernel
on the diagonal: 

\begin{prop}\label{SZEGO}   Suppose that $P$ is a Delzant polytope in $\R^m$.
Then:
\begin{enumerate}
\item[i)] For $z$ in the classically allowed region $\acal_P$, we have
$$\Pi_{|NP}(z,z) = \prod_{j=1}^m (Np+j) + R_N(z)
\;,\qquad \|R_N\|_{\ccal^k(K)}=O(e^{-\la_K N}) \quad \forall\ k\;,$$
for all
compact $K\subset \acal_P$, where $\la_K>0$.
\item[ii)] On each open forbidden region $\rcal_F^\circ$, 
$$\Pi_{|NP}(z,z)= N^{\frac{m+r}{2}} e^{-N b(z)}\big[c_0^F(z) +
c_1^F(z)N\inv+\cdots + c_k^F(z)N^{-k} +R_k^F(z)\big],$$ where $r=\dim F$ and
\begin{enumerate}
\item \ $c_j^F\in\ccal^\infty(\rcal_F^\circ)$ and $c_0^F>0$ on $\rcal_F^\circ$;
\item \ $\|R_k^F\|_{\ccal^j(K)}=O(N^{-k-1})$, for all compact $K\subset
\rcal_F^\circ$ and for all $j,k$;
\item \ $b>0$ on
$\C^{*m} \sm \overline{\acal_P}$;
\item \ $b$ is given by formula (\ref{b});
\item \  $b\in\ccal^1_\R(\C^{*m})$  (with  $b=0$ on
${\acal_P}$), and $b$ is $\ccal^\infty$ on
each closed region $\overline{\rcal_F}$.
\end{enumerate}
\end{enumerate}
\end{prop}

We remark that the asymptotics of Theorem~\ref{SZEGO} are uniform away from the
transition points only.  To take care of the transition points, we shall also prove
the following local uniform convergence result on all of $\C^{*m}$:  

\begin{prop}\label{convergence}
Let $P$ be a Delzant polytope. Then
$$\frac{1}{N} \log \Pi_{|NP}(z,z)
\to -b(z) $$ uniformly on all compact subsets of $\C^{*m}$.
\end{prop}

This section is devoted to the proofs of Propositions \ref{SZEGO} and
\ref{convergence}.  We begin with part (ii) of
the former, which is the major part of the section.   The proof begins by using Lemma~\ref{aha}
to obtain  an oscillatory integral formula for  the conditional
\szego kernel $\Pi_{|N P}(z, z)$.  We then observe that the integral is rapidly
decaying when $z$ is in the classically forbidden region, while for $z$ in the allowed
region, the phase is of positive type with a nondegenerate critical manifold, and hence
$\Pi_{|N P}(z, z)$ has an asymptotic expansion.  In \S \ref{s-forbidden}, we consider
$z$ in the forbidden region and we seek deformations of the contour of integration so
that the (analytic continuation of) the phase picks up critical points.   In \S
\ref{s-normalbundle}, we formulate necessary geometric conditions for the existence of
critical points where the phase has maximal real part and we show that there is a
contour for which these conditions are satisfied. In \S\S \ref{3.2.2}--\ref{3.2.3}, we
show that the conditions are sufficient and, with the added assumption that $z$ is
not a transition point, the resulting critical points are
nondegenerate (Lemma~\ref{crit}).  We complete the proof of
Proposition~\ref{SZEGO}(ii) in \S\S \ref{3.2.4}--\ref{s-decay}. In \S
\ref{s-decay}, we use \ref{SZEGO}(ii)
to prove Proposition~\ref{convergence}, which we then use to
prove \ref{SZEGO}(i) in \S \ref{3.2.6}. 

We start with a simple  integral formula for the conditional
\szego kernel:
\begin{lem}\label{PIP} For $y_1,y_2\in\C^{*m}$, we have
$$\Pi_{|N P} (y_1,y_2) = \int_{X_P} \int_{\T }  \Pi_{Np}^{\CP^m}
(t\cdot y_1,  y_2) \Pi_{N }^{M_P}( x, t \cdot x)\, d\vol_{\T}(t)\,
d\vol_{X_P}(x).$$
\end{lem}

\begin{proof} We consider the lattice sums
 \begin{equation}\label{lattice}
\chi_{NP}(t): = \sum_{\alpha \in NP}
t^\al\;,\quad t\in \T\;. \end{equation}
We easily see that
\begin{equation}\label{lattice1}\Pi_{|N P} (y_1,y_2) =  \int_{\T }
\Pi^{\CP^m}_{N\de} (t\cdot y_1,y_2) {\chi_{NP}(\bar
t)}\,d\vol_{\T}(t)\;.\end{equation} Recalling (\ref{Sztor}), we also have
\begin{equation}\label{lattice2}\int_{X_P}\Pi_N^{M_P}(x,t\cdot x)\,d\vol_{X_P}(x)
=\int_{X_P}
\sum_{\al\in NP}
\frac{\bar t^\al}{\|\wh \chi^P_\al\|_{X_P}^2}|\wh \chi^P_\al(x)|^2
\,d\vol_{X_P}(x)=\chi_{NP}(\bar t)\;. 
\end{equation}
The conclusion follows immediately from (\ref{lattice1})--(\ref{lattice2}).
\end{proof}

\bigskip
By Lemmas \ref{aha} and
\ref{PIP}, we have  for $z\in\C^{*m}$,
\begin{eqnarray} \Pi_{|NP}(z,z) &=& \frac{1}{(2\pi )^m} \int_{M_P} \int_{\T
}
\Pi_{Np}^{\CP^m}(e^{i\phi} \cdot z,z) \Pi_{N}^{M_P}( w,e^{i\phi} \cdot w)\,
d\phi\, d\vol_{M_P}(w)\nonumber \\&= & K_N(z)+ S_N(z)
\;,\label{szego1}\end{eqnarray} where
\begin{eqnarray}K_N(z) &=&
\frac{1}{(2\pi )^m} \int_{M_P}
\int_{\T }
\Pi_{Np}^{\CP^m}(e^{i\phi} \cdot z,z)\sigma_N(w) \wt\Pi_1^{M_N}(w, e^{i\phi} \cdot
w) \,
d\phi\, dw\;,\label{szego2}\\
S_N(z)&=&\frac{1}{(2\pi )^m} \int_{M_P} \int_{M_P}
\int_{\T }
\Pi_{Np}^{\CP^m}(e^{i\phi} \cdot z,z)R_N(w,\eta) \wt\Pi_1^{M_N}(\eta, e^{i\phi}
\cdot w) \,
d\phi\, dw\, d\eta\;,
\label{szego-rem}
\end{eqnarray}
where $\sigma_N,\ R_N$ are  as in Lemma \ref{aha}, and $\wt\Pi_1^{M_N}$ is given by
(\ref{Pi1}). (The above  are really functions of
$x=(z,\theta)$, but since they depend only on $z$, we
dropped the variable
$\theta$. In (\ref{szego-rem}), we replaced integrals over  $X_P$ with integrals
over $M_P$ by setting
$R_N(w,\eta)= \int_0^{2\pi} \int_0^{2\pi}R_N((w,\theta_1),(\eta,\theta_2))
e^{iN(\theta_1-\theta_2)}\frac{d\theta_1}{2\pi} \frac{d\theta_2}{2\pi}\,$.)

By (\ref{msum}), we have
$$\Pi_{Np}^{\CP^m}(x,y)=\frac{(Np+m)!}{(Np)!}\langle x,y\rangle^{pN}=
\frac{(Np+m)!}{(Np)!} \left(\sum_{|\al|\le p} \wh m_\al^{p\Si}(x) \overline{\wh m_\al^{p\Si}(y)}\right)^N\;.$$
Recalling (\ref{Pi1}),
we can then rewrite (\ref{szego2}):
\begin{equation} K_N(z)= \frac{1}{(2\pi )^m}\int_{M_P} \int_{\T }
e^{N\Psi(\phi,w; z)} a_N(w)\, d\phi\,
d\vol_{M_P}(w)\;,
\label{szego3}
\end{equation}
where the phase
is given by
\begin{equation} \label{PHASE2} \Psi(\phi,w; z) =
\log \sum_{\alpha \in P} e^{-i \langle \alpha, \phi
\rangle} |\wh m_\al^P(w)|^2 + \log  \sum_{|\al| \le\de} e^{i \langle
\alpha, \phi \rangle} |\wh m_\al^{p\Si}  (z)|^2\;, \end{equation}
and  $a_N(w)=\frac{(Np+m)!}{(Np)!}\sigma_N(w)$ is a symbol of order $2m$.

We fix $z\in\C^{*m}$.  It follows from the triangle inequality and (\ref{sum=1}) that
\begin{eqnarray}\Re\Psi(\phi,w;z)&=& \textstyle \log \left|\sum_{\alpha \in
P} e^{-i
\langle\alpha,\phi
\rangle} |\wh m_\al^P(w)|^2\right| + \log \left| \sum_{|\al| \le\de}
e^{i \langle \alpha, \phi \rangle} |\wt{m}_{\alpha}^{p\Sigma}  (z)|^2\right|
\nonumber\\ &\le & \Re\Psi(0,w;z) =0\;,\label{triangleineq}\end{eqnarray} with equality
if and only if
$\langle \alpha, \phi \rangle=\langle 0, \phi \rangle =0 $ for all $\al\in
\de \Si$, or equivalently, if $\phi=0$.  This says that $\Psi$ is a phase
function of positive type, and $\Re\Psi(\phi,w)=0$ if and only if $\phi=0$.
Hence we
need consider only critical points with
$\phi=0$.

Suppose that
$w=e^{i\theta}e^\rho\in\C^{*m}$.  We note that
$\Psi$ is independent of $\theta$; hence the critical set is $\T
$-invariant. So we
must find the set of
$w$ where   $$D_\phi \Psi(0,w)=0,\qquad D_\rho \Psi(0,w)=0\;.$$

We have
\begin{equation}\label{dphi}
iD_\phi\Psi=  \frac{\sum_{\alpha \in P} e^{-i \langle
\alpha, \phi
\rangle} |\wh m_\al^P(w)|^2 \alpha}{\sum_{\alpha \in P} e^{-i
\langle \alpha, \phi \rangle} |\wh m_\al^P(w)|^2} -
\frac{\sum_{\alpha \in p \Sigma} e^{i \langle \alpha, \phi
\rangle} |\wh m_\al^{p\Si}(z)|^2 \alpha}{\sum_{\alpha \in p \Sigma}
e^{i \langle \alpha, \phi \rangle} |\wh m_\al^{p\Si}(z)|^2}\;.\end{equation}
Hence the $D_\phi$ equation gives
\begin{equation} \frac{1}{p}\,
\mu_P(w)=\mu_\Si(z)\;.\label{dphi2}\end{equation}

Note that $$
D_\rho\Psi(0,w) = D_\rho 0
\equiv  0$$ automatically by (\ref{PHASE}).
However, we will need to compute $D_\rho\Psi$. To do this, we note that
$$\sum_{\alpha
\in P} e^{-i
\langle
\alpha,
\phi\rangle} |\wh m_\al^P(w)|^2 =\frac{\sum_{\alpha \in P} e^{-i
\langle
\alpha,
\phi\rangle} |{m}_{\alpha}(w)|^2}{\sum_{\alpha \in P}
|{m}_{\alpha}(w)|^2}\;.$$
Recalling that $w=e^{i\theta}e^\rho$, we then obtain
\begin{equation}\label{drho}\half D_\rho\Psi= \frac{\sum_{\alpha \in P}
e^{-i
\langle
\alpha,
\phi\rangle}|\wh m_\al^P(w)|^2 \alpha}{\sum_{\alpha \in P} e^{-i
\langle
\alpha,
\phi\rangle}|\wh m_\al^P(w)|^2 } - \mu_P(w)\;.\end{equation}

Let $\ccal_\Psi^0\subset \T \times M_P$ denote the set of critical points of
$\Psi$
where $\Re\Psi=0$. By (\ref{triangleineq}) and (\ref{dphi2}), we have:
\begin{equation} {\mathcal C}_{\Psi}^0 =\{ \phi = 0, w \in
\mu_P\inv(\mu_\Si(z))\}
\simeq \T .
\label{crit0}\end{equation}
Equation (\ref{crit0}) leads to the important conclusion:

\medskip
{\it There are no critical points with real part equal to 0 unless $z \in
(\mu_\Si)^{-1} (\frac{1}{p}P)$.}

\medskip\noindent
As we stated in the introduction, we refer to the region $\acal_P= (\mu_\Si)^{-1}
(\frac{1}{p}P^\circ)$ as the {\it classically allowed} region, and its complement
in $\C^{*m}$ as the {\it classically forbidden} region.  We now break up our
problem into these two
regions.

\subsection{The classically allowed region}
In this section, we illustrate our approach by showing that
\begin{equation}\label{pre-mass}\Pi_{|NP}(z,z) \sim N^m \big[c_0(z) +
c_1(z)N^{-1} + c_2(z)N^{-2} + \cdots \big],  \quad \mbox{for }\ z \in \acal_P
\;.\end{equation}
In \S \ref{s-forbidden}, we will derive our exponentially decaying, asymptotic
expansion of
$\Pi_{|NP}(z,z)$ on the forbidden region, which we will  use to obtain the precise
formula of Proposition~\ref{SZEGO}(i) on the allowed region.  (Although we do not
need  (\ref{pre-mass}) in the derivation of our precise formula in  \S
\ref{s-forbidden}, we shall use the methods and Lemma~\ref{A+} from this section in
our proof.) 

Suppose that $z \in \acal_P=(\mu_\Si)^{-1} (\frac{1}{p}P^\circ )$ is in the
classically allowed  region.   Since the phase function is
of
positive type and equals zero on the critical set, in order to apply the
method of
stationary phase it suffices to show that the Hessian of $\Psi(\phi,\rho)$
is
nondegenerate on the critical set.  We write the Hessian as a block matrix:
$$D^2\Psi|_{(0,\rho)}=\left(\begin{array}{cc}A&B\\B^t&C\end{array}\right)\;,
$$
where
$$A_{jk}=\frac{\d^2 \Psi}{\d\phi_j\d\phi_k}(0,\rho)\;,\quad
B_{jk}=\frac{\d^2 \Psi}{\d\phi_j\d\rho_k}(0,\rho)\;,\quad
C_{jk}=\frac{\d^2 \Psi}{\d\rho_j\d\rho_k}(0,\rho)\;.$$
Differentiating (\ref{dphi}) and recalling that $\mu_P(w)=\sum_{\alpha \in
P}|\wh m^P_\al(w)|^2 \al$, we obtain
$$A = \mu_P(w)\otimes   \mu_P(w) -\sum_{\alpha \in P} |\wh m_\al^P(w)|^2 \alpha
\otimes  \alpha +
\mu_{\de\Si}(z)\otimes   \mu_{\de\Si}(z) -\sum_{|\alpha|\le \de} |\wh m_{\alpha}^{p\Si}(z)|^2 \alpha
\otimes  \alpha\;.$$ (For a vector $v=(v_1,\dots,v_m)\in\R^m$, we identify
$v\otimes
v$ with the symmetric matrix $(v_jv_k)_{1\le j,k\le m}$.)

Thus on $\ccal_\Psi^0$, we have by (\ref{dphi2}),
\begin{equation}\label{A} A = 2p^2\, \mu_\Si(z)\otimes   \mu_\Si(z)
-\sum_{\alpha \in
P}|\wh m_\al^P(w)|^2 \alpha \otimes  \alpha -\sum_{|\alpha|\le p} |\wh m_\al^{p\Si}(z)|^2 \alpha
\otimes  \alpha\;.\end{equation} On the other hand, differentiating
(\ref{drho}), we
obtain
\begin{equation}\label{B} B= -2i\left( \sum_{\alpha
\in P} |\wh m_\al^P(w)|^2 \al\otimes  \al - \mu_P(w)\otimes
\mu_P(w)\right)\;.\end{equation}  Finally, since $\Psi(0,\rho)\equiv 0$, we
have
$$C=0\;.$$ Therefore
\begin{equation}\label{Hessdet}\det D^2\Psi|_{(0,\rho)}= (\det B)^2
\;.\end{equation}

In order to apply the method of stationary phase to
conclude that  the integral (\ref{szego3}) admits a complete asymptotic
expansion
in powers of $N^{-1}$, we must show that $ \det D^2\Psi|_{(0,\rho)}\ne 0$.
This
follows from (\ref{Hessdet}) and the following fact:

\begin{lem} \label{A+}The real symmetric matrix
$$S_w:= \sum_{\alpha
\in P} |\wh m_\al^P(w)|^2 \al\otimes  \al - \mu_P(w)\otimes
\mu_P(w)$$ is strictly positive definite for all $w\in \C^{*m}$.
Furthermore, if
$\mu_P(w)$ lies in a face $F$ of $P$, then $S_w$ is semi-positive and its $0$
eigenspace is $T_F^\perp$; in particular, $S_w$ is positive definite on
$T_F$.
\end{lem}

\begin{proof} We must show that $(\la\otimes\la,S_w) >0$ for $\la\in
(\R^m)'\sm\{0\}$. Consider the vectors $u,v\in\R^{\#P}$  given by $u_\al=|\wh m_\al^P(w)|$,
$v_\al=|\wh m_\al^P(w)|\la(\al)$. Since $\|u\|^2=\sum _{\al\in P} |\wh m_\al^P(w)|^2=1$
and $\mu_P(w)=\sum_{\al\in P}|\wh m_\al^P(w)|^2\al$, we have
\begin{eqnarray}\label{S}(\la\otimes\la,S_w) &=& \sum_{\al\in P} |\wh m_\al^P(w)|^2
\la(\al)^2 - \left(\sum
_{\al\in P} |\wh m_\al^P(w)|^2 \la(\al)\right)^2\nonumber \\&=& \|v\|^2- \langle u,
v\rangle^2 \ =\ \|u\|^2\|v\|^2- \langle u,
v\rangle^2\ \ge \ 0\;.\end{eqnarray}

If $w\in\C^{*m}$, then the Cauchy-Schwartz inequality in (\ref{S}) is strict
since $|\wh m_\al^P(w)|\ne 0$ for all
$\al\in P\cap\Z$ and $\la(\al)$ is not constant on $P\cap\Z$. On the other hand, if
$\mu_P(w)$ lies in a face $F$, then by Proposition~\ref{toricdivisors},
$|\wh m_\al^P(w)|\ne 0
\Leftrightarrow
\al\in \bar F$, and hence  $$(\la\otimes\la,S_w) = 0 \Leftrightarrow
\la|_F\equiv \,\mbox{const.} \Leftrightarrow \la\perp
T_F\;.$$
\end{proof}

The asymptotic expansion (\ref{pre-mass}) for $z$ in the
classically allowed region $\acal_P$ now follows from the
method of stationary phase (\cite[Theorem~7.7.5]{H})  applied to the complex oscillatory integral (\ref{szego3}) and the fact that $S_N(z)$ is rapidly decaying.
To be precise, the asymptotics of  (\ref{szego3}) are determined by
 the  component of the critical point set of the phase  on which $\Re \Psi$ is maximal.
Note that $\Re \Psi$ and $\Im \Psi$ are constant
on components of the critical set, and that on its maximal component
our phase satisfies $\Re \Psi = 0$. The critical point set is  a manifold
and we have just shown that the  phase is non-degenerate in the normal directions.
Although  (\cite[Theorem~7.7.5]{H})
assumes the critical point with
$\Re
\Psi = 0$ is isolated, the proof in \cite{H} can be modified, precisely as in the real
case, to apply to phases with non-degenerate critical manifolds. 

Hence the expansion follows from the complex stationary phase method for phases
of positive type with non-degenerate critical manifolds.  Since the critical manifold is
of (real) codimension $2m$
and the amplitude $a_N$ in (\ref{szego3}) is or order $2m$, it follows that the leading
term of the expansion contains $N^{2m-m}=N^m$, as claimed.

\subsection{The classically forbidden region} \label{s-forbidden}

For $z$ in  this region of $\C^{*m}$, there are no critical points of the
phase with real part 0, and the integral is rapidly
decaying by (\ref{triangleineq}). Perhaps surprisingly, this implies that $\partial
\bar{\partial}
\log \Pi_{|NP}(z,z)$ {\it will} contribute here to the distribution of
zeros.

To determine the exponential decay rate of the integral (\ref{szego3}), we
shall deform the contour of integration to pick up critical points with
maximal real
part along the contour. To accomplish this, we first complexify the real
torus $\T $
to
$\C^{*m}$ with variables
$\zeta_j =  e^{\tau_j + i \phi_j}$. The integral (\ref{szego3}) may be
written in
terms of these variables as
\begin{eqnarray} K_N(z) &=&\frac{1}{(2\pi i)^m} \int_{M_P} \int_{\T
}
e^{N\Psi_\C(\zeta,w; z)} \sigma_N(w)\, \prod_{j = 1}^m \frac{d
\zeta_j}{\zeta_j}\; d\vol_{M_P}(w),\label{PIP2}\end{eqnarray}
where
\begin{eqnarray}\Psi_\C(\zeta,w; z) &=&
\log \sum_{\alpha \in P} \zeta^{-\al}|\wh m_\al^P(w)|^2 + \log
\sum_{|\al| \le\de} \zeta^\al |\wh m_\al^{p\Si}
(z)|^2\nonumber \\ &=& \log \sum_{\alpha \in P} e^{- \langle \alpha, \tau +i\phi
\rangle} |\wh m_\al^P(w)|^2 + \log  \sum_{|\al| \le\de} e^{ \langle
\alpha, \tau +i\phi \rangle} |\wh m_{\alpha}^{\de
\Sigma}  (z)|^2\;.  \label{PHASEcx}\end{eqnarray}
Since the integrand is holomorphic in $\zeta\in
\C^{*m}$, we  can deform
the contours in (\ref{PIP2}) and instead integrate over
\begin{equation}\label{contour}\int_{M_P}
\int_{|\zeta_1| = \tau_1, \cdots,  |\zeta_m| = \tau_m} \;.\end{equation}
We therefore look for  complex critical points of $\Psi_\C$.
We note that as before,
\begin{equation}\label{maxatphi0}\Re
\Psi_\C(e^{\tau +i\phi},w) <\Re\Psi_\C(e^{\tau},w)\;,\qquad \mbox{for\ \ }
\phi\ne
0\;.\end{equation} Thus the  critical points $\Psi$ with maximal real values on
the contour
(\ref{contour}) are those with
$\phi=0$.

We now write $\Psi_\C(\tau,\phi,w)=\Psi_\C(e^{\tau+i\phi},w;z)$.
The
complexification of (\ref{dphi}) becomes
\begin{eqnarray}
iD_\phi\Psi_\C &=&  \frac{\sum_{\alpha \in P}
e^{- \langle\alpha, \tau+i\phi\rangle} |\wh m_\al^P(w)|^2
\alpha}{\sum_{\alpha
\in P} e^{-\langle\alpha, \tau+i\phi\rangle} |\wh m_\al^P(w)|^2} -
\frac{\sum_{\alpha \in p \Sigma} e^{\langle\alpha, \tau+i\phi\rangle} |\wh m_\al^{p\Si}(z)|^2 \alpha}{\sum_{\alpha \in p\Sigma} e^{\langle\alpha,
\tau+i\phi\rangle} |\wh m_\al^{p\Si}(z)|^2}\nonumber\\
&=&  \frac{\sum_{\alpha \in P}
e^{- \langle\alpha, \tau+i\phi\rangle} | m_{\alpha}(w)|^2
\alpha}{\sum_{\alpha
\in P} e^{-\langle\alpha, \tau+i\phi\rangle} | m_{\alpha}(w)|^2} -
\frac{\sum_{\alpha \in p \Sigma} e^{\langle\alpha, \tau+i\phi\rangle} |
m_{\alpha}(z)|^2 \alpha}{\sum_{\alpha \in p \Sigma} e^{\langle\alpha,
\tau+i\phi\rangle} |
m_{\alpha}(z)|^2}\nonumber\\
&=&  \frac{\sum_{\alpha \in P}
e^{- i\langle\alpha, \phi\rangle} |\wh m_\al^P(e^{-\tau/2}\cdot w)|^2
\alpha}{\sum_{\alpha
\in P} e^{- i\langle\alpha, \phi\rangle} |\wh m_\al^P(e^{-\tau/2}\cdot
w)|^2
} -
\frac{\sum_{\alpha \in p \Sigma} e^{i\langle\alpha, \phi\rangle} |\wh m_\al^{p\Si}(e^{\tau/2}\cdot z)|^2 \alpha}{\sum_{\alpha \in p \Sigma}
e^{i\langle\alpha, \phi\rangle} |\wh m_\al^{p\Si}(e^{\tau/2}\cdot z)|^2}\;,\label{dphicx0}\end{eqnarray}
for $w\in\C^{*m}$.  Hence by continuity, we see that
\begin{equation} iD_\phi\Psi_\C =  \frac{\sum_{\alpha \in P}
e^{- i\langle\alpha, \phi\rangle} |\wh m_\al^P(e^{-\tau/2}\cdot w)|^2
\alpha}{\sum_{\alpha
\in P} e^{- i\langle\alpha, \phi\rangle} |\wh m_\al^P(e^{-\tau/2}\cdot
w)|^2
} -
\frac{\sum_{\alpha \in p \Sigma} e^{i\langle\alpha, \phi\rangle} |\wh m_\al^{p\Si}(e^{\tau/2}\cdot z)|^2 \alpha}{\sum_{\alpha \in p \Sigma}
e^{i\langle\alpha, \phi\rangle} |\wh m_\al^{p\Si}(e^{\tau/2}\cdot z)|^2}\;,\label{dphicx}\end{equation} for all
$w\in
M_P$.  In particular,
\begin{equation} iD_\phi\Psi_\C (\tau,0,w)=
\mu_P(e^{-\tau/2}\cdot w)-\mu_{\de \Si}(e^{\tau/2}\cdot
z)\;.\label{ii}\end{equation}

 Hence the equation $D_\phi\Psi_\C(\tau,0,w)=0$ is
equivalent to
\begin{equation} \frac{1}{p}\,
\mu_P(e^{-\tau/2}\cdot w)=\mu_\Si(e^{\tau/2}\cdot
z)\;.\label{ii''}\end{equation}

To compute the $w$ derivative, we first consider the case where
$w=e^{\rho+i\theta}\in \C^{*m}$. Recalling that $$\sum_{\alpha
\in P} e^{-
\langle
\alpha,\tau+i
\phi\rangle} |\wh m_\al^P(w)|^2 =\frac{\sum_{\alpha \in P} e^{-i
\langle
\alpha,
\phi\rangle} |{m}_{\alpha}(e^{-\tau/2}\cdot w)|^2}{\sum_{\alpha \in P}
|{m}_{\alpha}(w)|^2}
 =\frac{\sum_{\alpha \in P} e^{-i \langle
\alpha,
\phi\rangle} |{m}_{\alpha}(e^{-\tau/2+\rho})|^2}{\sum_{\alpha \in P}
|{m}_{\alpha}(e^\rho)|^2}\;,$$ we obtain
\begin{equation}\label{drhocx}\half D_\rho\Psi_\C= \frac{\sum_{\alpha \in P}
e^{-i
\langle
\alpha,
\phi\rangle}|\wh m_\al^P(e^{-\tau/2}\cdot w)|^2 \alpha}{\sum_{\alpha \in
P}
e^{-i
\langle
\alpha,
\phi\rangle}|\wh m_\al^P(e^{-\tau/2}\cdot w)|^2 } -
\mu_P(w)\;.\end{equation}
Setting $\phi=0$,
\begin{equation}\label{GRADIENT} D_\rho\Psi_\C(\tau,0,w)= \mu_P(e^{-\tau/2}
w)
- \mu_P(w)\;.\end{equation} Hence the equation $D_\rho\Psi_\C|_{\phi=0}=0$
gives
\begin{equation} \mu_P(e^{-\tau/2}\cdot w) =
\mu_P(w)\;.\label{iii''}\end{equation}

We must show that (\ref{iii''}) also holds for critical points $w\not\in
\C^{*m}$. (These are points with $\mu_P(w)\in\d P$.)  To do this we
introduce
the function
\begin{equation}\label{HP}H_P(\tau,w):=\frac{\sum_{\alpha \in P} e^{-
\langle
\alpha, \tau \rangle}|m_\al(w)|^2}{ \sum_{\alpha \in P}
|m_\al(w)|^2} =\sum_{\al\in P}e^{-\langle \al,\tau\rangle}|\wh m^P_\al(w)|^2\;.\end{equation} Note that $H_P(\tau,\cdot)$ is a
well-defined function on $M_P$ that factors through the moment
map $\mu_P$. We also have
\begin{equation}\label{Hd} H_{\de\Si}(\tau, z) = H_{\Si}(\tau, z)^\de
=\left[\frac{1+\|e^{-\tau/2}\cdot z\|^2}{1+\|
z\|^2}\right]^\de\;.
\end{equation}

We have
\begin{equation}\label{phaseH}\Psi_\C(\tau,0,w;z)=\log
H_P(\tau,w)+\log H_{\de\Si}(-\tau,z)\;,\end{equation} and hence
$$D_w\Psi_\C(\tau,0,w;z)=D_w\log H_P(\tau,w)\;.$$
Furthermore the vector fields \begin{equation}\label{Xj}\textstyle X_j: =
\frac{\partial}{\partial \rho_j}=\Re( w_j \frac{\partial}{\partial
w_j})\end{equation} extend to
$M_P$ as the generators of the $\R_+^m$ action (induced by the $\C^{*m}$
action). Hence (\ref{GRADIENT}) is equivalent to the equation
\begin{equation} \label{INVGRAD} (X_1, \dots, X_m) \log H_P(\tau,
w) = \mu_P(e^{-\tau/2} w) - \mu_P(w), \end{equation} which is
invariantly defined on all of $M_P$.

Thus we have seen that critical
points of $H_P(\tau, \cdot)$ must satisfy  (\ref{iii''}).  We shall show later
that  (\ref{iii''}) is also a sufficient condition for $w$ to be a critical
point of $H_P(\tau, \cdot)$, and hence (\ref {ii''}) and (\ref{iii''}) give
necessary and sufficient conditions for $(0,w)$ to be a critical point of
$\Psi_\C(\tau,\cdot,\cdot)$.

The following fact is useful in describing the
critical set.

\begin{lem}\label{sta} Let $\tau\in \R^m$, $w\in M_P$.  Then the following are
equivalent:
\begin{enumerate}
\item[\rm a)] \ $\mu_P(e^\tau \cdot w) = \mu_P(w)$;
\item[\rm b)] \  $e^{r\tau} \cdot w = w$, for all $r\in\R$;
\item[\rm c)] \ $\tau$ is perpendicular to $T_{F^w}$, where $F^w$ is the
face of $P$ containing $\mu_P(w)$.
\end{enumerate}\end{lem}

\begin{proof} Let $w=(w_\al)_{\al\in P}\in
M_P\subset\CP^{\#P-1}$. By Proposition \ref{toricdivisors}, we see that $w_\al=0$ if and
only if $\al\not\in \overline{F^w}$, and thus for any
$t\in\C^{*m}$, \begin{equation}\label{stabiliz}
t\cdot w = w\ \Leftrightarrow\ \exists c \mbox{\
\ such that\ \ }t^\al w_\al = c w_\al\ \forall \al\in P \ \Leftrightarrow\
t^\al  = c
\ \forall
\al\in \overline{F^w}\cap\Z^m\;.\end{equation} Suppose (a) holds. Then there
exists $e^{i\phi}\in \T$ such that
$e^{\tau+i\phi}\cdot w=w$. Hence by (\ref{stabiliz}), $e^{\langle
\tau+i\phi,
\al\rangle} =c$ and thus  $e^{\langle
r\tau,
\al\rangle} =|c|^r$, for all $\al\in \overline{F^w}\cap\Z^m$. Therefore (b)
holds.

Furthermore, it follows from (\ref{stabiliz}) that (b) is equivalent to
\begin{equation}\label{c}\langle
\tau, \al-\al'\rangle = 0 \qquad\forall\ \al ,\al'\in
\overline{F^w}\cap\Z^m\;,\end{equation} which is a restatement of (c).
\end{proof}

\begin{rem} Let $\tau\in
\R^m$. It follows from Lemma \ref{sta}
that the set $A_\tau$ of points in $M_P$ satisfying (\ref{iii''}) consists of
the inverse image  under $\mu_P$ of the union of faces of P that are orthogonal
to
$\tau$.  The set $A_\tau$ is a (possibly empty) smooth algebraic subvariety of
$M_P$. It then follows from (\ref{HP}) and (\ref{c}) that
$H_P(\tau,\cdot)$ is constant on each connected component of $A_\tau$.  Indeed,
on a component $A^0$ of  $A_\tau$, we have $H_P(\tau,\cdot)\equiv e^{-\langle
\al^0, \tau\rangle}$, where $\al^0$ is a lattice point in $\mu_P(A^0)$.
\end{rem}

\subsubsection{The normal bundle}\label{s-normalbundle}  Although equations
(\ref{ii''}) and (\ref{iii''}) are necessary and sufficient conditions for
$(0,w)$ to be a critical point
of the phase on the $\tau$-contour (as we shall show), they do not guarantee
that the real part of the phase is maximal at
that point.  To obtain maximality, we need to strengthen condition
(\ref{iii''})
 to require that $-\tau$ lies in the  `normal bundle' of $P$ at
$w_z$.

Recall that the normal cone $C_F$ of a face $F$ of $P$ is given by
\begin{equation}\label{ncone}C_F =\{u\in\R^m:\langle u,
x\rangle =\sup _{y\in P}  \langle u,y\rangle,\ \forall x\in F\}\;.\end{equation}
We now consider compact polytopes in $\R^m$ with arbitrary vertices that are not
necessarily integral or even rational.  We define the faces of such polytopes as in
the integral case, and we use  (\ref{ncone}) to define the normal cones.

\begin{defin} The {\it normal bundle\/}
$\ncal(Q)$ of a (not necessarily integral) convex polytope $Q\subset\R^m$ is the
subset of
$T_{\R^m}=\R^m\times \R^m$ consisting of pairs
$(x,v)$, where
$x\in Q$ and $v$ is in the normal cone $C_F$
of the face $F$ of $Q$ that contains $x$. (Note that $\ncal(Q)$
is not a fiber bundle over $Q$.)
\end{defin}

The normal bundle
$\ncal(Q)$ is a piecewise smooth submanifold of
$T_{\R^m}$; it is homeomorphic to $\R^m$ via the `exponential map'
$$\ecal_Q:\ncal(Q)\to
\R^m\;, \quad \ecal_Q(x,v)= x+v,\qquad x\in F,\ c\in C_F\qquad (F\
\mbox{a face of} \ Q)\;.$$ It is easily seen that $\ecal_Q$ is a
homeomorphism and is a
$\ccal^\infty$ (in fact,
linear) diffeomorphism on each of the `pieces' $\bar F\times C_F\subset
\ncal(Q)$.

We shall use the following elementary, but not so well known, fact about the
invertibility of Lipschitz maps. We let $B_\ep(x_0)=\{x\in\R^m:|x-x_0|<\ep\}$
denote the $\ep$-ball about $x_0\in\R^m$.

\begin{lem}\label{lipeo} {\rm \cite{Fan}}
Let $f:U\to\R^m$ be a Lipschitz map, where $U$ is open in $\R^m$.  Then $f$
has a
local orientation-preserving Lipschitz inverse with Lipschitz constant $L$ at a
point
$x_0\in U$ if and only if there exists $\ep>0$  such that
\begin{enumerate}
\item[i)] \  $\liminf_{v\to 0} |f(x+v) - f(x)|/|v|  \ge 1/L$ for all $x \in
B_\ep(x_0)$,
\item[ii)] \ $\det f'(x)>0$ for all $x \in
B_\ep(x_0)$ such that $f$ is differentiable at $x$,
\item[iii)] \ $f(x_0)\not\in f(\d B_\ep(x_0))$ and $\deg\big[f:\d B_\ep(x_0) \to
\R^m\sm\{f(x_0)\}\big] = 1$.
\end{enumerate}\end{lem}

\noindent By the degree in (iii), we mean the degree of
$f_*(\d B_\ep(x_0)) \in H_{m-1}(
\R^m\sm\{f(x_0)\},\Z)$.  The hypotheses (i) and (iii) in \cite{Fan} differ
slightly from those above, but Fan's proof of sufficiency uses only (i)--(iii).
(Necessity is obvious.) Hypothesis (i) is given as a lemma in
\cite{Fan}. Hypothesis (iii) above is
replaced in \cite{Fan} by the condition that the index of $f$ at $x_0$ is 1,
which is easily seen to be equivalent to (iii) under the assumptions (i) and
(ii).

We now let
$$\ncal^\circ=\textstyle\ncal(\frac{1}{\de}P)
\cap(\Si^\circ\times\R^m)\;.$$

\begin{lem} For each point $z\in \C^{*m}$, there exists a
unique
$(x,v)\in
\ncal^\circ$ such that $$\mu_\Si(e^{v/2}\cdot
z)=x\;.$$\label{claim}\end{lem}

\begin{proof} We note that the $\R_+^m$ action on $\C^{*m}$ descends via
the moment
map
$\mu_\Si$ to an $\R_+^m$ action on $\Si^\circ$, and we  consider the map
$$\Phi:\ncal^\circ \to \Si^\circ\;,\qquad \Phi(x,v)=
\mu_\Si(e^{v/2}\cdot
x)\;.$$ It suffices to show that $\Phi$ is a bijection.

Let $\lcal:\Si^\circ \rightarrow\R^m$ be the diffeomorphism given by
\begin{equation}\label{lcal}\lcal(x) =\left(\log \frac{x_1}{1-\sum x_j},\ldots,
\log
\frac{x_n}{1-\sum x_j}\right)\;,\end{equation} so that
$$\lcal \circ \mu_\Sigma (z) = (\log |z_1|^2, \ldots, \log |z_m|^2)\;.$$ Thus,
writing $x = \mu_\Sigma(z)$ we have
\begin{eqnarray}\lcal \circ \Phi (x,v) &=& \lcal \circ\mu_\Sigma (e^{v/2}\cdot
z)\nonumber\\ &=&(v_1 +\log |z_1|^2,\ldots, v_m +\log |z_m|^2)\nonumber\\
&=&v+\lcal (x)\;.\label{lcalphi}\end{eqnarray} We first observe that $\lcal
\circ \Phi$ is proper:  suppose on the contrary that the sequence $\{(x^n,
v^n)\}$ is unbounded in $\ncal^\circ$, but $\lcal \circ\Phi (x^n,
v^n)\rightarrow a\in \R^m$.  By passing to a subsequence, we can assume that
$x^n\rightarrow x^0 \in
\frac{1}{p}P$.  Then $x^0 \in \partial \Si$, since otherwise
$v^n\rightarrow a-\lcal (x^0)$.  Write $v^n = r^n u^n$, where $r^n>0$,
$|u^n|=1$. We can assume without loss of generality that $u^n\rightarrow u^0$.
We consider the case where
$\Si x^0_j<1$, i.e., $\mu_\Sigma^{-1}(x^0)$ is not in the hyperplane at
infinity.  If $x^0=0$, then $u^0_j\geq 0$ $(1\leq j\leq m)$, since otherwise
$\lcal \circ\Phi (x^n, v^n)$ would diverge.  But $u^0\in C_{\{0\}}$; hence
$\frac{1}{p} P\subset \{x:\langle u^0, x\rangle \leq 0\}\subset
\R^m\backslash \Si^\circ$, a contradiction. Now suppose that
$x^0_1=\cdots=x^0_k=0$,
$x_l>0$ for $k<l\leq m$.  Then we conclude as before that $u^0_j\geq 0$
for
$1\leq j\leq k$, and
$u^0_l = 0$ for $k<l\leq m$, and we again obtain a contradiction.
Finally, if $\Si x^0_j=1$, we can change coordinates (permute the homogeneous
coordinates in $\C\PP^m$) to reduce to the previous case.

Let $\ecal=\ecal_{\frac{1}{\de}P}:  \ncal(\frac{1}{\de}P)
\buildrel{\approx}\over\to \R^m$, and let $U=\ecal (\ncal^\circ)\subset
\R^m$. We consider the map
$f:U\rightarrow \R^m$ given by $f\circ \ecal|_{\ncal^{\circ}}=\lcal \circ
\Phi$, i.e. by the commutative diagram:

$$\begin{array}{ccc} \ncal^\circ  &\buildrel{\Phi} \over \longrightarrow
& \Sigma^\circ \\ {\scriptstyle\approx}\!\downarrow {\scriptstyle\ecal} &
& {\scriptstyle\approx}\!\downarrow {\scriptstyle\lcal}  \\ U
&\buildrel{f} \over \longrightarrow & \R^m \end{array}$$
Since $\lcal\circ\Phi$ is a proper map,
$f$ is also proper.  Hence to show that $\Phi$ is a bijection, it suffices to
show that
$f$ is a local homeomorphism and is therefore a (global) homeomorphism.

To describe the map $f$, for each $x\in \Si^\circ$, we let $q_x$ denote the
closest point in $\frac{1}{p} P$ to $x$.  We note that $q_x\in \Si^\circ$;
if $x\not\in\frac{1}{p}P$, then $q_x\in \partial (\frac{1}{p}P)$; if
$x\in \frac{1}{p}P$, then $q_x=x$.  Furthermore, $\ecal^{-1}(x) = (q_x,
x-q_x)$ and hence
$$f(x) = x-q_x +\lcal (q_x)\;.$$

We shall show that $f$ satisfies the hypotheses of Lemma \ref{lipeo}. Let $F$
be an arbitrary face of
$\frac{1}{p}P\cap \Si^\circ$. To verify
(i) and (ii), it suffices to show that
$\det Df>0$ on the (noncompact) polyhedron $$U_F:=\ecal(\bar F\times C_F)\cap
U=\ecal[(\bar F\cap
\Si^\circ)\times C_F]\;.$$
To compute the determinant, we let
$T_F\subset
\R^m$,
$N_F=T_F^\perp\subset \R^m$ denote the tangent  and normal spaces,
respectively,  of $F$.  Let $x^0\in F$ be fixed.  For $y\in U_F$, we have
$$y=x + v \ \buildrel{f}\over \mapsto \ \lcal(x)+v\;, \qquad x-x^0\in T_F,\
v\in N_F\;.$$  Choose orthonormal bases $\{Y_1,\dots, Y_r\}$, $\{Y_{r+1},\dots,
Y_m\}$ of $T_F,\ N_F$, respectively. We let $Df$ denote the matrix of
the derivative
$(f|_{U_F})'$ with respect to the basis
$\{Y_1,\dots,Y_m\}$ of $\R^m$.  We have:
\begin{equation}\label{Df}Df=\left(\begin{array}{cc}
T^t\lcal'(x)T&0\\ *&I\end{array}\right)\;,
\end{equation} where $T$ is the $m\times r$ matrix $[Y_1 \cdots Y_r]$.
We have by (\ref{lcal}),
\begin{equation}\label{lcal'}\Big(\lcal'(x)\Big)_{jk} = \frac{1}{x_0} +
\delta^k_j
\frac{1}{x_j}\;,\quad x_0=1-\sum_{j=1}^m x_j>0\;,\end{equation} for
$x\in\Si^\circ$. Hence
$\lcal'(x)$ is a positive definite symmetric matrix, it being the sum of a
semipositive matrix (all of whose entries are
$\frac{1}{x_0}$) and a positive definite diagonal matrix.  Therefore,
$T^t\lcal'(x)T$ is positive definite, and hence
$\det Df(x)=\det (T^t\lcal'(x)T) >0$, completing the proof that hypotheses (i)
and (ii) of Lemma \ref{lipeo} are satisfied.

Note that (\ref{lcal'}) implies that the eigenvalues of
$\lcal'(x)$ are $>1$ for $x\in\Si^\circ$ and hence by (\ref{Df}), $Df(x)$ is a
diagonalizable matrix whose eigenvalues are   real and $\ge 1$.

We verify (iii) by a homotopy argument:  Choose a point $q^0\in
(\frac{1}{\de}P)^\circ$. We contract $\frac{1}{\de}P$ to $q^0$; i.e., for
$0\le t\le 1$, we let
$$Q_t= (1-t)\{q^0\} + \frac {t}{\de}P\;,$$ so that $Q_0=\{q^0\},\ Q_1=
\frac {1}{\de}P$. For $0\le t<1$, $Q_t\subset \Si^\circ$, and hence we have a
map
$$\Phi_t:\ncal(Q_t)\to \Si^\circ\;.$$ For $0\le t <1$, we define $f_t:\R^m\to
\R^m$ by the commutative diagram:
$$\begin{array}{ccc} \ncal(Q_t)  &\buildrel{\textstyle\Phi_t} \over
\longrightarrow & \Sigma^\circ \\[10pt] \ {\scriptstyle\approx}\!\downarrow
\ecal_{Q_t} & & {\scriptstyle\approx}\!\downarrow {\lcal}
\\ \R^m  &\buildrel{\textstyle f_t} \over \longrightarrow & \R^m \end{array}$$
As before we have
\begin{equation}\label{ft}f_t(x) = x-q^t_x +\lcal (q^t_x)\;,\end{equation}
where $q^t_x$ is the (unique) point of $Q_t$ closest to $x$. The above argument
shows that the maps $f_t$ also satisfy (i) and (ii) of Lemma \ref{lipeo}

We write
$$F:\big(\R^m\times [0,1)\big)\cup (U\times\{1\})\to \R^m\;,\qquad (x,t)\mapsto
f_t(x)\;,$$ where
$f_1=f:U\to\R^m$.  One easily sees
that $$\left|q^t_x -q^{t'}_{x'}\right| \le |x-x'|+|t-t'|\;,$$  and hence $F$ is
continuous. Furthermore, $F$ is uniformly continuous on $\R^m\times[0,t]$,
for each
$t<1$.

Let $H$ denote the set of
$t\in [0,1)$ such that  $f_t:\R^m\to\R^m$ is an orientation
preserving homeomorphism.   We note that
$0\in H$ since
$$f_0(x)=x - q+\lcal(q)= x + \mbox{constant}\;.$$
Note that  (\ref{Df}) also holds for the maps $f_t$ and hence as above,
the eigenvalues of
$Df_t$ are   real and $\ge 1$. Hence $f_t$ satisfies conditions (i) and (ii) of
Lemma \ref{lipeo} with $L=1$, for all $t\in [0,1]$.

First we show that $H$ is open in $[0,1)$.
Suppose $t_0\in H$. Then by Lemma \ref{lipeo},
Lip$(f_{t_0}\inv)\le 1$ and hence
$$|x-x_0|=1 \Longrightarrow |f_{t_0}(x)-f_{t_0}(x_0)|\ge 1\;.$$  By the uniform
continuity of
$F$ on $\R^m \times [0, \frac{t_0+1}{2}]$, we can choose $\ep>0$ such that for
all
$x_0\in\R^m$, we have
\begin{equation}\label{x}|f_t(x)-f_t(x_0)|\ge \half\;, \qquad \mbox{for} \  x\in
\d B_1(x_0),\ |t-t_0|<\ep\;.\end{equation}
To simplify notation, we shall write
$$\deg(f,x_0,r):=\deg\big[f:\d B_r(x_0) \to
\R^m\sm\{f(x_0)\}\big] \;.$$ We conclude from (\ref{x}) that for $|t-t_0|<\ep$,
$$\deg(f_{t},x_0,1) = \deg(f_{t_0},x_0,1)) = 1\;,$$ and hence
$f_t$ is a local homeomorphism at $x_0$,  by Lemma
\ref{lipeo}.  By the very first part of the
argument that $f:U\to \R^m$ is proper, we easily see that
$f_t$ is proper. Since
$x_0\in\R^m$ is arbitrary, it follows that
$f_t:\R^m\to \R^m$ is a covering map and therefore is a homeomorphism.

Next we show that $H$ is closed in $[0,1)$ and hence $H=[0,1)$.  Let $t_n\in
H$ such that
$t_n\to t_0\in [0,1)$. Since $f_{t_0}$ satisfies condition (i) of Lemma
\ref{lipeo}, we can choose $\ep>0$ so that $f_{t_0}(x_0)\not\in
f_{t_0}(\d B_\ep(x_0))$. Then for
$n$ sufficiently large,
$$\deg(f_{t_0},x_0,\ep) = \deg(f_{t_n},x_0,\ep) = 1\;.$$  It follows as
above that $f_{t_0}:\R^m\to \R^m$ is  a homeomorphism.

We have shown that $f_t$ is a homeomorphism for $0\le t <1$.
To complete the proof of the lemma, we must show that $f=f_1$ is a local
homeomorphism.  So we let $x_0\in U$ be arbitrary, and we choose $\ep>0$ such
that $\overline {B_\ep(x_0)}\subset U$ and $f(x_0)\not\in
f(\d B_\ep(x_0))$. Then for $t<1$ sufficiently close
to 1, we have as before $$\deg(f,x_0,\ep) = \deg(f_t,x_0,\ep) = 1\;.$$
Thus by Lemma \ref{lipeo}, $f$ is a local homeomorphism.\end{proof}

As a consequence of Lemmas \ref{sta} and \ref{claim}, we immediately obtain:

\begin{lem}\label{claim2-3} For each $z\in\C^{*m}$, there exists
$\tau_z\in
\R^m,\ w_z\in M_P$ so that \begin{itemize}

\item $\ \mu_\Si(e^{\tau_z/2}\cdot z)=\frac{1}{p}\,
\mu_P(w_z)$;
\item
$\ (\frac{1}{\de}\mu_P(w_z),-\tau_z)\in\ncal^\circ$.
\end{itemize}  Furthermore, $\tau_z,\, w_z$
are unique
(modulo the
$\T$ action on $M_P$) and also satisfy (\ref{ii''}) and (\ref{iii''}).\end{lem}

We shall show in \S\S \ref{3.2.2}--\ref{3.2.3} that along the contour given by
$\tau=\tau_z$, the phase
$\Psi_\C$ has critical points and  maximal real part at $(0,w_z)$.

We recall the concept of transition points mentioned in the introduction: 

\begin{defin} A point $z\in \C^{*m}$ is said to be a {\it transition
point\/} if $-\tau_z$ is in the boundary of the normal cone of $P$ at
$\mu_P(w_z)$, or equivalently if $z$ is in the boundary of one of the regions
$\rcal_F$ described in the introduction.
\end{defin}

\begin{rem} It follows from (\ref{lcalphi}) that $2\,\Log(\rcal_F)=\lcal(F)+C_F$, where
$\Log(z)=(\log|z_1|,\dots,\log |z_m|)$. Hence we can decompose
$\R^m$ as the disjoint union of the sets $\lcal(F)+C_F$.  If $z$ is a transition
point, then $2\Log(z)$ must lie  in the common boundary
of at least two of the sets
$\lcal(F)+C_F$.
Figure~\ref{L} below shows the transition points in log coordinates (as solid lines)
for the case where
$P$ is the polytope of Figure~\ref{f-fan}.

\begin{figure}[htb]
\centerline{\includegraphics*[bb= 0.7in 6.9in 4in 9.6in]{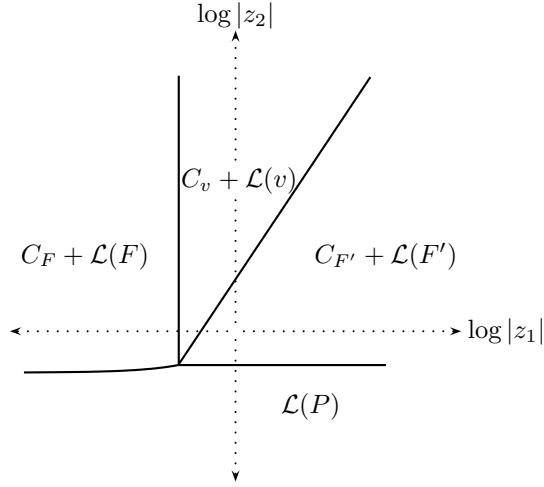}}
\caption{Transition points and regions in log coordinates}
\label{L}\end{figure}
  Note that although $C_{F'}$ is orthogonal to $F'$ in Figure~\ref{f-fan}, $C_{F'}$ is
not orthogonal to $\lcal(F')$.  Similarly, in Lemma~\ref{claim2-3}, if $\mu_P(w_z)\in
F'$, then the orbit $\mu_\Si(e^{r\tau_z}\cdot z)$ is not orthogonal to $F'$. \end{rem}

\subsubsection{Computation of the Hessian}\label{3.2.2} We consider
$$\ccal_z:=\{e^{i\theta}\cdot w_z: e^{i\theta}\in
\T \}=\{e^{i\theta}\cdot w_z:\theta\in
T_F\}\subset M_P\;,$$ which is a (totally real) submanifold of $M_P$ of
dimension
$r$.  We are going to show  that
$\{0\}\times \ccal_z\subset  \T
\times M_P$ is the set of critical points of $\Psi_\C$ with maximal real part
along the $\tau_z$ contour.
In order to apply the method of stationary phase as in the case
where $z$ is in the classically allowable region,  one would also like to show that
$\{0\}\times \ccal_z$ is a  nondegenerate critical submanifold; i.e., the normal
Hessian of the phase  is nondegenerate along $\{0\}\times \ccal_z$.   But the normal
Hessian turns out to be degenerate whenever $z$ is a transition point. For example, if
$z\in\mu\inv_{\Si}(F)$, where
$F$ is a codimension
$r$ face of
$\frac{1}{\de}P$, then the corresponding normal Hessian has an $r$-dimensional
nullspace (as our computations will show) and hence the normal Hessian is degenerate
whenever
$\mu_{\Si}(z)\in
\d (\frac{1}{\de}P)$. However, we shall show:

\begin{lem} Let $z\in\C^{*m}$.  Then $\{0\}\times \ccal_z$ is the (unique) component
of the critical set of
$\Psi_\C(\tau_z,\cdot,\cdot)$ where $\Re \Psi_\C$ attains its maximum on $\T\times
M_P$.  If
$z$ is not a transition point, then $\{0\}\times \ccal_z$ is a nondegenerate critical
submanifold.
\label{crit}\end{lem}

Recall that in the previous section, we verified the proposition for the case where $z$
is in the classically allowed region.  So we now fix a point
$z$ in the classically forbidden region. Let $F_z$ be the face of
$P$ containing
$q_z:=\mu_P(w_z)$, and let
$r$ denote the dimension of $F_z$.  Note that $F-q_z$ an open subset of the
tangent space $T_F$ of
$F\subset
\R^m$.

Let $Y=\mu_P\inv(F_z)$, which is an $r$-dimensional complex submanifold, and
thus
has real dimension $2r$.  We let $V\subset T_{M_P,w_z}$
be the tangent space to the $\R_+^m$ orbit of $w_z$. Recalling (\ref{Xj}), we
see that
$$V=\mbox{span}\{X_1(w_z),\dots, X_m(w_z)\}= \left\{\sum
a_j X_j(w_z): (a_1,\dots,a_m)\in T_F\right\}\;,$$
and hence $\dim V =r$.
(Note that $V=JT^v$, where $T^v$ is the tangent space to the fiber of
$\mu_P$.)

The normal bundle ${\mathbf N}$ to
$\{0\}\times \ccal_z$ can be
decomposed as follows:
\begin{equation}\label{N} {\mathbf N}=T_{\T ,0} \oplus
V \oplus N_Y\;.\end{equation}

To prove Lemma~\ref{crit}, we need the following:

\begin{lem}\label{Hlemma} The normal Hessian of $\Psi_\C(\tau_z,\cdot,\cdot)$ to
$\{0\}\times \ccal_z$ with respect
to a basis for $\mathbf N$ giving the decomposition (\ref{N}) is of the form
$$D_{\mathbf N}^2\Psi_\C(\tau_z,0,w_z)=
\left(\begin{array}{ccc}A&B&0\\
B^t&0&0\\ 0&0& C\end{array}\right)\;,$$
where
\begin{enumerate}
\item[i)]$A$ is an $m\times m$ negative-definite, real symmetric matrix;
\item[ii)] $B$ is an $m\times r$ pure
imaginary matrix of rank $r$;
\item[iii)] $C$ is a $(2m-2r)\times (2m-2r)$ real symmetric matrix, which  is
negative definite if $z$ is not a
transition point.
\end{enumerate}
\end{lem}

\begin{proof} Let
$$D_{\mathbf N}^2\Psi_\C(\tau_z,0,w_z)=\left(\begin{array}{ccc}A&B&\star\\
B^t&\star''&\star'\\ \star^t&\star'{}^t& C\end{array}\right)\;.$$
To simplify notation,
we now write $\tau=\tau_z$.

(i) To compute $A$, we see from (\ref{dphicx}) that
\begin{equation} iD_\phi\Psi_\C|_{\phi=0} =  {\sum_{\alpha \in P}
 |\wh m_\al^P(e^{-\tau/2}\cdot w)|^2
\alpha} + \ \mbox{ a function of } z\;,\label{dphicx2}\end{equation}

As in our previous computation, we then obtain
\begin{eqnarray*}A &=& \mu_P(e^{-\tau/2}\cdot w)\otimes
\mu_P(e^{-\tau/2}\cdot w)
-\sum_{\alpha
\in P} |\wh m_\al^P(e^{-\tau/2}\cdot w)|^2
\alpha \otimes \alpha \\ && +\
\mu_{\de\Si}(e^{\tau/2}\cdot z)\otimes  \mu_{\de\Si}(e^{\tau/2}\cdot z)
-\sum_{\alpha\le p} |\wh m_{\alpha}^{p\Si}(e^{\tau/2}\cdot z)|^2 \alpha\otimes
\alpha\;.
\end{eqnarray*}
By Lemma \ref{A+}, $A$ is the sum of a negative semi-definite term and a
strictly negative term, and hence is negative definite.

(ii) To compute $B$, we similarly obtain
\begin{equation}\label{B2} iD_\rho D_\phi\Psi_\C|_{\phi=0}=
\mu_P(e^{-\tau/2}\cdot w)\otimes
\mu_P(e^{-\tau/2}\cdot w) -\sum_{\alpha
\in P} |\wh m_\al^P(e^{-\tau/2}\cdot w)|^2
\alpha \otimes \alpha \;,\end{equation}
where $D_\rho=(\d/\d\rho_1,\dots,\d/\d\rho_m)=(X_1,\dots,X_m)$. By Lemma
\ref{A+},
$$\big(iD_\rho D_\phi\Psi_\C|_{\phi=0}, \la\otimes\la\big) > 0,
\quad\mbox{for}\ \la\in
T_F\sm\{0\}\;,$$ and hence $B$ has rank $r$.

(iii) To determine $C$, we need to change variables.
Let us recall that $M_P$ is covered by affine coordinate charts centered at
the
vertices of $P$ described as follows:  Let $v=\al^0$ be a vertex of $P$.
Since $P$
is simplicial, there are exactly $m$ edges (1-dimensional faces) incident to
$v$.
Choose $\al^1,\dots,\al^m$ on the respective edges. Let $\be^j=\al^j-\al^0$;
then
$\{\be^1,\dots,\be^m\}$ form a basis for $\R^m$. Since we are assuming
that
$P$ is Delzant, we can choose  the $\al^j$ such that
$\{\be^1,\dots,\be^m\}$
generate $\Z^m$.  We consider the $m\times m$ matrix of integers $$\De=
\big(\De_{jk}\big)= \big(\be^j_k\big)\;.$$  Then $\det \De=\pm 1$ and the
matrix
$\Ga:=\De\inv$ has integer entries.  Consider the (Zariski) neighborhood of
$v$:
$$U_v:= \{w\in M_P:\chi_v(w)\neq 0\}= \bigcup\{F: F \ \mbox {is a face
of} \ P,\
v\in \bar F\}\;.$$ (The $U_v$  cover $M_P$.) Let $\eta_j=w^{\be^j}$,
$j=1,\dots,m$.  We note that the monomial mapping
$$\eta=(\eta_1,\dots,\eta_m):\C^{*m}\to \C^{*m}$$ is bijective, its
inverse
given by $$w_j=\eta_1^{\Ga_{j1}}\cdots \eta_m^{\Ga_{jm}}\;.$$

Let $F^1,\dots,F^m$ denote the facets of $P$ incident to $v$, indexed so
that
$\al_j\not\in F^j$. Let $x\in F^j$ be arbitrary. We then have
$$\eta_j=\frac{\chi_{\al^j}(w)}{\chi_{\al^0}(w)}\to
\frac{\chi_{\al^j}(x)}{\chi_{\al^0}(x)}= 0 \quad
\mbox{as}\ \mu_P(w)\to x\in F^j\;.$$
(Since $\al^0\in\overline {F^j}$, $\wh \chi_{\al^0}(x)\neq 0$.) Therefore
$\eta$ extends to an isomorphism
$\eta:(U_V,v)
\buildrel{\approx}\over\to
(\C^m,0)$. Furthermore, the facets incident to $v$ are the images of the
divisors
$\{\eta_j=0\}$ under the moment map.  (This shows that $M_P$ is nonsingular
if the
Delzant condition is satisfied.)

We note that the standard unit vectors $e^j=\Ga \beta^j = \Ga \al^j -\Ga v$.
Consider the
nonhomogeneous linear map
$$\wt\Ga:\Z^m\to \Z^m\;,\quad u\mapsto \Ga u -\Ga v\;.$$
Then $\wt\Ga(\al^0)=0,\ \wt\Ga(\al^j)= e^j$ ($1\le j\le m$). Consider the
polytope
$Q=\wt\Ga(P)$. Under the map $\wt\Ga$, the closed facets $\overline {F^j}$
correspond to polytopes
in the coordinate hyperplanes $\{x_j=0\}$.  Hence the polytope $Q$ lies in
the upper
quadrant
$\{x_j\ge 0, 1\le j\le m\}$. We also let $$\wt\De = \wt \Ga\inv:\Z^m\to
\Z^m\;,\quad
u\mapsto \De u +v\;.$$

Now let us rewrite the phase in terms of the $\eta$ coordinates.  We first
note that
$$\chi_\al(w)=\chi_{\wt\Ga(\al)}(\eta),\quad w\in \C^{*m},\quad
\eta_j=w^{\be^j}\;.$$

Hence, \begin{equation}\label{newcoord}H_P(\tau,w)=\frac{\sum_{\alpha \in P}
e^{- \langle
\alpha, \tau \rangle}|c_{\al}\chi_{\wt\Ga(\al)}(\eta)|^2}{ \sum_{\alpha \in
P}
|c_{\al} \chi_{\wt\Ga(\al)}(\eta)|^2}=\frac{\sum_{\ga\in Q}e^{- \langle
\wt\De(\ga), \tau \rangle}| c_{\wt\De(\ga)}\chi_\ga(\eta)|^2}{ \sum_{\ga\in
Q}|c_{\wt\De(\ga)} \chi_\ga(\eta)|^2} \;.\end{equation}

We first consider the case where $w_z=v$; i.e., $r=0$.  Expanding
(\ref{newcoord}) in
powers of
$\eta_j$, we have
\begin{eqnarray*}H_P(\tau,w)&=& \frac{|c_v|^2e^{-\langle
v,\tau\rangle}+\sum_{j=1}^m
|c_{\al^j}|^2e^{-\langle
\al^j,\tau\rangle} |\eta_j|^2 + \cdots}{|c_v|^2+\sum_{j=1}^m
|c_{\al^j}|^2 |\eta_j|^2 + \cdots}\\
&=& e^{-\langle v,\tau\rangle} \left[1+\sum_{j=1}^mc'_j\left(e^{-\langle
\be^j,\tau\rangle} -1\right) |\eta_j|^2 +\cdots\right]\;,\end{eqnarray*}
where $c'_j=|{c_{\al^j}}/{c_{\al^0}}|^2>0$. Since  $-\tau$ is in the normal cone
at $v$,we have  $\langle -\tau,v\rangle \ge \langle -\tau,\al_j\rangle$ and
hence
$\langle
\be_j,\tau\rangle \ge 0$, for
$1\le j\le m$.
{\bf Now suppose that $z$ is not a
transition point}; i.e., $-\tau$ is not in the boundary of the normal cone at
$\{v\}$.  Then  $\tau$ is not
perpendicular to
any of the rays along the edges incident to $v$; so $\langle
\be_j,\tau\rangle \ne 0$, and therefore
$\langle \be_j,\tau\rangle >0$. Hence
$$\log H_P(\tau,w) = -\langle v,\tau\rangle +\sum_{j=1}^m b_j|\eta_j|^2 + \cdots,
\quad b_j= c'_j\left(e^{-\langle
\be^j,\tau\rangle} -1\right)< 0\;.$$ This leads to two conclusions: (i) $d\log
H_P|_{w=w_z}=0$ and hence $(0,w_z)$ is a critical point of
$\Psi_\C(\tau,\cdot,\cdot)$; (ii)
$C=D_w^2\log H_P|_{w=w_z}$ is a  diagonal matrix with negative eigenvalues
$b_1,b_1,\dots,b_m,b_m$. Note that $C$ is singular whenever $z$ is in the
boundary of the normal polyhedron of $v$.

Now consider the case $1\le r\le m-1$.  (If $r=m$, then $z$ is in the
allowable region.)
We can assume without loss of generality that $Y=\mu_P\inv(F_z)$ is given by
\begin{equation}\eta_1=\cdots = \eta_{m-r}=0, \qquad \eta_{m-r+1}\ne 0,\
\dots\ ,\ \eta_m\ne
0\;.\label{facecoord}\end{equation}  In order to  decompose the sum in
(\ref{newcoord}), we consider the sets of lattice points $$E_j: = \{\ga\in
Q: \ga_k=\delta_k^j \
\mbox{for}\ 1\le k\le m-r\}=Q\cap\textstyle\left\{e^j+\sum_{l=m-r+1}^m s_l
e^l\right\}\;,$$
for $1\le j\le m-r$. (These are the indices of the only terms that
contribute to the normal
Hessian at $w_z$.) We also let $E_0=Q\cap (\{0_{m-r}\}\times \Z^r)$. Since
$\langle\be^l,
\tau\rangle = 0$ for $l> m-r$, it follows that
$$\langle
\wt\De(\ga), \tau \rangle= \langle
\wt\De (e^j), \tau \rangle = \langle \al^j, \tau \rangle\qquad \forall
\ga\in E_j\;,$$ and
similarly
$$\langle
\wt\De(\ga), \tau \rangle=  \langle \wt\De(0), \tau \rangle= \langle v, \tau
\rangle\qquad
\forall \ga\in E_0\;.$$

We write $$\eta=(\eta',\eta''),\quad \eta'=(\eta_1,\dots,\eta_{m-r}),\
\eta''=(\eta_{m-r+1},\dots,\eta_{m})\;.$$ Let
$E_j''\subset \Z^m$, for $j=0,1,\dots,m-r$, be given by
$$\{0_{m-r}\}\times E_0''=E_0\,,\quad
\{0_{m-r}\}\times E_j''=E_j-e^j \quad (1\le j\le m-r)\;.$$
Thus we rewrite
(\ref{newcoord}):
\begin{eqnarray}H_P(\tau,w)&=&\frac{e^{-
\langle v, \tau \rangle} \sum_{\ga\in E_0}|c_{\wt\De(\ga)}\chi_\ga(\eta)|^2
+\sum_{j=1}^{m-r}\left[e^{-
\langle \al^j, \tau \rangle}\sum_{\ga\in
E_j}|c_{\wt\De(\ga)}\chi_\ga(\eta)|^2\right]+\cdots}{\sum_{\ga\in
E_0}|c_{\wt\De(\ga)}\chi_\ga(\eta)|^2 +\sum_{j=1}^{m-r}\left[\sum_{\ga\in
E_j}|c_{\wt\De(\ga)}\chi_\ga(\eta)|^2\right]+\cdots}\nonumber\\
&=&e^{-\langle v, \tau \rangle}\frac{\sum_{\nu\in
E_0''}c_{0\nu}'|\chi_\nu(\eta'')|^2+\sum_{j=1}^{m-r}\left[\sum_{\nu\in
E_j''}c_{j\nu}'|\chi_\nu(\eta'')|^2\right]e^{-
\langle \be^j, \tau \rangle}|\eta_j|^2+\cdots}{\sum_{\nu\in
E_0''}c_{0\nu}'|\chi_\nu(\eta'')|^2+\sum_{j=1}^{m-r}\left[\sum_{\nu\in
E_j''}c_{j\nu}'|\chi_\nu(\eta'')|^2\right]|\eta_j|^2+\cdots}\nonumber\\
&=&e^{-\langle v,\tau\rangle}
\left[1+\sum_{j=1}^{m-r}\frac{\la_j(\eta'')}{\la_0(\eta'')}
\left(e^{-\langle
\be^j,\tau\rangle} -1\right) |\eta_j|^2
+O(\|\eta'\|^4)\right]\;,\label{newcoord2}\end{eqnarray}
where $c_{j\nu}'\in\R_+$ and  $$\la_j(\eta'')=\sum_{\nu\in
E_j''}c_{j\nu}'|\chi_\nu(\eta'')|^2 = c_{j0}' + \sum_{\nu\in
E_j''\sm\{0\}}c_{j\nu}'|\chi_\nu(\eta'')|^2 >0 \qquad (0\le j\le m-r)\;.$$

(Note that (\ref{newcoord2}) implies that
$H_P(\tau,w)= e^{-\langle v,\tau\rangle}$ for $w\in
F_z$, as we observed earlier, and  that $(0,w_z)$ is a critical point of
$\Psi_\C$.)
{\bf Suppose as before that $z$ is not a
transition point}; i.e., $-\tau$ is not in the boundary of the normal cone at
$F_z$. Then $\langle
\be_j,\tau\rangle\ne 0$ and in fact $\langle\be_j,\tau\rangle > 0$, for
$1\le j\le m-r$,  since $-\tau$ is in the cone at $F_z$.
It now follows from (\ref{newcoord2}) that
$C=D^2_{\eta'}\log
H_P|_{\{\eta'=0\}}$ is diagonal (with respect to the $\eta'$ coordinates) with negative
eigenvalues.

(iv) To complete the proof of the lemma, we now show that the other blocks of
the normal Hessian
$D_{\mathbf N}^2\Psi_\C$ vanish.  We first note that since
$H_P(\tau_z,\cdot)$ is constant on
$Y$, the matrix
$\star''$
vanishes. We also immediately  conclude from (\ref{newcoord2}) that
$$\star'=
D_{\rho''}D_{\eta'}\log H_P|_{\{\eta'=0\}}=0\;,$$
where $\rho''=(\Re \eta_{m-r+1},\dots, \Re \eta_m)$.

It remains to show that $\star=0$.  Assuming (\ref{facecoord}), we expand
(\ref{dphicx2}) as
before:
\begin{eqnarray*} iD_\phi\Psi_\C|_{\phi=0} &=& \frac{\sum_{\ga\in Q}e^{-
\langle
\wt\De(\ga), \tau \rangle}| c_{\wt\De(\ga)}\chi_\ga(\eta)|^2\wt\De(\ga)}{
\sum_{\ga\in
Q}|c_{\wt\De(\ga)} \chi_\ga(\eta)|^2} + f(z)\\
&=& e^{-\langle v,\tau\rangle}
\left[\vec g_0(\eta'')+\sum_{j=1}^{m-r}|\eta_j|^2 \vec g_j(\eta'')
+O(\|\eta'\|^4)\right] + f(z)\;.
\end{eqnarray*} Hence $$\star = D_{\eta'}D_\phi\Psi_\C|_{\phi=0,\eta'=0}\
=0\;.$$
\end{proof}

\medskip
\subsubsection{Proof of Lemma \ref{crit}}\label{3.2.3}

It follows from (\ref{newcoord2}) that $(0,w_z)$ is a critical point of
$\Psi_\C(\tau,\cdot,\cdot)$ and hence $\{0\}\times \ccal_z$ is contained in the
critical set. We conclude from Lemma
\ref{Hlemma} that  the real
$r\times r$ matrix  $B^t A\inv B= (iB)^t(-A\inv)(iB)$
is positive definite, and hence
\begin{equation}\label{det}\det D_{\mathbf N}^2\Psi_\C=  \det A\; \det (B^t
A\inv B)\;
\det C
\neq 0\;,\end{equation} if $z$ is not a transition point.
Therefore $\{0\}\times \ccal_z$ is a
nondegenerate critical manifold of $\Psi_\C$ whenever $z$ is not a
transition
point.

Recalling (\ref{phaseH}), we set
\begin{equation} \label{CRITVAL}b(z):=-\Psi_{\C}(\tau_z,0,w_z)  =- \de
\log  H_\Si(-\tau_z, z)-\log H_P(\tau_z,w_z)\;.\end{equation}
We now show that
the maximum of
$\Re\Psi_\C$ on the
$\tau$-contour is attained on $\{0\}\times \ccal_z$ and hence the maximum equals
$-b(z)$. We recall that by (\ref{maxatphi0}), it suffices to show that
$\Re\Psi(\tau,0,w)$ attains its maximum at $w=w_z$, and by (\ref{phaseH}), this is
equivalent to establishing that $\log H(\tau,\cdot)$, whose Hessian is $C$, takes its
maximum at $w_z$. We first suppose that $z$ is not a transition point. Note that by the
above,
$H(\tau,\cdot)$ is constant and has a local maximum along
$Y=\mu_P\inv(F_z)$.  Suppose on
the contrary that there is a point $\wt w\in M_P$ such that
\begin{equation}\label{bigger}H(\tau,\wt w)=\sup_{w\in M_P} H(\tau,w) >
H(\tau,w_z)\;.\end{equation} Let $\wt F$ denote the face that contains
$\mu_P(\wt w)$. Since $H(\tau,\cdot)$ is constant on $F_z$, by (\ref{bigger})
the face
$\wt F$ is not contained in
$\overline{F_z}$, and hence $C_{\wt F}\cap C_{F_z}^\circ=\emptyset$, where
$C_F^\circ$ denotes the nonboundary points of the cone $C_F$. Since
$-\tau\in  C_{F_z}^\circ$ (by the assumption that $z$ is not a transition
point), it follows that
$-\tau\notin C_{\wt F}$.   Then we conclude from
(\ref{newcoord2}) that at least one of the eigenvalues of $C(\wt w)$ is
positive, contradicting (\ref{bigger}). On the other hand, if $z$ is a
transition point, then $-\tau\in  C_{F'}^\circ$, where $F'$ is a face whose
closure contains $F_z$.  Then $H(\tau,w_z)=H(\tau,w')$ for $w'\in F'$ and
$H(\tau,\cdot)$ attains a local maximum along $F'$.  Then by the above
argument, $H(\tau,w')$ is a global maximum. Recalling (\ref{CRITVAL}), we then
have
\begin{equation}\label{max-b}-b(z)=\Psi_\C(\tau_z,0,w_z) = \sup_{w\in M_P}
\Psi_\C(\tau_z,0,w_z) =\sup_{\phi\in\T , w\in M_P}  \Re\Psi_\C(\tau_z,\phi,w_z)
\;.\end{equation}  This completes the proof of Lemma \ref{crit}.\qed

\subsubsection{Asymptotic expansions on the classically forbidden region.}\label{3.2.4}
We continue the proof of Proposition~\ref{SZEGO}(ii).  Equation (\ref{max-b})
says that
$\Psi +b(z)$ is a phase function of positive type that takes the value 0 on the
critical set
$\{0\}\times \ccal_z$.  {\bf If $z$ is not a transition point,} the points of
$\{0\}\times \ccal_z$ are nondegenerate critical points.  As before, we obtain by the
method of stationary phase  an asymptotic expansion
\begin{equation}\label{asympt}K(z) =N^{\frac{m+r}{2}} e^{-N b(z)}[c_0^F(z) +
c_1^F(z)N\inv+ c_2^F(z)N^{-2}+\cdots]\end{equation} valid over each open region
$\rcal_F^\circ$. (Since the critical submanifold $\{0\}\times\ccal_z$ has codimension
$3m-r$, the leading term of the expansion of $e^{N b(z)}K(z)$ contains
$N^{2m-\half(3m-r)}= N^{\frac{m+r}{2}}$.) To show that
$\Pi_{|NP}(z,z)$ also has the asymptotic expansion (\ref{asympt}), we must show that
\begin{equation}\label{verysmall} S_N(z) = O(N^{-k}  e^{-N b(z)})\qquad \forall\
k\in\N\;,\end{equation} where the bound is locally uniform in $\ccal^j$ for all $j$.
We deform the $\T$ integral in  (\ref{szego-rem}) to an integral over
$\{|\zeta|=\tau_z\}$ as before, so that we have
\begin{eqnarray} S_N(z)&=&\frac{1}{(2\pi)^m} \int_{M_P} \int_{M_P}
\int_{\T } \wt
R_N(w,\eta)\left( \sum_{|\al| \le \de} e^{ \langle
\alpha, \tau_z +i\phi \rangle} |\wt{m}_{\alpha}^{\de
\Sigma}  (z)|^2\right)^N\nonumber\\ &&\hspace{1.2in} \times\ \left(\sum_{\al\in
P} e^{-
\langle
\alpha, \tau_z +i\phi \rangle}
\wh m_\al^P(\eta)\overline{\wh m_\al^P(w)}\right)^N\, d\phi\, dw\,
d\eta,\end{eqnarray} where $\wt R_N=\frac{(Np+m)!}{(Np)!}R_N$ is
rapidly decaying. We first note that
$$ \left|\sum_{|\al| \le \de} e^{ \langle
\alpha, \tau_z +i\phi \rangle} |\wh{m}_{\alpha}^{\de
\Sigma}  (z)|^2\right| \le \sum_{|\al| \le \de} e^{ \langle
\alpha, \tau_z  \rangle} |\wh{m}_{\alpha}^{\de
\Sigma}  (z)|^2=  H_{p\Si}(-\tau_z,z)\;.$$  By the
Cauchy-Schwartz inequality and the fact established above that $H_P(\tau_z,w)$
takes its maximum at $w=w_z$, we have
\begin{eqnarray*}\left|\sum_{\al\in P}
e^{-
\langle
\alpha, \tau_z +i\phi \rangle}
\wh m_\al^P(\eta)\overline{\wh m_\al^P(w)}\right|
&\le& \sum_{\al\in P}
e^{-
\langle
\alpha, \tau_z/2 \rangle}
|\wh m_\al^P(\eta)|\;e^{-
\langle
\alpha, \tau_z/2 \rangle}|\wh m_\al^P(w)|\\
&\le& \left(\sum_{\al\in P}
e^{-
\langle
\alpha, \tau_z \rangle}
|\wh m_\al^P(\eta)|^2\right)^\half\left(\sum_{\al\in P}
e^{-
\langle
\alpha, \tau_z \rangle}
|\wh m_\al^P(w)|^2\right)^\half\\
&=& H_P(\tau_z,\eta)^\half  H_P(\tau_z,w)^\half \ \le \ H_P(\tau_z,w_z)\;.
\end{eqnarray*}
Therefore, we can write
$$S_N(z)=\frac{1}{(2\pi)^m} \int_{M_P} \int_{M_P}
\int_{\T } \wt
R_N(w,\eta) e^{N\hat\Psi(\phi,w,\eta;z)}\,d\phi\,dw\,d\eta\;,$$
where the phase satisfies
$$\Re\hat\Psi \le \log \big[H_{p\Si}(-\tau_z,z) H_P(\tau_z,w_z)\big]
=e^{-b(z)}\;.$$
The estimate (\ref{verysmall}) now follows from the fact that $\wt R_N$ is rapidly
decaying.

Thus, we have verified the
asymptotic formula (ii) of Proposition~\ref{SZEGO}, together with (a) and (b),
and we have shown that
$b(z)$ is given by (\ref{CRITVAL}). 

\begin{rem} It follows from (\ref{szego3}), (\ref{maxatphi0}), and (\ref{max-b}) that
in a neighborhood
$U\subset\!\subset\C^{*m}$ of a transition point $z^0$, we have uniform bounds of the
form
\begin{equation}\label{upperbound}\Pi_{|NP}(z,z)\le C_{U}N^{2m}
e^{-Nb(z)}\;.\end{equation} But a sharper inequality should hold;  for example, we know
that if
$z^0\in\d\acal_P$, then $b(z^0)=0$ and $$ \Pi_{|NP}(z^0,z^0)\le
\Pi_{Np}^{\CP^m}(z^0,z^0) = [p^m+o(1)]N^m\;.$$
It is an open question whether we have asymptotic expansions of the form (\ref{asympt})
at transition points.  See \S\ref{airy-tale} for a discussion of the issues involved in
studying the asymptotics at transition points.
\end{rem}

\subsubsection{The decay function} \label{s-decay} In this section we verify (c)--(e)
of Proposition~\ref{SZEGO}(ii) and then prove Proposition~\ref{convergence}. We first
note that when 
$z\in{\acal_P}$, we have  
$\tau_z=z$ and hence by (\ref{CRITVAL}), $b=0$ on $\acal_P$. Furthermore, since $z\mapsto
\tau_z$ is easily seen to be continuous on $\C^{*m}$ and  $\ccal^\infty$ on each
$\overline{\rcal_F}$, the same holds for $b(z)$.

We now
show that
$b$ is
$\ccal^1$  by computing its derivative.
Recalling that
$b(z)=-\Psi_\C(\tau_z,0,w_z;z)$, we have
$$ - D_z b = D_\tau \Psi_\C|_{(\tau_z,0,w_z,z)} \cdot
D_z\tau_z + D_w \Psi_\C|_{(\tau_z,0,w_z,z)} \cdot D_z w_z + D_z
\Psi_\C|_{(\tau_z,0,w_z,z)}\;. $$
Here, $D_z\tau_z$ and $D_zw_z$ are only piecewise continuous, being
discontinuous at transition points. By (\ref{dphicx0})--(\ref{ii''}) and
the fact that
$(0,w_z)$ is a critical point of $\Psi_\C(\tau_z,\cdot,\cdot)$, we have
$$D_\tau \Psi_\C|_{(\tau_z,0,w_z,z)} = -iD_\phi
\Psi_\C|_{(\tau_z,0,w_z,z)}=0\;.$$
The critical point condition also says that $D_w \Psi_\C|_{(\tau_z,0,w_z,z)}=0$.
Recalling (\ref{phaseH}), we thus have
\begin{equation}\label{Db} D_zb= - D_z
\Psi_\C(\tau,0,w;z)|_{\tau=\tau_z, w=w_z} = -\de D_z \log
H_\Si(-\tau,z)|_{\tau=\tau_z}\in \ccal^0(\C^{*m})\;.
\end{equation} Therefore, $b\in \ccal^1(\C^{*m})$, verifying (e).

To verify (c), we need to find a more explicit formula for the derivative,
using the log coordinates
$\zeta_j=\rho_j + i\theta_j=\log z_j$. We consider the map
$\xi=(\xi_1,\dots,\xi_m):\C^{*m}\to \C^{*m}$ given by
\begin{equation}\label{xi}\xi(z)=e^{\tau_z/2}\cdot z= e^{\tau_z/2
+\rho+i\theta}\;,\end{equation} so that
$\mu_\Si
\circ
\xi (z)=\frac{1}{\de}\mu_P(w_z)$. By (\ref{Db}),
\begin{eqnarray} D_\rho b &=& -\de D_\rho \log
H_\Si(-\tau,e^\rho)|_{\tau=\tau_z}\nonumber\\ &=&-p D_\rho
\log \left.\left(1+\sum
e^{\tau_j+2\rho_j}\right)\right|_{\tau=\tau_z}+ pD_\rho\log \left(1+\sum
e^{2\rho_j}\right)
\label{Db1}\;.\end{eqnarray}
We have \begin{eqnarray}D_\rho
\log \left.\left(1+\sum
e^{\tau_j+2\rho_j}\right)\right|_{\tau=\tau_z} &=&
\left(
\frac{ 2e^{\tau_1+2\rho_1}}{1+\sum
e^{\tau_j+2\rho_j}},\dots, \frac{ 2e^{\tau_m+2\rho_m}}{1+\sum
e^{\tau_j+2\rho_j}}\right)\nonumber \\
&= &\left(
\frac{2|\xi_1|^2}{1+\|\xi\|^2},\dots \frac{2|\xi_m|^2}{1+\|\xi\|^2}\right) \ =
\ 2\,\mu_\Si\circ\xi\;,
 \label{Dlog}\end{eqnarray} and similarly, $D_\rho \log \left(1+\sum
e^{2\rho_j}\right) = \ 2\,\mu_\Si$. Hence
\begin{equation}\label{Db2} D_\rho b = 2p(\mu_\Si - \mu_\Si\circ\xi)\;.\end{equation}

Let $z\in \C^{*m}\sm\overline{\acal_P}$.  If
$z$  is not a transition  point,  then the
asymptotic expansion holds at $z$, and as we observed at the beginning of
\S\ref{s-forbidden}, $\Pi_{|NP}(z,z)$ is rapidly decaying; hence $b(z)>0$. Now
suppose on the contrary that
$z$ is a transition point and $b(z)=0$.  Then $z$ would be a
critical point for $b$, which is impossible since $Db|_z\neq 0$ whenever
$z\notin\overline{\acal_P}$ by (\ref{Db2}).  This completes the proof of (c).

It remains to verify formula
(\ref{b}). First we establish (\ref{b-action}): Fix a point $z$ in the classically
forbidden region, and make the change of variables $\rho_j=
\sigma_j/2 + \log |z_j|$. Then by (\ref{Db2}),
\begin{equation}b(e^{\tau_z/2}\cdot z)) - b(z) = \int_0^{\tau_z} D_\sigma
b(e^{\sigma/2}\cdot z)
\cdot d\sigma = p \int_0^{\tau_z} \big[\mu_\Si(e^{\sigma/2}\cdot z) -
\mu_\Si\circ\xi(e^{\sigma/2}\cdot z) \big] \cdot d\sigma\;.\label{FTC}\end{equation}
By definition, $q(z)=\mu_P(w_z) = p\mu_\Si\circ\xi(z)$, and hence
$p\mu_\Si\circ\xi(e^{\sigma/2}\cdot z) = q(e^{\sigma/2}\cdot z)$. Furthermore
$b(e^{\tau_z/2}\cdot z))= 0$ since $e^{\tau_z/2}\cdot z=\xi(z)\in \d \acal_P$, and
formula (\ref{b-action}) follows from (\ref{FTC}).

Note that the integral in (\ref{b-action}) is independent of the path from $0$ to
$\tau_z$. Choosing the path $\sigma(r)=r\tau_z$, $0\le r\le 1$, we have
\begin{equation}\label{b1}\int_0^{\tau_z}q(e^{\sigma/2}\cdot z)\cdot d\sigma
=\int_0^{\tau_z}q(z)\cdot d\sigma = \langle q(z),\tau_z\rangle\;.\end{equation}
On the other hand, since $D_\sigma\log(1+\|e^{\sigma/2}\cdot z\|^2)=\mu_\Si(e^{\sigma/2}\cdot
z)$, we have
\begin{equation}\label{b2}\int_0^{\tau_z} \mu_\Si(e^{\sigma/2}\cdot
z) \cdot d\sigma = \log(1+\|e^{\tau_z/2}\cdot z\|^2) - \log(1+\| z\|^2)
\;.\end{equation} Substituting (\ref{b1})--(\ref{b2}) in (\ref{b-action}), we obtain
formula (\ref{b}). This completes the proof of part (ii) of Proposition~\ref{SZEGO}.

\bigskip We note that $b$ fails to be $\ccal^2$ at transition points.  This
follows from (\ref{Db2}), (\ref{xi}), and the fact that $z\mapsto\tau_z$ is
piecewise $\ccal^1$ with discontinuous derivative at the interfaces between
the regions $\rcal_F$.

\begin{rem} An alternate derivation of (\ref{b}) is as follows:
By an argument similar to the above, we obtain \begin{equation}
\log H_P(\tau,w) = \int_0^{\tau} \mu_P(e^{\sigma/2}\cdot w)
\cdot d \sigma\;,\end{equation} and therefore
\begin{equation}\label{HP-simple}\log H_P(\tau_z,w_z)= \langle
-\mu_P(w_z),\tau_z\rangle\;.\end{equation}\end{rem}
Formula (\ref{b}) is an immediate consequence of (\ref{CRITVAL}) and (\ref{HP-simple}).

 We now prove Proposition \ref{convergence}: By Proposition~\ref{SZEGO}(ii) and 
(\ref{CRITVAL}),
\begin{equation}\label{limit}\frac{1}{N} \log \Pi_{|NP}(z,z)
\to -b(z) = \de
\log  H_\Si(-\tau_z, z)+\log H_P(\tau_z,w_z)\;,
\end{equation} for all non-transition points $z$.
Let \begin{equation}\label{uN}u_N(z)=\frac{1}{N}\log \Pi_{|NP}(z,z) +\de
\log(1+\|z\|^2)=\frac{1}{N} \log \sum_{\al\in NP}
{\textstyle{N\de \choose \al}} e^{2\langle
\al,\rho\rangle}\;,\quad z\in\C^{*m}\;.\end{equation}
Thus for all non-transition points $z$, we have
\begin{equation}\label{u0}
u_N(z)\to u_\infty(z):= \de
\log(1+\|z\|^2) -b(z)= \de\log (1+\|e^{\tau_z/2}\cdot z\|^2)
+\log H_P(\tau_z,w_z)\;,\end{equation}

We must show that convergence of (\ref{u0}) also holds at the transition points
and is uniform on compact subsets of $\C^{*m}$.   We again use the log coordinates
$\zeta_j=\rho_j + i\theta_j=\log z_j$, so that
\begin{equation}D_\rho u_N = \frac{1}{N}D_\rho \log \sum_{\al\in NP}
{\textstyle{N\de \choose \al}} e^{2\langle
\al,\rho\rangle}= \frac{\sum_{\al\in NP}
{\textstyle{N\de \choose \al}}e^{2\langle
\al,\rho\rangle}2\al}{N\sum_{\al\in NP}
{\textstyle{N\de \choose \al}}e^{2\langle
\al,\rho\rangle}}
=\frac{2}{N}\mu_{NP}\;,
\label{DuN}\end{equation}
and therefore $\|D_\rho u_N\| \le 2\de$.  Since $\{u_N(z^0)\}$ converges for any
non-transition point  $z^0$, it follows that $\{u_N\}$ is uniformly bounded and
uniformly equicontinuous on compact sets.  Therefore it converges uniformly
on compact sets in   $\C^{*m}$.\qed

\subsubsection{Precise asymptotics on the classically allowed region.} \label{3.2.6}
In this section, we prove part (i) of Proposition~\ref{SZEGO}. Since
$$\Pi_{Np\Si}(z,z)\equiv
\frac{\dim H^0(\CP^m,\ocal(Np))}{\vol(\CP^m)}={Np+m\choose m}\,m!
=\prod_{j=1}^m (Np+j)\;,$$
we have by (\ref{szego-proj}) and (\ref{mchi}), $$\Pi_{|NP}(z,z) = \prod_{j=1}^m (Np+j)
+ R_N(z)\;,$$ where
\begin{equation} \label{Piallow}R_N(z)= -\sum_{\al\in
Np\Si\sm NP}\
\frac{1}{\|\chi_\al\|^2}
|\wh\chi_\al(z)|^2
= -\left[\prod_{j=1}^m (Np+j)\right] \sum_{\al\in
Np\Si\sm NP}\
|\wh m_\al^{Np\Si}(z)|^2\;.\end{equation}
Since $\#(Np\Si\sm NP) = O(N^m)$,  Proposition~\ref{SZEGO}(i) is an immediate
consequence of  (\ref{Piallow}) and the following uniform asymptotic decay estimate
for monomials.

\begin{lem}\label{mondecay} Suppose that $K\subset \C^{*m}$ is compact, and let
$U\subset\Si$ be an open neighborhood of $\mu_\Si(K)$. Then there exist  positive
constants $C_l=C_{l,K,U},\ \la=\la_{K,U}$ such that $$\left\|\wh
m_\al^{N\Si}\right\|_{\ccal^l(K)}\le C_le^{-\la N} \qquad \mbox{whenever }
\frac{\al}{N}\in
\Si\sm U\;,\ N\ge 1\;, l\ge 0\;.$$
\end{lem}

\begin{proof} We first prove the $\ccal^0$ estimate. Choose a positive integer
$k>{2m^2}/\mbox{dist}\left(\Si\sm U,\mu_\Si(K)\right)$. We consider the finite
collection of polytopes
$$Q_\be:=2m\Si +\be\subset k\Si\;,\qquad \be\in\Z^m\cap(k-2m)\Si\;.$$  
For each $\be$ as above, we let $b_\be(z)$ be
the function on $\C^{*m}$ given by formula (\ref{b}) with $P$ replaced by $Q_\be$ and
with $p$ replaced by $k$.   We claim that the $\ccal^0$ estimate holds with
$$\la=\frac{1}{k+\ep} \inf\left\{b_\be(z): z\in K,\ \be\in\Z^m\cap(k-2m)\Si,\ \frac 1k
Q_\be\subset
\Si\sm U\right\}>0\;.$$

Suppose on the contrary that the estimate does not hold with this value
of $\la$.  Then we can choose an increasing sequence of integers $\{N_j\}$ and a
sequence
$\al^j\in N(\Si\sm U)\cap\Z^m$ such that \begin{equation}\label{diverges}e^{\la
N_j}\sup_K\left|\wh m_{\al^j}^{N_j\Si}\right|^2 \to +\infty\;.\end{equation}  By
passing to a subsequence, we can assume that $\frac{\al^j}{N_j}\to x\in \Si\sm U$. 
Choose
$\be\in\Z^m\cap(k-2m)\Si$ such that
$$x\in
\mbox{interior}_\Si\left(\frac 1k Q_\be\right)\;,$$
and hence $k{\al^j}\in N_j Q_\be$ for $j\gg 1$. 
Since dist$(x,K)\ge \frac{2m^2}{k}\ge \mbox{diam}(\frac1k Q_\be)$, it follows that
$(\frac 1k Q_\be)\cap K = \emptyset$.  Hence by Proposition~\ref{convergence} applied
to the polytope $Q_\be$, there is a positive constant
$A$ such that 
\begin{equation}\Pi_{|NQ_\be}(z,z) \le A  \textstyle\exp\left(-\frac k{k+\ep}
b_\be(z)N\right)
\le A e^{-\la  kN } \qquad
\forall z\in K,\ \forall N\ge 1\;.\label{Qdecay}\end{equation}

We shall use the multinomial inequality:
\begin{equation}\label{multi} {N\choose\al}^k \le {kN\choose k\al}\;.\end{equation}
One way to verify (\ref{multi}) is to regard ${N\choose\al}$ as the partition
function $\pcal_N(\al)$ for  the polytope $\Si$; i.e., the number of ways to
write $\al=\be^1+\cdots +\be^N$, where $\be\in\Si$ (see \S \ref{s-exact}). Now,
given $k\al\in kN\Si$ we first write $k\al = \al + \cdots + \al$ ($k$ times) and then
decompose each $\al$ as above, resulting in ${N\choose\al}^k$ different
decompositions of $k\al$ as a sum of $kN$ elements of $\Si$. This gives us the lower
bound for ${kN\choose k\al}$ in (\ref{multi}).

We thus have
$$|\wh m_{\al^j}^{N_j\Si}|^{2k}= {N_j\choose \al^j}^k |\wh \chi_{\al^j}|^{2k}
\le {kN_j\choose k\al^j} |\wh \chi_{k\al^j}|^2 = |\wh m_{k\al^j}^{kN_j\Si}|^2\;.$$
Since $k{\al^j}\in N_j Q_\be$, it then follows from (\ref{Qdecay}) that
$$|\wh m_{\al^j}^{N_j\Si}(z)|^{2}\le  |\wh m_{k\al^j}^{kN_j\Si}|^{2/k} \le
\Pi_{|N_jQ_\be}(z,z)^{1/k} \le A^{1/k} e^{- \la  N_j}\;,$$ for all $z\in K$, which
contradicts (\ref{diverges}).  Hence the $\ccal^0$ estimate holds.

We now verify the $\ccal^1$ estimate:  We have
$$|\wh m_\al(z)| = {N\choose \al}^\half h(z)^{\frac N2}|z^\al|\;,\qquad h(z)=\frac
1{1+\|z\|^2}\;.$$  Differentiating, we obtain
 
$$D_\rho |\wh m_\al(z)| = {N\choose \al}^\half h(z)^{\frac N2}|z^\al| \left[ \al
+\frac N2 D_\rho\log h\right] = |\wh m_\al(z)| \left[ \al
+\frac N2 D_\rho\log h\right]
\;,$$ and hence the $\ccal^1$ estimate follows  from the $\ccal^0$ estimate. 
Differentiating repeatedly, we obtain all the $\ccal^l$ estimates.
\end{proof}

The proof of Proposition~\ref{SZEGO} is now complete.

\section{Distribution of Zeros}\label{DZ}

In order to deduce Theorem
\ref{main} from  Proposition~\ref{SZEGO}, we need to relate
$\E_{|N P}(Z_s)$ to  the \szego kernel $\Pi_{|N P}$.  The idea is
essentially the same as in  \cite[Prop.~3.1]{SZ}, but neither the
statement nor the proof there cover the application we need. Hence we
provide a more general statement for the zeros of general linear systems on
compact \kahler manifolds (Proposition~\ref{EZ}).

\subsection{Expected distribution of zeros and \szego kernels}\label{EDZ} As in
\S\ref{s-szego}, we let
$L\to M$ be a holomorphic line bundle over a compact
\kahler manifold, and we consider a linear
subspace $\scal\subset H^0(M,L^N)$, which we  endow with a Hermitian inner
product. (The set of zero divisors $\{Z_s:s\in\scal\}\equiv \PP\scal$
is called a {\it linear system\/} on $M$.) We  let  $e_L$ be a local
holomorphic
frame over a trivializing chart $U$. We choose  an
orthonormal basis $\{S_j:1\le j\le k\}$ for
$\scal$ and we write $S_j = f_j e_L$ over $U$. Any  section $s\in\scal$ may
then be written as
$$s = \langle c, F \rangle e_L^{\otimes N}\;, \quad \mbox{where\ \ \ }
F=(f_1,\dots,f_k)\;,\quad\langle c,F
\rangle = \sum_{j = 1}^k c_j f_j\;.$$

The Kodaira map is given by:
\begin{equation}\label{kodaira}\Phi_\scal: M \dashrightarrow \CP^{k -
1},\quad
\Phi_\scal(z) = [f_1(z),
\dots, f_k(z) ]\;.\end{equation}
(We use the symbol `$\dashrightarrow$' in (\ref{kodaira}) to
indicate that the map $\Phi_\scal$ is a rational map; it is holomorphic if
and
only if
$\scal$ has no basepoints other than fixed components.) Under a change of
orthonormal basis, one gets a unitarily equivalent map.  The current of
integration over the zeros of $s =
\langle c,F
\rangle e_L^{\otimes N}$ is then given locally by the Poincar'e-Lelong
formula:
\begin{equation} Z_s =
\frac{\sqrt{-1}}{ \pi } \partial \bar{\partial}\log | \langle c,F
\rangle|\;.
\label{Zs} \end{equation} It is of course independent of
the choice of local frame $e_L$ and  basis
$\{S_j\}$.

We  write a section locally as $s = f e_L$ and we write its
Hermitian norm as $\|s(z)\|_h = a(z)^{-\half}|f(z)|$ where $a(z)
= \|e_L(z)\|_h^{-2}$. The {\it curvature form\/}
$c_1(L,h)$ of $L$ is given locally by $$c_1(L,h)=\frac{\sqrt{-1}}{2
\pi}\d\dbar\log a\;.$$

\begin{prop}\label{EZ}Let $(L,h)$ be a Hermitian line bundle on a compact complex
manifold
$M$. Let
$\scal$ be a subspace of $ H^0(M,L)$ of dimension $\ge 2$.  We
give
$\scal$ an
inner product and we let $\ga$ be the associated  Gaussian probability
measure.
Then the expected zero current of a random section
$s\in\scal$ is given by

\begin{eqnarray*}\E_\ga(Z_s)  &=&\frac{\sqrt{-1}}{2\pi}
\partial
\bar{\partial} \log \Pi_{\scal}(z, z)+c_1(L,h)\\
&=&\Phi^*_\scal\om_\FS + D_0\;,\end{eqnarray*}
where $D_0$ is the fixed component of the linear system $\PP\scal$.
\end{prop}

\begin{proof} The Gaussian measure $\ga$ is the standard Gaussian for
an inner product on $\scal$.  Let
$\{S_j\}$ be  an orthonormal basis for this inner product. As above, we
choose a
local nonvanishing section
$e_L$ of
$L$ over $U\subset M$, and we write
$$s=\sum_{j=1}^{k}c_jS_j=\langle c,F\rangle e_L\;,$$ where
$S_j=f_j
 e_L,\ F=(f_1,\dots,f_{k})$.  As in the proof of \cite{SZ},
Proposition 3.1, we then write $F(x)= |F(x)| u(x)$ so that
$|u| \equiv 1$ and
\begin{equation}\label{2terms}\log  | \langle c, F \rangle| = \log |F| +
\log  |
\langle c, u
\rangle|\;.\end{equation}

By (\ref{Zs}), we have
\begin{eqnarray*}\big(\E_\ga(Z_s),\phi\big)&=&\frac{\sqrt{-1}}{ \pi}
\int_{\C^{k}}
\left(\partial
\bar{\partial}\log  | \langle c,F \rangle| , \phi\right) d\gamma(c)\\
&=&\frac{\sqrt{-1}}{ \pi} \int_{\C^{k}}
\left(\partial
\bar{\partial}\log  |F| , \phi\right) d\gamma(c)
+\frac{\sqrt{-1}}{ \pi}\left( \int_{\C^{k}}
\partial\bar{\partial}\log  | \langle c, u
\rangle| d\gamma(c) , \phi\right)\;,
\end{eqnarray*}
for all test forms
$\phi\in\dcal^{m-1,m-1}(U)$.
Upon integration in $c$, the second term of (\ref{2terms}) becomes constant
in
$z$ and the derivatives kill it. The first term is independent of $c$ so we
may remove the Gaussian integral.  Thus
\begin{equation}\label{Egamma}\E_\ga(Z_s)=\frac{\sqrt{-1}}{2 \pi} \partial
\bar{\partial}\log  |F|^2  =\frac{\sqrt{-1}}{2 \pi} \partial
\bar{\partial}\left(\log \sum\|S_j\|^2+\log a\right)\;.\end{equation}
Recalling that $\Pi_{\scal}(z, z)=\sum\|S_j(z)\|^2$ and that
$c_1(L,h)=\frac{\sqrt{-1}}{2 \pi}\d\dbar\log a$, the first identity of
the proposition follows.

Before completing the proof of the second identity, we recall the definition
of
the pull-back of a form by a rational map
$X\buildrel{g}\over\dashrightarrow Y$: we let $\wt\Ga$ be the
desingularization
of the graph
$\Ga
\subset X\times Y$ of $g$ so that we have
$g\circ \pi_1 = \pi_2$, where $\pi_1,\pi_2$ are the projections from
$\wt\Ga$ to
$X,Y$ respectively.  Then for any smooth form $\eta$ on $Y$ we define the
pull-back current $g^*\eta:=\pi_{1*}\pi_2^*\eta$, which is easily seen to have
$\lcal^1$ coefficients with singularities on the indeterminacy locus $I_g$ of $g$.

We can suppose that $U$ is chosen so that there is a $\psi\in\ocal(U)$ with
$\mbox{Div}(\psi)=D_0$ on $U$.  Write $f_j=\psi g_j$.  Then by
(\ref{Egamma}) we
have
$$\E_\ga(Z_s)=\frac{\sqrt{-1}}{2 \pi} \partial
\bar{\partial}\log  |F|^2 =\frac{\sqrt{-1}}{2 \pi} \partial
\bar{\partial}\log |\psi|^2 +\frac{\sqrt{-1}}{2 \pi} \partial
\bar{\partial}\log\textstyle\sum|g_j|^2\;.$$
The first term equals $D_0$. By definition, the second term equals
$\Phi^*_\scal\om_\FS$ outside the indeterminacy locus
$I_{\Phi_\scal}$.  Since $I_{\Phi_\scal}$ has real codimension $\ge 4$ and
$\log
\sum|g_j|^2$ is (pluri-)subharmonic, the coefficients of
$\frac{i}{2\pi}\d\dbar
\log
\sum|g_j|^2$ are measures that do not give mass to $I_{\Phi_\scal}$; i.e.,
the
coefficients are in $\lcal^1(U)$.  Since the same is true for the
coefficients
of $\Phi^*_\scal\om_\FS$, the current equality extends to all of $U$.
\end{proof}

We note that the expected zero current $\E_\ga(Z_s)$ is a smooth form outside
the base point set of $\scal$. We shall use the following result on simultaneous
expected zeros (away from base points).

\begin{prop}\label{EZsimult}  Let $M$ be a compact complex manifold, and let
$(L_1,h_1),\dots ,(L_k,h_k)$ be Hermitian line bundles on $M$ ($1\le k\le \dim
M$). Suppose we are given subspaces $\scal_j\subset H^0(M,L_j)$ with inner
products $\langle,\rangle_j$ and let $\ga_j$ denote the associated Gaussian
probability measures on the $\scal_j$ ($1\le j\le k$).  Let $U$ be an open
subset of $M$ on which $\scal_j$ has no base points for all $j$. Then the
expected simultaneous-zero current of random sections
$s_1\in\scal_1,\;\dots,\;s_k\in\scal_k$ is given over $U$ by
$$\E_{\ga_1\times\cdots\times\ga_k}(Z_{s_1,\dots,s_k}) = \bigwedge_{j=1}^k
\E_{\ga_j}(Z_{s_j}) = \bigwedge_{j=1}^k
\left[\frac{\sqrt{-1}}{2\pi}
\partial
\bar{\partial} \log \Pi_{\scal_j}(z, z)+c_1(L_j,h_j)
\right] =\bigwedge_{j=1}^k \Phi^*_{\scal_j}\om_\FS\;.$$
\end{prop}

\begin{proof} We first note that the expect current is well defined, since for
all choices of the $s_j$ the total mass of the zero current $Z_{s_1,\dots,s_k}$
equals the Chern number $c_1(L_1)\cdots c_1(L_k)$, and hence for each test
form $\phi\in\dcal^{m-k,m-k}(U)$, the function $(Z_{s_1,\dots,s_k},\phi)$ is in
$\lcal^\infty(\scal_1\times\cdots\times \scal_k)$.

Next we note the following functorial property of the expected zero current:
Suppose that $\rho:Y\to M$ is a holomorphic mapping of complex manifolds and
let $\rho^*\ga_j$ be the pulled-back Gaussian measure on $\rho^*\scal_j\subset
H^0(Y,\rho^*L)$.  Then
\begin{equation}\label{functorial}\E_{\rho^*\ga_j}(Z_{t_j}) = \rho^*
\E_{\ga_j}(Z_{s_j})
 \quad\mbox{on }\ \rho\inv(U) \qquad (t_j\in \rho^*\scal_j)\;.\end{equation}
Indeed, (\ref{functorial}) follows immediately from Proposition~\ref{EZ} and
its proof.  As a special case, suppose that $Y$ is a submanifold of $M$ and
$\rho $ is the inclusion.  Then (\ref{functorial}) becomes
\begin{equation}\label{Y} \E_{\ga_j}([Y]\wedge Z_{s_j})=[Y]\wedge
\E_{\ga_j}(Z_{s_j})
\quad\mbox{on }\ U\;.\end{equation} (In fact, (\ref{Y}) holds for any
subvariety $Y$, but we don't need this fact here.)

We now verify the current identity by induction on $k$.  For $k=1$, this is
Proposition \ref{EZ}, so assume that $k>1$ and the proposition has been
verified for $k-1$.  Consider $Y= |Z_{s_1,\dots,s_{k-1}}|$. By Bertini's
Theorem, we know that $Y$ is smooth (and of codimension $k-1$) for almost all
$s_1,\dots,s_{k-1}$.  Therefore by (\ref{Y}),
$$ \int  Z_{s_1,\dots,{s_k}}\,d\ga_k(s_k) =
\int  Z_{s_1,\dots,s_{k-1}}\wedge Z_{s_k}\,d\ga_k(s_k)
=Z_{s_1,\dots,{s_{k-1}}}\wedge \E_{\ga_k}(Z_{s_k}) \  \mbox{for a.a. }
s_1,\dots,s_{k-1}.$$
Integrating over $s_1,\dots,s_{k-1}$ and using the inductive hypothesis, we
obtain
$$\E_{\ga_1\times\cdots\times\ga_k}(Z_{s_1,\dots,s_k}) =
\E_{\ga_1\times\cdots\times\ga_{k-1}}(Z_{s_1,\dots,s_{k-1}})\wedge
\E_{\ga_k}(Z_{s_k}) =\E_{\ga_1}(Z_{s_1}) \wedge\cdots\wedge
\E_{\ga_k}(Z_{s_k})\;.$$
The other equalities follow from Proposition \ref{EZ}.
\end{proof}

If we consider the case where $M = \CP^m$, we obtain a formula for the expected
simultaneous zero current for random polynomials with given Newton polytopes:

\begin{cor} \label{E-cond}The expected zero current of  $k$ independent random
polynomials
$f_j\in H^0(\CP^m,\ocal(p_j),P_j)$, $1\le j\le k$, is given
over $\C^{*m}$ by
$$\E_{|P_1,\dots,P_k}(Z_{f_1,\dots,f_k})  = \bigwedge_{j=1}^k\left(
\frac{\sqrt{-1}}{2\pi}
\partial
\bar{\partial} \log \Pi_{| P_j}(z, z)+p_j
\omega_{\FS}^{\CP^m}\right)\qquad \mbox{\rm on}\ \ \C^{*m}\;.
$$
\end{cor}

\begin{proof} Since  $H^0(\CP^m,\ocal(p_j),P_j)$ has no base points
in $\C^{*m}$, the conclusion follows from Proposition \ref{EZsimult} and the fact
that
$c_1(\ocal(p_j),h_\FS)=p_j\om_\FS$.
\end{proof}

\begin{rem} By Proposition \ref{EZ}, the identity of currents
\begin{equation}\label{E-cond1}\E_{|P}(Z_f)  =
\frac{\sqrt{-1}}{2\pi}
\partial
\bar{\partial} \log \Pi_{| P}(z, z)+p
\omega_{\FS}^{\CP^m}\end{equation}
holds on all of $\CP^m$.

We say that  {\it $p$ is minimal for $P$} if
\begin{equation} \min _{\alpha
\in P}\al_j=0 \quad (1\le j\le n)\,,\quad \max_{\alpha
\in P} |\alpha| = \de\;.\end{equation} (Otherwise, we can choose $\de$ to be
smaller, after translating the polytope if necessary.)
The hypothesis that $p$ is minimal for $P$ guarantees that
$H^0(\CP^m,\ocal(p),P)$ has no fixed divisors.  In this case, (\ref{E-cond1})
holds as an identity of $(1,1)$-forms with $\lcal^1_{\rm {loc}}$ coefficients.
\end{rem}

\subsection{Asymptotics of zeros}\label{s-proofs} In this section, we prove 
Theorems \ref{main} and
\ref{simultaneous} and give some additional asymptotic formulas for zeros.
Theorem~\ref{probK} is a special case  of Corollaries
\ref{zerovolumes} and \ref{subtler} of Theorem~\ref{simultaneous}. 

To prove Theorem
\ref{main}, we let
$u_N,\ u_\infty$ be as in the proof of Proposition~\ref{convergence}. By Corollary
\ref{E-cond} (recalling that $\om_\FS= \frac{\sqrt{-1}}{2\pi}\ddbar \log
(1+\|z\|^2)$), we have
\begin{equation}\label{ddbaruN}\frac{\sqrt{-1}}{2\pi}\d\dbar
u_N=\frac{1}{N}\E_{|NP}(Z_f) \qquad \mbox{on }\ \C^{*m} \;,\end{equation} and
hence $u_N$ is plurisubharmonic.  By (\ref{u0}) and
(\ref{ddbaruN}), we conclude that
\begin{eqnarray} \frac{1}{N}\E_{|NP}(Z_f) \to  \frac{\sqrt{-1}}{2\pi}
\ddbar u_\infty=\de\om_\FS -
\frac{\sqrt{-1}}{2\pi}
\ddbar b\;,\label{psi0}\end{eqnarray} where differentiation is in
the distribution sense.

We now show that the current $\ddbar b\in \dcal'{}^{1,1}(\C^{*m})$ is given by a
$(1,1)$-form with piecewise smooth coefficients; in fact $\ddbar b$ is
$\ccal^\infty$ on each of the regions
${\rcal_F}$ given by (\ref{RF}). (Equivalently, the Radon measure
$\ddbar b\wedge \om_\FS$ does not charge the set of transition points, and hence
formula (\ref{psi0}) can be interpreted as differentiation in the ordinary
sense on the regions
$\rcal_F$.) By Proposition~\ref{SZEGO},
$b\in\ccal^1(\C^{*m})$; i.e., if $z^0\in\d{\rcal_F}\cap
\d\rcal_{F'}$, then the values
of $d b(z^0)$ computed in the two regions $\overline{\rcal_F}$ and
$\overline{\rcal_{F'}}$ agree.  Then for a test form $\phi\in\dcal^{m-1,m-1}(\C^{*m})$,
we have
\begin{equation}\label{ddbarb1}(\d\dbar b, \phi) = \sum_F\int_{\rcal_F}b\d\dbar
\phi =\sum_F\int_{\rcal_F^\circ}
\ddbar b\wedge \phi - \sum_F\int_{\d\rcal_F}(\dbar b\wedge \phi
+b\wedge\d\phi)\;.\end{equation}
Note that
$\d\rcal_F$ consists of
$\ccal^\infty$ real hypersurfaces (consisting of those points of
$\d{\rcal_F}$ that are contained in the boundary of only one other region
$\rcal_{F'}$) together with submanifolds of real codimension
$\ge 2$, and hence  Stokes' Theorem applies (see e.g. \cite[4.2.14]{Fe}).
Since $b$ and $\dbar b$ are continuous on $\C^{*m}$, the boundary integral
terms  in (\ref{ddbarb1}) cancel out and we obtain
\begin{equation}\label{ddbarb2}(\d\dbar b, \phi)  =\sum_F\int_{\rcal_F^\circ}
\ddbar b\wedge \phi  =\int_{\C^{*m}\sm E}\ddbar b\wedge \phi\;,\end{equation}
where
$E:=\bigcup_F \d\rcal_F$ is the set of transition points.
Formula (\ref{ddbarb2}) says that  the current $\ddbar b$ is  a $(1,1)$-form with
piecewise smooth coefficients obtained by differentiating $b$ on the regions
$\rcal_F^\circ$.

We now let
\begin{equation}\label{psi} \psi_P=\frac{\sqrt{-1}}{2\pi}
\ddbar u_\infty=\de\om_\FS -
\frac{\sqrt{-1}}{2\pi}
\ddbar b\end{equation}
on each of the regions $\rcal_F$. By
(\ref{ddbarb2}), the current $\psi_P$ is also a piecewise smooth $(1,1)$-form; by
(\ref{psi0})
\begin{equation}\label{weaklim0}N\inv\E_{|NP}(Z_f)
\to \psi_P\qquad \mbox {(weakly)}\;.\end{equation}

We shall show $\lcal^1_{\rm loc}$ convergence of (\ref{weaklim0}) when we
prove Theorem \ref{simultaneous} below.  Continuing with the proof of Theorem
\ref{main}, we observe that (ii) is an immediate consequence of (\ref{psi}),
since $b=0$ on the classically allowed region.

To  verify (iii), we again use the log coordinates $\rho_j +i\theta_j=\log z_j$, so that
\begin{equation}\label{uinfty}u_\infty=p\log\left(1+\sum
e^{2\rho_j}\right)-b(e^\rho)\;.\end{equation} Since
$u_\infty$ depends only on $(\rho_1,\dots,\rho_m)$, we have
\begin{equation}\label{d2rho}\psi_P=\frac{\sqrt{-1}}{2\pi}\ddbar u_\infty =
\frac{\sqrt{-1}}{8\pi}\sum_{jk}
\frac{\d^2 u_\infty}{\d\rho_j\d\rho_k} d\zeta_j\wedge d\bar\zeta_k\ge
0\;.\end{equation}  Thus we must show that the Hessian of $u_\infty(\rho)$ has
rank $r$ at points $z^0\in\rcal^\circ_F$, where $r=\dim F$. {From}
(\ref{Db1})--(\ref{Dlog}) and (\ref{uinfty}), we obtain
\begin{eqnarray} D_\rho u_\infty&=&2p\,\mu_\Si\circ\xi\;,
 \label{dbaru}\end{eqnarray} where $\xi(z)=e^{\tau_z/2}\cdot z$ as before.
Hence $\frac{1}{2p} D^2_\rho u_\infty$ equals the Jacobian of
$\mu_\Si\circ\xi$.  Since $\mu_\Si
\circ
\xi (z)=\frac{1}{\de}\mu_P(w_z)$, we easily see that
\begin{equation}\label{submersion}\mu_\Si\circ\xi:\rcal^\circ_F
\to F\;,\end{equation} so rank$\,D(\mu_\Si\circ\xi) \le r$.  In fact,
(\ref{submersion}) is a submersion, so that the rank equals $r$.  To see that
it is a submersion, we recall from the proof of Lemma
\ref{claim} that the `lipeomorphism' $\Phi\inv:\Si^\circ\to
\ncal^\circ$ restricts to a diffeomorphism
$$\Si^\circ\supset \mu_\Si(\rcal^\circ_F) \buildrel{\approx}\over \to F\times
C_F\subset \ncal^\circ\;, \qquad \mu_\Si(z)\mapsto
\left(\textstyle\frac{1}{\de}\mu_P(w_z),-\tau_z\right)
=(\mu_\Si\circ\xi(z),-\tau_z)\;,$$ and hence (\ref{submersion}) is a
submersion, completing the proof of (iii).

We now prove Theorem \ref{simultaneous}; the case $k=1$ of the theorem
will then yield part  (i) of Theorem~
\ref{main}.

Let $P_1,\dots, P_k$ be Delzant polytopes, as in Theorem \ref{simultaneous}.
We first show that
\begin{equation}\label{weaklim} N^{-k}\E_{|NP}(Z_{f_1,\dots,f_k})
\to \psi_k:=\psi_{P_1}\wedge \cdots\wedge \psi_{P_k}\qquad \mbox
{(weakly)}.\end{equation}
For $1\le j\le k$, we let $$u_N^j=\frac{1}{N}\log \Pi_{|NP_j}(z,z) +\de_j
\log(1+\|z\|^2)\;,$$ so that, recalling (\ref{u0}),
$$u_N^j(z)\to u_\infty^j(z):= \de_j
\log(1+\|z\|^2) -b(z)= \de_j\log (1+\|e^{\tau_z^j/2}\cdot z\|^2)
+\log H_{P_j}(\tau_z^j,w_z^j)\;,$$ where $\tau_z^j,w_z^j$ are as in Lemma
\ref{claim2-3} with $P=P_j$.  By (\ref{psi}), $$\frac{\sqrt{-1}}{2\pi}
\ddbar u_\infty^j = \psi_{P_j}\;.$$

Recalling Corollary \ref{E-cond}, we introduce the $(k,k)$-forms
\begin{equation}\label{w1}\kappa_N:=\left(\frac{\sqrt{-1}}{2\pi}
\ddbar u_N^1\right) \wedge \cdots \wedge \left(\frac{\sqrt{-1}}{2\pi}
\ddbar u_N^k\right)
=N^{-k}\E_{|NP_1,\dots,NP_k}(Z_{f_1,\dots,f_k})\;.\end{equation}
Since $u_N^j\to u_\infty^j$ locally uniformly, it follows from the Bedford-Taylor
Theorem \cite{BT,Kl} that
\begin{equation}\label{w2}\kappa_N \to \left(\frac{\sqrt{-1}}{2\pi}
\ddbar u_\infty^1\right) \wedge \cdots \wedge \left(\frac{\sqrt{-1}}{2\pi}
\ddbar u_\infty^k\right) \qquad
\mbox {(weakly)}.\end{equation}  In fact, the limit current in (\ref{w2})
is absolutely continuous, and hence
is equal to the locally bounded (piecewise smooth)
$(k,k)$-form
$\psi_k$.  To see this, we note that by the Bedford-Taylor Theorem,
$$\bigwedge_{j=1}^k(\psi_{P_j} *\phi_\epsilon)=\bigwedge_{j=1}^k
\left[\frac{\sqrt{-1}}{2\pi}\ddbar (u_\infty^j *
\phi_\epsilon)\right] \to \bigwedge_{j=1}^k\left(\frac{\sqrt{-1}}{2\pi}\ddbar
u_\infty^j\right)\;,$$ where $\phi_\epsilon$ denotes an
approximate identity. Since the currents $\psi_{P_j}$ have coefficients in
$\lcal^\infty_{\rm{loc}}$, the forms $\bigwedge_{j=1}^k(\psi_P
*\rho_\epsilon)$ have locally uniformly bounded coefficients; absolute
continuity of the limit current follows. The weak limit (\ref{weaklim}) now
follows from (\ref{w1})--(\ref{w2}).

To complete the proof of Theorem \ref{simultaneous}, we must show that
$\kappa_N \to \psi_k$ in
$\lcal^1(K)$ for all compact $K\subset\C^{*m}$.  Let $\ep>0$ be arbitrary,
and choose a nonnegative function $\eta\in \dcal(\C^{*m})$ such that $\eta\equiv
1$ on a neighborhood of $(E_1\cup\cdots\cup E_k)\cap K$, where $E_j$ is the
set of transition points for $P_j$, and
$$\int_{\C^{*m}}
\psi_k\wedge
\eta\om_\FS^{m-k}<\ep\;.$$
By (\ref{weaklim})--(\ref{w1}),
$$\int_{\C^{*m}}
\kappa_N\wedge
\eta\om_\FS^{m-k}<2\ep\qquad \mbox{for }\ N\gg 0\;.$$
(Note that the above integrands are nonnegative.)
Therefore,
$$\left\|\kappa_N - \psi_k\right\|_{\lcal^1(K)}  \le
\left\|(1-\eta)(\kappa_N - \psi_k)\right\|_{\lcal^1(K)} +
\left\|\eta\kappa_N\right\|_{\lcal^1(K)}+
\left\|\eta\psi_k\right\|_{\lcal^1(K)}\;.
$$ Since $\kappa_N$ is a positive $(k,k)$-form,
$$\left\|\eta\kappa_N\right\|_{\lcal^1(K)}=
\int_K\kappa_N\wedge
\eta\om_\FS^{m-k}<2\ep \;.$$ Similarly,
$\left\|\eta\psi_k\right\|_{\lcal^1(K)}=\int_K
\psi_k\wedge
\eta\om_\FS^{m-k}<\ep$. Since
$$(1-\eta)(\kappa_N - \psi_k)\to 0\qquad \mbox{uniformly on }\ K\;,$$
it follows that
$$\left\|\kappa_N - \psi_k\right\|_{\lcal^1(K)} \to 0\;,$$
completing the proof of Theorems \ref{main} and \ref{simultaneous}.\qed

\bigskip
We note that the conclusion of Theorem \ref{simultaneous} can be replaced with an
asymptotic expansion away from transition points:

\begin{theo}\label{expansion} Let $P_1,\dots,P_k$ be Delzant polytopes.  Let
$U$ be a relatively compact domain in $\C^{*m}$ such that $\overline U$ does not
contain transition points for any of the $P_j$. Then we have a complete
asymptotic expansion of the form
$$\frac{1}{N^{k}}\E_{|NP_1,\dots,NP_k} (Z_{f_1, \dots, f_k}) \sim
\psi_{P_1}\wedge\cdots\wedge\psi_{P_k} +\frac{\phi_1}{N} + \cdots +
\frac{\phi_n}{ N^n} +\cdots \qquad \mbox
{on } \ U\;,$$ with uniform
$\ccal^\infty$ remainder estimates, where the $\phi_j$ are smooth
$(k,k)$-forms on $U$.
\end{theo}

\begin{proof} By (\ref{uN}) and (\ref{w1}) we have
$$\frac{1}{N^{k}}\E_{|NP_1,\dots,NP_k} (Z_{f_1, \dots, f_k})=
\bigwedge_{j=1}^k\frac{\sqrt{-1}}{2\pi}
\ddbar \left[\frac{1}{N}\log \Pi_{|NP_j}(z,z) +p_j
\log(1+\|z\|^2)\right]\;.$$ The conclusion now follows from the
asymptotic expansion of Proposition~\ref{SZEGO}.\end{proof}

The proof of Theorem \ref{main}(iii) gives us some more information about the
expected zero current in the classically forbidden region, which we state in
the theorem below.  We first define the {\it complexified normal cones\/}
$$\wt C_F = \{\tau+i\theta: \tau \in C_F,\ \theta\in T_F^\perp\}\;.$$  (Recall
that $C_F\subset T_F^\perp$.)
 We note that $\wt C_F$ is a semi-group, which
 acts on
$\rcal_F$ by the rule $\eta(z )= e^\eta \cdot z$; we call this action the
`(joint) normal flow.'  The (maximal) orbits of the normal flow are of the form
$\wt C_F\cdot z^0=\{e^{\eta}\cdot z^0:\eta\in\wt C_F\}$, where  $z^0\in
\mu_\Si\inv(F)$. We note that the orbit $\wt C_F\cdot z^0$ is a complex
$(m-r)$-dimensional submanifold (with boundary) of $\C^{*m}$.  (Indeed,
$(\wt C_F\cdot z^0)\cap \rcal^\circ_F$ is a submanifold without boundary in
$\rcal^\circ_F$.)

\begin{theo}\label{more} Let $P$ be a Delzant polytope and let $\psi_P$ be
the limit expected zero current of Theorem \ref{main}. Then
$\psi_P$ vanishes along the
orbits of the normal flow.
\end{theo}

\begin{proof} Let $$O:=\wt C_F\cdot z^0=\{e^{\tau+i\theta}\cdot z^0: \tau
\in C_F,\
\theta\in T_F^\perp\}$$ be a maximal orbit of the normal flow, where
$\mu_\Si(z^0)\in F$. For $z= e^{\tau+i\theta}\cdot z^0\in O$, we have
$$\mu_\Si\circ\xi (z)= \frac{1}{\de}\mu_P(w_z)=
\frac{1}{\de}\mu_P(w_{e^{\tau}\cdot z^0})=\frac{1}{\de}\mu_P(w_{z^0})
=\mu_\Si(z^0)$$ and
hence $\mu_\Si\circ\xi$ is constant on $O$. It then follows from
(\ref{d2rho})--(\ref{dbaru}) that $\psi_P|_O=0$.
\end{proof}

\begin{rem}
Theorem \ref{more} tells us about the behavior of tangent
directions to typical zero sets $\left|Z_f\right|$ of polynomials
$f\in H^0(\CP^m,\ocal(N\de),NP)$ as $N\to\infty$. Roughly speaking,  the
tangent spaces of
$\left|Z_f\right|$ are highly likely to be close to containing the tangent
spaces of the orbits of the normal flow, for large $N$.  E.g., let us consider a
region $\rcal_F^\circ$, where
$F$ is an edge (i.e., $\dim F =1$).    Let
$T^{\rm flow}_z$ denote the holomorphic tangent space to the normal flow through
a point $z\in\rcal_F^\circ$.  Then $\dim_\C T^{\rm flow}_z = \dim_\C Z_f=m-1$,
and for any compact set
$K\subset
\rcal^\circ_F$, we have
\begin{equation}\label{limtangents}  \frac{1}{N}\E_{|NP}\left(
\int_{Z_f\cap K}\mbox{dist}(T_{Z_f,z}^{1,0},\, T^{\rm flow}_z)^2
\,d\vol_
{Z_f}(z)\right)\longrightarrow 0\;,
\end{equation}
where `dist' means the distance in the projective space of complex
$(m-1)$-dimensional subspaces of $T^{1,0}_{CP^m,z}$.

To verify (\ref{limtangents}), we let $V$ be a compactly
supported vector field on
$\rcal^\circ_F$ such that for all $z\in\rcal^\circ_F$, $V_z$
is a $(1,0)$-vector tangent to the normal flow through $z$.  Let
$\phi=*(iV\wedge
\bar V)\in
\dcal^{m-1,m-1}(\rcal^\circ_F)$ be given by $$\la\wedge\phi_z = (\la, iV_z\wedge
\bar V_z)\, \textstyle\frac{1}{m!}\om^m_\FS\;,\qquad \la\in T^{*1,1}_z\C^m\;,\
z\in
\rcal^\circ_F\;.$$
By Theorem \ref{more},  we have $$(\psi_P, \phi) = \int \psi_P\wedge\phi =
\int (\psi_P , iV_z\wedge
\bar V_z)\, \textstyle\frac{1}{m!}\om^m_\FS=0\;.$$ We let $\tcal_f\in
\ccal^\infty(Z_f^{\rm reg}, T^{m-1,m-1}_{Z_f})$ denote the dual to the volume
form on $Z_f^{\rm reg}$, so that
\begin{equation}\label{tangents} \int_{Z_f}\|\tcal_f
\wedge V\wedge \bar V\|\,d\vol_
{Z_f}=\int_{Z_f}(\tcal_f
, \phi)\,d\vol_
{Z_f} = (Z_f,\phi)\;.
\end{equation} Hence by Theorem \ref{main}(i),
\begin{equation}\label{limtangents1}\frac{1}{NP}\E_{|NP}\left(
\int_{Z_f}\|\tcal_f
\wedge V\wedge \bar V\|\,d\vol_
{Z_f}\right)=\frac{1}{N}\E_{|NP}(Z_f,\phi)\to(\psi_P,\phi)=0\;,
\end{equation}
which yields (\ref{limtangents}). \end{rem}

\medskip
As mentioned in the introduction, we can apply Proposition \ref{EZsimult} (as
in \cite{SZ} for the case $k=1$) to the asymptotic expansion of the \szego kernel
in
\cite{Z} to immediately obtain the following asymptotics using toric norms:

\begin{prop} \label{toriczero} Let $(M_P, L_P)$ be
as above, and  let $I_j,\theta_j$ be the action-angle variables of
the  moment map $\mu_P: M_P \to \R^m$ of the $\T $-action on
$M_P$. Then we have:
$$ \frac{1}{N^k } E_{\ga_N^{M_P}}  (Z_{f_1, \dots, f_k}) =\om_P^k +
O\left(\frac{1}{N}\right)=
\left(\sum_{j=1}^m dI_j \wedge d\theta_j\right)^k +
O\left(\frac{1}{N}\right)\;.$$
\end{prop}

\begin{proof} Let $\omega_P=\sum_{j=1}^m dI_j \wedge d\theta_j$ be the natural
\kahler form on $M_P$. We recall the asymptotic formula for the \szego kernel
in \cite{Z}:
$$\frac{1}{N}\Phi_{L_P^N} ^* \omega_{\FS} =
\omega_P + O(1/N)\;. $$ By Proposition \ref{EZsimult} with $M=M_P, \ L= L_P^N$,
we then have
$$\E_{M_P} (Z_{s_1,\dots,s_k})= (\E_{M_P} (Z_{s}))^k=\Phi_{L_P^N} ^*
\omega_{\FS}^k =N^k\left[ \om_P^k+O(1/N)\right]\;. $$ \end{proof}

\subsection{Amoebas in the plane}\label{AM}  The term `amoeba' was introduced by Gelfand,
Kapranov and  Zelevinsky \cite{GKZ} to refer to the image under the moment map of a
zero set $Z_{f_1,\dots,f_k}$ of polynomials, and have been studied in various contexts
(see \cite{FPT,GKZ,M1,PR} and the references in the survey article by Mikhalkin
\cite{M}).  The image of a zero set under the moment map
$\mu_\Si$ is called a {\it compact amoeba\/}, while the image under the map
$$\Log:\C^{*m}\to\R^m, \qquad (z_1,\dots,z_m)\mapsto (\log|z_1|,\dots, \log|z_m|)$$ is
 a {\it noncompact amoeba\/}, or simply an {\it amoeba\/}.  Note that $\Log$ is the
moment map for the
$\T$ action with respect to the Euclidean symplectic form $\sum dx_j\wedge dy_j$, and
$\Log= \half \lcal\circ \mu_\Si$, where $\lcal:\Si^\circ\approx \R^m$ is the
diffeomorphism given by (\ref{lcal}).

To illustrate what our statistical results can say about amoebas, we
consider zero sets in $\C^{*2}$. An amoeba in $\R^2$ is the image of a plane algebraic
curve under the Euclidean moment map Log. An example of an amoeba of the form
$\Log(Z_f)$, where $f$ is a quartic polynomial in two variables with (full)
Newton polytope
$4\Si$, is given in the illustration from \cite{Th} reproduced in
Figure~\ref{fig-amoeba} below.\footnote{The authors would like to thank T. Theobald for
giving us permission to use this figure from his paper.}

\begin{figure}[htb]
\centerline{\includegraphics*[bb=194 274 417 469]{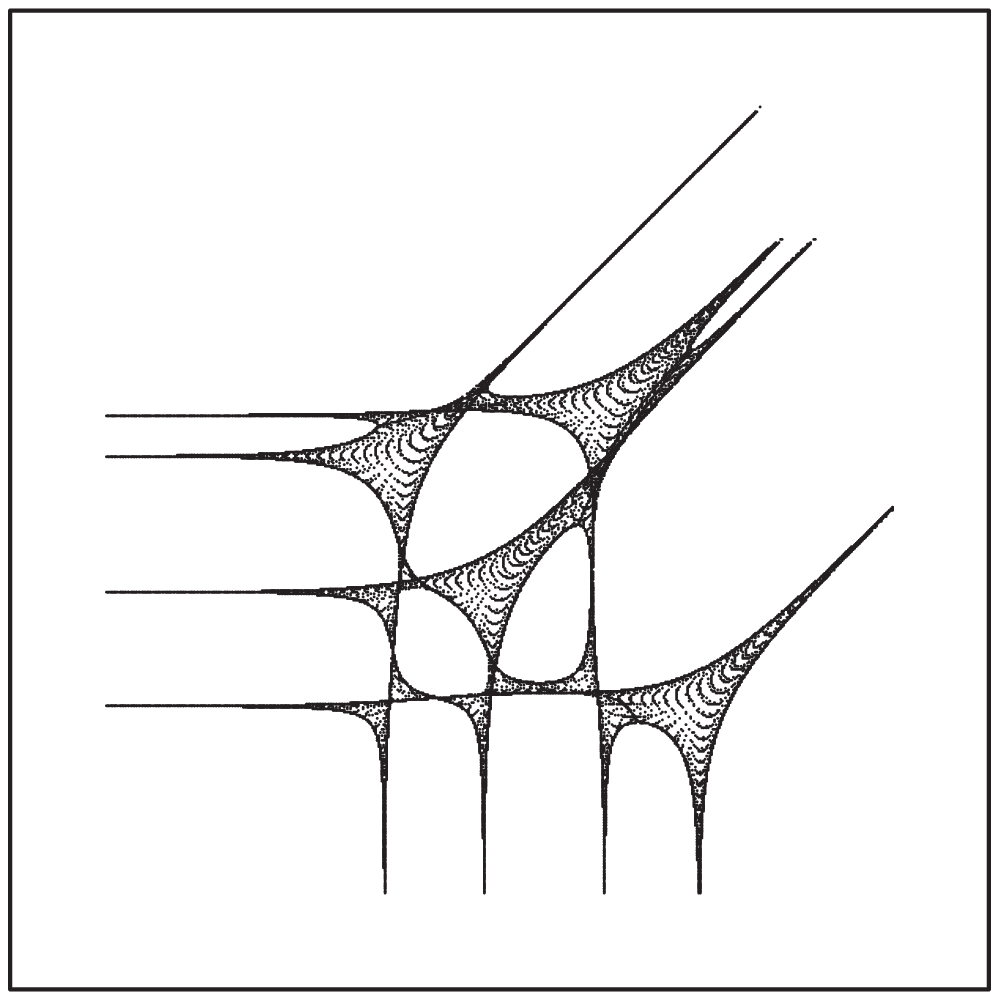}}
\caption{An amoeba with polytope $4\Si$}
\label{fig-amoeba}\end{figure}

One notices that this amoeba contains 12 `tentacles'. By definition, a
tentacle on a compact amoeba
$A$ is a connected component of a small neighborhood in $A$ of
$A\cap\d\Si$; the tentacles of a noncompact amoeba correspond to those of the
compact amoeba under the diffeomorphism  $\Si^\circ\approx \R^m$.

There is a natural injective map from the set of connected components
(which are convex sets) of the complement of a noncompact amoeba $A$ to the
set $P\cap\Z^2$ of lattice points of the polytope, and that there are amoebas
for which each lattice point is assigned to a component of the complement.
(This fact is also valid in higher dimensions; see
\cite{FPT, M1}.)  For a generic 2-dimensional
noncompact amoeba with polytope
$P$, each lattice point in $\d P$ corresponds to a distinct  unbounded
component of
$\R^2\sm A$, and adjacent lattice points correspond to adjacent unbounded components.
(The correspondence is given by \cite[\S 3.1]{M1} or \cite{FPT}.) Hence the number of
tentacles of
$A$ equals the number of points of
$\d P\cap\Z$; this number is called the {\it length\/} of $\d P$. (For example,
the 12 tentacles in Figure \ref{fig-amoeba} correspond to the 12 lattice points
of $\d(4\Si)$.) Each tentacle corresponds to a segment connecting 2 adjacent
lattice points on
$\d P$.

We can decompose $\d P$ into two pieces: $\d^\circ P=\d P\cap p\Si^\circ$ and
$\d^e P=P\cap \d (p\Si)$. Tentacles corresponding
to segments of $\d^\circ P$ end (in the compact picture $\Si$) at a vertex of
$\Si$, and tentacles corresponding to segments of $\d^e P$ are free to end
anywhere on the face of $\Si$ containing the segment.  We call the latter {\it
free tentacles\/}, and we say that a free tentacle is a {\it classically allowed
tentacle\/} if its end is in the classically allowed region $\acal_P$. For an
amoeba $A$, we let
$\nu_{\rm AT}(A)$ denote the number of classically allowed tentacles of $A$.
It is clear from the above that
$$\nu_{\rm AT}(A) \le \#\{\mbox{free tentacles}\}=\mbox{Length}(\d^e P)$$ and
that this bound can be attained for any polytope $P$.  Here, `Length' means the
length in the above sense; i.e., the diagonal face of $p\Si$ is scaled to have
length $p$. As a consequence of Theorem
\ref{main} (for
$m=1$), we conclude that this maximum is asymptotically the average:

\begin{cor}\label{amoeba} For a Delzant polytope $P$, we have
$$\frac{1}{N}\E_{|NP}\Big(\nu_{\rm AT}\big(\Log
(Z_f)\big)\,\Big) \to \mbox{\rm Length}(\d^e P)\;.$$
\end{cor}

\begin{proof} Let $F_1,F_2,F_3$ denote the facets of $p\Si$ and apply Theorem
\ref{main}(ii) to the 1-dimensional polytopes $P\cap \bar F_j$,
$j=1,2,3$.\end{proof}

\section{Examples}

We now compute $b(z)$ and $\psi_P$ for some  examples of toric surfaces
($m=2$).
Recalling (\ref{b}) and (\ref{u0}),
$$b (z) =  q(z)\cdot \tau_z + \de
\log\left(\frac{1+\|z\|^2}{1+\|e^{\tau_z/2}\cdot z\|^2}\right)\;,\quad
q(z)=\de\,\mu_\Si(e^{\tau_z/2}\cdot z) \in\d P\;,$$
$$u_\infty = -b(z) +p
\log (1+\|z\|^2)=p\log (1 +\|e^{\tau_z/2}
\cdot
z\|^2) - \langle q(z),
\tau_z\rangle\;.$$

Hence by (\ref{psi}),\begin{equation}\label{psi-formula} \psi_P =
\frac{\sqrt{-1}}{2\pi }\partial\bar\partial u_\infty =
\frac{\sqrt{-1}}{2\pi}\partial\bar\partial\left[p \log (1 +\|e^{\tau_z/2}
\cdot
z\|^2) - \langle q(z),
\tau_z\rangle\right]\;.\end{equation}

\subsection{The square}\label{example1} For our first example, we let
$P$ be the unit
square with vertices $\{(0,0),\,(1,0),\,(0,1),\,(1,1)\,\}$
so that
$p = 2$ and
$M_P =
\CP^1\times
\CP^1$.
Recalling that
$$\mu_\Si(z_1,z_2)=\left(\frac{|z_1|^2}{1+|z_1|^2 +|z_2|^2},
\frac{|z_2|^2}{1+|z_1|^2 +|z_2|^2}\right)\;,$$ we see that the classically
allowed region is given by
$$|z_1|^2 -1 < |z_2|^2< |z_1|^2 +1\;,$$ as illustrated in Figure~\ref{allowed-square}
in the introduction. The forbidden region consists of two subregions:
$$\begin{array}{ll}\rcal_F=\{(z_1,z_2):|z_2|^2 \ge |z_1|^2 +1\}\;,\qquad &
 F=\{(x_1,\half):0\le x_1
\le
\half\}\;,\\[8pt]\rcal_{F^*}=\{(z_1,z_2):|z_2|^2 \le |z_1|^2 -1\}\;,\qquad &
 F^*=\{(\half,x_2):0\le x_2
\le
\half\}
\;.\end{array}$$

Suppose that $z$ is a point in the upper forbidden region $\rcal_F$.
Write $\tau_z=(\tau_1,\tau_2)$; then $\tau_1=0$ since $\tau\perp T_F$.
Let $$ q(z)= 2 \mu_\Si (e^{\tau_z/2}\cdot z) = (a,1)\in 2F\subset \d
P\;.$$
  Writing
$$ |z_1|^2 =s_1,\quad |z_2|^2 = s_2\;,\quad |e^{\tau_2/2}
z_2|^2=e^{\tau_2} s_2 =\tilde s_2\;,$$
we have $$\frac{ s_1}{1+s_1+\tilde s_2} =\frac{a}{2}\;,\qquad \frac{\tilde
s_2}{1+s_1 +\tilde s_2}=\frac{1}{2}\;.$$
Therefore
$$\begin{array}{c}\di s_1 =\frac{a}{1-a},\quad \tilde s_2
=\frac{1}{1-a}= \frac{s_1}{a},\quad a= \frac{s_1}{1+s_1} =
\frac{|z_1|^2}{1+|z_1|^2},\\[12pt]\di e^{\tau_2} =\tilde s_2/s_2
=\frac{s_1}{as_2} =
\frac{1+|z_1|^2}{|z_2|^2}\;.\end{array}$$
 We have $$
\log (1+\|e^{\tau_z/2}\cdot z\|^2) =
 \log\left (1+|z_1|^2 + \frac{1+|z_1|^2}{|z_2|^2}
|z_2|^2\right)= \log (1 +|z_1|^2) +\log 2\;,$$
$$\langle q(z), \tau_z\rangle =\left\langle\left(\frac{|z_1|^2}{1+|z_1|^2}
,1\right), \left(0,\log \frac{1+|z_1|^2}{|z_2|^2}\right)\right\rangle
= \log (1+|z_1|^2) -\log |z_2|^2\;.$$
Therefore $$u_\infty= \log|z_2|^2 + \log
(1+|z_1|^2) +\log 4\;.$$

We conclude that $$  \psi_P=
\left\{\begin{array}{ll}
\frac{\sqrt{-1}}{2\pi} \partial\bar\partial \log
(1+|z_1|^2)\qquad &
\mbox{for }\ |z_1|^2 +1\leq |z_2|^2\\[14pt]
2\om_\FS = \frac{\sqrt{ -1}}{\pi} \partial\bar\partial \log
(1+|z_1|^2+|z_2|^2)\qquad &
\mbox{for }\  |z_1|^2-1\le |z_2|^2\leq |z_1|^2+1\\[14pt]
 \frac{\sqrt{ -1}}{2\pi} \partial\bar\partial \log
(1+|z_2|^2)\qquad &
\mbox{for }\ |z_2|^2\leq |z_1|^2-1\end{array}\right.$$
(where the third case is by symmetry).  Note that $\psi_P$ has constant rank 1
in both of the forbidden regions $\rcal_F,\ \rcal_{F^*}$, as indicated in
Theorem \ref{main}(iii).

By the above, we also obtain
$$ e^{-b(z)} = \frac{4|z_2|^2(1+|z_1|^2)}{(1+|z_1|^2+|z_2|^2)^2} \qquad
\mbox{for }\ |z_1|^2 +1\leq |z_2|^2\;,$$ and similarly for
$|z_2|^2\leq |z_1|^2-1$.

\begin{rem}
On the boundary $\{|z_1|^2 = |z_2|^2-1$\}, we have $e^{-b(z)}=1$ as expected.
On
$\{|z_1| = c\}$, we have the growth rate
$e^{-b(z)} \sim 1/|z_2|^2$ as $z_2\rightarrow\infty$.\end{rem}

\subsection{The Hirzebruch surfaces}  We now consider the trapezoidal polytope of
 Figure~\ref{F2} below.

\begin{figure}[htb]
\centerline{\includegraphics*[bb= 1.9in 4.9in 4in 6.9in]{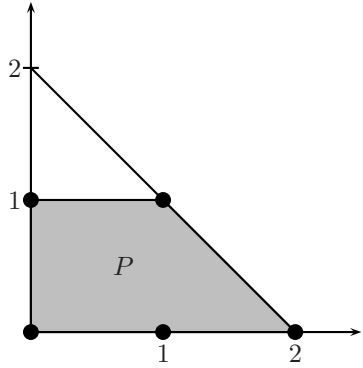}}
\caption{$P=\{(0,0),\,(1,0),\,(2,0),\,(0,1),\,(1,1)\}$}
\label{F2}\end{figure}

\noindent Again $p = 2$, but this time
$M_P\stackrel{\pi}{\rightarrow} \PP^2$ is the blow up of $(0,0,1)$, i.e.,
$M_P$
is the Hirzebruch surface $F_1$ (see \cite{F}).  Comparing with the case of the square,
we see that the classically allowed region is given by
$|z_2|^2 < |z_1|^2 +1$, and the forbidden region coincides with the upper
forbidden region
$\rcal_F=\{|z_2|^2
\ge |z_1|^2 +1\}$ from \S \ref{example1}.
Thus, the map
$z\mapsto (\tau_z,w_z)$ is the same as before when $z$ is in the
forbidden region
$\rcal_F$. Hence,
$e^{-b(z)}$ and
$\psi_P$ are also the same as in \S \ref{example1} on $\rcal_F$.  On the
classically allowable region, $e^{-b(z)}=1$ and $\psi_P=2\om_\FS$.

\subsubsection{The Hirzebruch surfaces $F_n$ ($n\ge 2$)}\label{s-ex3}
Let $n\ge 2$ and consider the polytope
$$P=\{(0,0),\,(1,0),\,\dots,(n+1,0),\,(0,1),\,(1,1)\}\;,$$
so that $M_P=F_n$ (see \cite{F}).
Here $\de=n+1$ and $\frac{1}{\de}P$ has an interior vertex $v=(\frac{1}{n+1},
\frac{1}{n+1})\in\Si^\circ$.

We see that the classically
allowed region $\mu_\Si\inv(\frac{1}{\de}P^\circ)$ is given by
$$|z_2|^2< \min\left\{\frac{|z_1|^2 +1}{n},\, \frac{1}{n-1}\right\}\;.$$  This
time, the forbidden region consists of three subregions: $\rcal_F,\ \rcal_v,\
\rcal_{F'}$, where $$\textstyle F=\{(x_1,\frac{1}{n+1}) :0\le
x_1<\frac{1}{n+1}\}
\;,\qquad F'=\{(x_1,x_2): x_2=\frac{1}{n}(1-x_1),\ \frac{1}{n+1} <x_1\le
1\}\;.$$  (See Figure~\ref{Fn} below.)

\begin{figure}[htb]
\centerline{\includegraphics*[bb= 1.5in 4.9in 4.4in 6.8in]{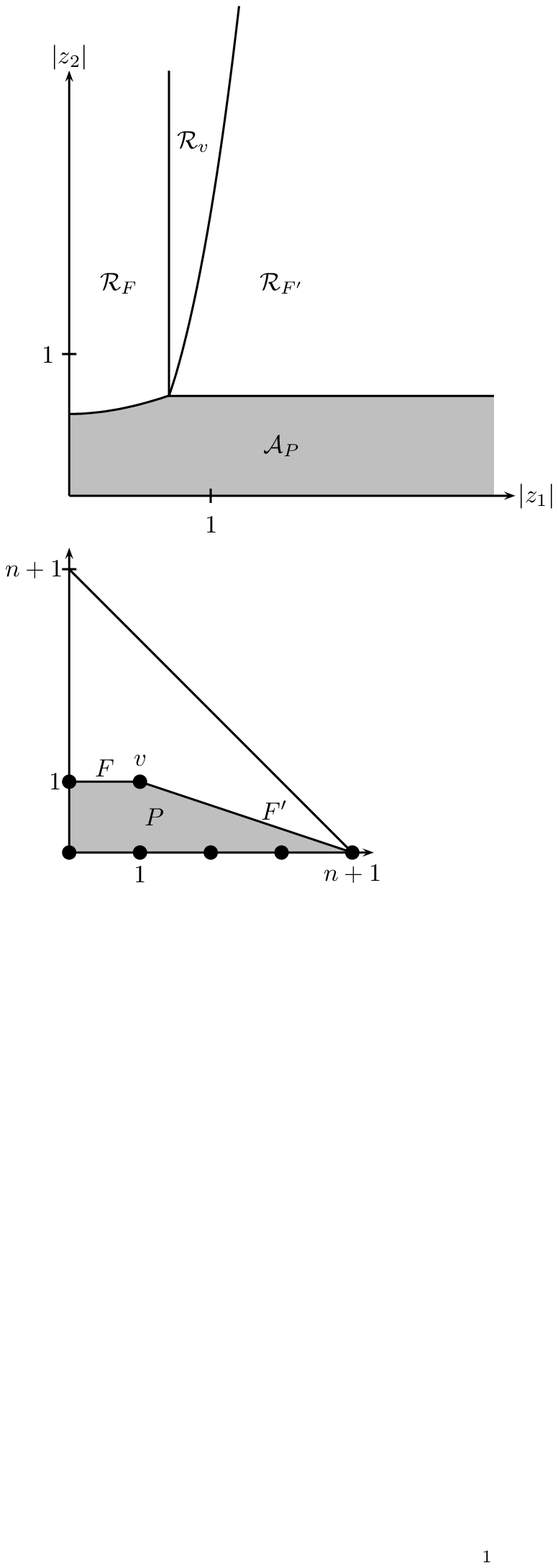}
\includegraphics*[bb= 1.6in 6.8in 5in 9.6in]{f1a.ps}}
\caption{$P=\{(0,0),\,(1,0),\,\dots,(n+1,0),\,(0,1),\,(1,1)\}$}
\label{Fn}\end{figure}

Suppose that $z$ is a point in the region $\rcal_F$.
Then $\tau_z=(0,\tau_2)$, as before.
Let $$ q(z)= (n+1) \mu_\Si (e^{\tau_z/2}\cdot z) = (a,1)\in (n+1)F \subset
\d P \qquad (0<a<1)\;.$$
Again writing
$$ |z_1|^2 =s_1,\quad |z_2|^2 = s_2\;,\quad |e^{\tau_2/2}
z_2|^2=e^{\tau_2} s_2 =\tilde s_2\;,$$
we have $$\frac{ s_1}{1+s_1+\tilde s_2} =\frac{a}{n+1}\;,\qquad \frac{\tilde
s_2}{1+s_1 +\tilde s_2}=\frac{1}{n+1}\;.$$
Therefore
$$\begin{array}{c}\di s_1 =\frac{a}{n-a},\quad \tilde s_2
=\frac{1}{n-a}= \frac{s_1}{a},\quad a= \frac{ns_1}{1+s_1} =
\frac{n|z_1|^2}{1+|z_1|^2},\\[12pt]\di e^{\tau_2} =\tilde s_2/s_2
=\frac{s_1}{as_2} =
\frac{1+|z_1|^2}{n|z_2|^2}\;.\end{array}$$
In particular,
$$ a<1 \ \Leftrightarrow\ |z_1|^2 < \frac{1}{n-1}$$
and therefore
$$\rcal_F=\left\{(z_1,z_2): |z_2|^2\ge \frac{|z_1|^2 +1}{n},\ |z_1|^2
<\frac{1}{n-1}\right\}\;.$$

We have $$
\log (1+\|e^{\tau_z/2}\cdot z\|^2) =
 \log\left (1+|z_1|^2 + \frac{1+|z_1|^2}{n|z_2|^2}
|z_2|^2\right)= \log (1 +|z_1|^2) +\log \frac{n+1}{n}\;,$$
$$\langle q(z), \tau_z\rangle =\left\langle\left(\frac{n|z_1|^2}{1+|z_1|^2}
,1\right), \left(0,\log \frac{1+|z_1|^2}{n|z_2|^2}\right)\right\rangle
= \log (1+|z_1|^2) -\log |z_2|^2-\log n\;.$$
Therefore $$u_\infty= \log|z_2|^2 +n \log
(1+|z_1|^2) +(n+1)\log (n+1) -n\log n\;.$$
Hence, $$\left.\begin{array}{rcl} \psi_P&=&\di n
\frac{\sqrt{-1}}{2\pi}
\partial\bar\partial
\log (1+|z_1|^2)\\[14pt] e^{-b(z)} &=&\di \frac{(n+1)^{n+1}}{n^n}
\frac{|z_2|^2(1+|z_1|^2)^n}{(1+|z_1|^2+|z_2|^2)^{n+1}}
\end{array}\right\} \qquad \mbox{for }\ z\in \rcal_F\;.$$

Now suppose that $z$ is a point in $\rcal_{F'}$.  Since $\tau_z \perp T_{F'}$,
we can write $\tau_z=(\tau_1, n\tau_1)$. Let
\begin{equation}\frac{1}{n+1} q(z)=  \mu_\Si (e^{\tau_z/2}\cdot z) =
\left(x_1,
\frac{1}{n}(1-x_1)\right)\in F'\;.\label{slant}\end{equation}
As before, we write
$$s_1=|z_1|^2,\quad s_2=|z_2|^2,\qquad \tilde s_1=|e^{\tau_1/2}z_1|^2=e^{\tau_1}
s_1,\quad \tilde s_2=|e^{n\tau_1/2}z_2|^2=e^{n\tau_1}s_2\;.
$$  By (\ref{slant}), we have
$$\frac{\tilde s_1}{1+\tilde s_1 +\tilde s_2} = x_1\;,\qquad
\frac{\tilde s_2}{1+\tilde s_1 +\tilde s_2} = \frac{1}{n}(1-x_1)\;.$$
Solving for $\tilde s_1,\tilde s_2$, we obtain
\begin{equation}\label{stilde}\tilde s_1 = \frac{n}{n-1}\; \frac{x_1}{1-x_1}\;,
\qquad
\tilde s_2 =
\frac{1}{n-1}\;.\end{equation}  Therefore, $e^{n\tau_1}= {\tilde s_2}/{s_2}=
{|z_2|^{-2}}/(n-1)$, so we have $$\tau_1= -\frac{1}{n}\log
(n-1)|z_2|^2\;,\qquad \tilde s_1=\frac{|z_1|^2}{(n-1)^{1/n}
|z_2|^{2/n}}\;.$$ Thus, $$
\log (1+\|e^{\tau_z/2}\cdot z\|^2) =
 \log (1+\tilde s_1 +\tilde s_2)=\log \left(\frac{n}{n-1} +
\frac{|z_1|^2}{(n-1)^{1/n}
|z_2|^{2/n}}\right)\;.$$

By (\ref{slant}),
$$\big\langle q(z), \tau_z\big\rangle= (n+1)
\left\langle\left(x_1,
\frac{1}{n}(1-x_1)\right), (\tau_1,n\tau_1)\right\rangle = (n+1)\tau_1
= -\frac{n+1}{n}\log
(n-1)|z_2|^2\;.$$
Hence by (\ref{psi-formula}),
$$\psi_P=(n+1)\frac{\sqrt{-1}}{2\pi}\ddbar\log \left(\frac{n}{n-1} +
\frac{|z_1|^2}{(n-1)^{1/n}
|z_2|^{2/n}}\right) \qquad \mbox{for }\ z\in \rcal_{F'}\;.$$
(Note that $\psi_P$ has constant rank 1 on $\rcal_{F'}$ as indicated by Theorem
\ref{main}(iii), since $\psi_P =(n+1)g^*\om_{\CP^1}$  on $\rcal_{F'}$,
where
$g$ is the multi-valued holomorphic map to $\CP^1$ given by $g(z_1,z_2)= (c_n
z_1, z_2^{1/n})$.)

By Theorem \ref{main}(iii), we know that $\psi_P=0$ on $\rcal_v$.  To complete
the description of $\psi_P$, it remains to describe the regions $\rcal_{F'}$
and $\rcal_v$. We note that a forbidden point $z$ lies in $\rcal_{F'}$ if and
only if $x_1 >\frac{1}{n+1}$. By (\ref{stilde}), this is equivalent to $\tilde
s_1 >
\frac{1}{n-1}$, or $$|z_2|^2 < (n-1)^{n-1} |z_1|^{2n}\;.$$
Therefore
$$\rcal_{F'}=\left\{(z_1,z_2): \frac{1}{n-1} \le |z_2|^2 < (n-1)^{n-1}
|z_1|^{2n},\ \ |z_1|^2 >\frac{1}{n-1}\right\}\;.$$
This leaves us with
$$\rcal_v=\left\{(z_1,z_2):  |z_2|^2 \ge  (n-1)^{n-1}
|z_1|^{2n},\ \ |z_1|^2 \ge\frac{1}{n-1}\right\}\;.$$

To summarize:
$$  \psi_P=
\left\{\begin{array}{ll}
(n+1) \frac{\sqrt{ -1}}{2\pi} \partial\bar\partial \log
(1+|z_1|^2+|z_2|^2) &
\mbox{for }\  |z_2|^2< \min\left\{\frac{|z_1|^2 +1}{n},\,
\frac{1}{n-1}\right\}\\[14pt]
n
\frac{\sqrt{-1}}{2\pi}
\partial\bar\partial
\log (1+|z_1|^2) &
\mbox{for }\  |z_2|^2\ge \frac{|z_1|^2 +1}{n},\ |z_1|^2
<\frac{1}{n-1}\\[14pt](n+1)
\frac{\sqrt{-1}}{2\pi}\ddbar\log \left(\frac{n}{n-1} +
\frac{|z_1|^2}{(n-1)^{1/n}
|z_2|^{2/n}}\right) &
\mbox{for }\  \frac{1}{n-1} \le |z_2|^2 < (n-1)^{n-1}
|z_1|^{2n}\\[14pt]
 0 &
\mbox{for }\ |z_2|^2 \ge  (n-1)^{n-1}
|z_1|^{2n},\ \ |z_1|^2 \ge\frac{1}{n-1}\end{array}\right.$$

\subsection{Remarks on the behavior at transition points}\label{airy-tale} So far, we
have studied the uniform asymptotics of  $\Pi_{|NP}(z, z)$ on compact subsets of the
classically allowed region and on compact sets without transition points in the
forbidden region.
We now address a few words to the  asymptotics at the
boundary of  the classically allowed region, which is its caustic set.  For background
on oscillatory integrals and  caustics,  we refer to
\cite{AGV}.

We thus  consider $\Pi_{|NP}(z, z)$ as an oscillatory integral (denoted $K_N(z) $ in
(\ref{szego3})) on $(\C^*)^m$. Since it is ${\bf T}^m$ invariant we may  regard it more
simply as an oscillatory integral in $z = e^{\rho}$ on $\R_+^m$ or on $p \Sigma,$
with the phase $\Psi$ of (\ref{PHASE2})
regarded as a function on $ {\bf T}^m  \times M_P \times \R_+^m$.
The asymptotics of $\Pi_{|NP}(z, z)$ are determined by the component
\begin{equation} C^0 = \{(0, w; e^{\rho} ) \in  \R_+^m \times M_P
 :  \mu_{p \Sigma} (e^{\rho})  =  \mu_P(w) \},
\end{equation}
 of the
critical set on which $\phi = 0$, so we may (and will) ignore
other critical points which may occur.  The map $(0, w; e^{\rho}) \to w$ gives
a  natural identification of
$C^0$ with  $M_P$.

As is standard in the theory of oscillatory integrals
depending on parameters, we define
the map
\begin{equation} \iota_{\Psi} : C^0 \to T^* \R^m_+,\;\;\;\;
\iota_{\Psi}(0, w; e^{\rho}) = (e^{\rho} , d_{\rho} \Psi),
\end{equation}
and we denote the image by $\Lambda^0 =  \iota_{\Psi}
(C^0)$. Over the set of parameters $z$ for which
the critical point set $C^0_z$ is a non-degenerate critical manifold,  $\iota_{\Psi}$ is
a fibration and  $\Lambda^0$ is a
Lagrangean manifold.  The set of $z = e^{\rho}$ for which the phase is degenerate is known
as the caustic set.
Under the identification of $C^0 \sim M_P$, we have
$\frac{1}{2} d_z \Psi ( w; e^{\rho}) \equiv 0\; \mbox{on}\;
C^0$, hence  $ \Lambda^0 $ is the zero section of $T^*P$.
Identifying $\iota_{\Psi}$  with the map $\mu_P \times 0,$
we find that  the caustic set is the boundary of the classically allowed region.

Uniform asymptotic expansions of $\Pi_{| NP}(z,z)$ in the vicinity of the caustic
could be obtained with a careful study of the singularities of the moment map $\mu_P$.
The uniform expansion would give a smooth transition in a boundary layer from the $N^m$ mass density in
the interior of $P$ (the illuminated or allowed region)  to the $e^{- N b(z)}$ mass density
in the exterior (the shadow or forbidden region).  The order of the asymptotics
along the caustic depends on two competing effects: on the one hand, oscillatory
integrals whose (analytic) phase has a  degenerate critical point decay more slowly
than in the non-degenerate case (by the index of the singularity at the critical point,
see \cite{AGV}). On the other hand, the principal term of the  amplitude may vanish
on the caustic, causing more rapid decay.  We will be content to illustrate
this in the simplest example.

\subsubsection{Example: $P = [0, 1], p = 2$}

We put  $m = 1$,  $P = [0, 1], p = 2,
p \Sigma = [0, 2].$ Thus, $M_P = \CP^1$. The only interior facet is $F =
\{1\}$; recalling (\ref{PHASE2}), our phase is:
\begin{equation} \label{PHASE4} \Psi(w,\phi; z) =
\log \frac{ 1+ e^{-i \phi} |w|^2 }{1 +  |w|^2}  + 2\log \frac{1
+  e^{i \phi} |z|^2} {1 + |z|^2} \;.
\end{equation}
We are interested in the asymptotics of $\Pi_{| N P}(z,z)$ when
$$\mu_{2 \Sigma}(z) =  \frac{  2 |z|^2}{1 +  |z|^2} = 1 \iff  |z| = 1.$$
We claim that $\Pi_{|N P}(z,z) \sim C N^{1} $ when $|z| = 1.$

 Since $\mu_{P}(w) = \frac{|w|^2}{1 + 2 |w|^2},$ we
have
$\mu_P(w) = \{1\} \iff w = \infty.$
 To determine the asymptotics of $\Pi_{| N P}(z,z)$ at
$|z| = 1$ , we put  $\eta = \frac{1}{w} \in \C$ and consider the Taylor
expansion of the phase at $\eta = 0, \phi = 0$. We change to polar coordinates
$ \eta = r e^{i \theta}$ around the fixed point of the toric $S^1$ action.
The phase and amplitude are independent of $\theta$ and the principal part of
the Taylor series (at $|z|=1$) has the form:
$$\Psi(w, \phi;1) \sim i \phi r^2 -\frac 14 \phi^2 + O(|\phi|^3+|\phi|^2
|r|^2 + |\phi||r|^4)\;. $$
By a general result \cite{AGV}, the asymptotics of (generic) oscillatory integrals are determined by
 the principal part of the Taylor expansion of the phase at  the critical point,
which gives the
 model integral
$$I_N=   \int_{-\pi}^\pi \int_0^{\infty}  e^{ N (i \phi r^2 -\frac 14 \phi^2)}
A(r,
\phi) r dr d\phi,$$ where $A$ is an amplitude decaying sufficiently
rapidly at infinity. Here, we also used that the amplitude of $\Pi_{| N P}(z,z)$
is simply the volume density in our local coordinates $(\phi, r, \theta).$
We observe that the  phase $f(r, \phi) = i \phi r^2 - \frac 14
 \phi^2$ is {\it weighted homogeneous} in the sense that
$$f(N^{1/4} r, N^{1/2} \phi) = N f(r, \phi). $$
We also observe that the amplitude vanishes to order one as $r \to 0$.
We change variables as indicated to obtain that
$$\Pi_{|NP}(1,1)\sim N^2 I_N = N^2 \int_{\R} \int_{\R}N\inv e^{  (i \phi r^2 + C(z)
\phi^2)}
\rho(N^{-1/4} r, N^{-1/2} \phi) r dr d\phi \sim C N \;,  $$ as claimed.

Alternatively, we observe that  the index of oscillation   of the phase with $|z| = 1$
equals $-3/4$. This agrees with the general result (under
assumptions satisfied here)   that the index of oscillation equals
the `remoteness' of the Newton polygon of the phase (\cite{AGV},
\S II.6, Theorem~4), i.e.\ the reciprocal of the value of $t$
such that the ray $(t,t), t> 0$ intersects the boundary of the
Newton polygon. The index of the singularity is therefore $\delta
= 1/4$, and therefore $I_N\sim CN\inv$.

\end{document}